\documentclass[11pt, reqno]{amsart}

\usepackage{amsmath}
\usepackage{fullpage}
\usepackage{amssymb}
\usepackage{amsthm}
\usepackage{bbm}
\usepackage{subfigure}
\usepackage{float}
\usepackage{graphicx}
\usepackage{color}
\usepackage{epstopdf}
\usepackage{enumerate}
\usepackage[title]{appendix}

\newtheorem{lemma}{Lemma}[section]
\newtheorem{theorem}[lemma]{Theorem}
\newtheorem{proposition}[lemma]{Proposition}

\newtheorem{definition}[lemma]{Definition}
\allowdisplaybreaks

\newcommand{\fe}{{\rm e}}
\DeclareMathOperator{\re}{Re}
\numberwithin{equation}{section}
\newcommand{\msc}[2][2000]{
\let\@oldtitle\@title
\gdef\@title{\@oldtitle\footnotetext{#1 \emph{Mathematics subject classification.} #2}}}

\title[Explicit splitting for the original ZS]{Low-regularity error estimates of a filtered Lie-Trotter splitting scheme for the Zakharov system in arbitrary dimensions}

\author{Lun Ji}
\address{SKLMS, Academy of Mathematics and Systems Science, Chinese Academy of Sciences, Beijing 100190, China (L.~Ji)}
\address{School of Mathematical Sciences, University of Chinese Academy of Sciences, Beijing 100049, China (L.~Ji)}
\address{Current address: Department of Applied Mathematics, The Hong Kong Polytechnic University, Hung Hom, Hong Kong SAR, China (L.~Ji)}
\email{lun-422.ji@polyu.edu.hk}

\author{Hang Li}
\address{LJLL (UMR 7598), Sorbonne Universit\'e, UPMC, 4 place Jussieu, 75005, Paris, France (H.~Li)}
\email{li727263992@gmail.com}

\author{Chunmei Su}
\address{Yau Mathematical Sciences Center, Tsinghua University, 100084 Beijing, China (C.~Su)}
\email{sucm@tsinghua.edu.cn}

\thanks{L. Ji was partially supported by the National Natural Science Foundation of China (No.~92470119) and Research Grants Council (GRF project No.~15306123). H. Li is supported by the European Research Council (ERC) under the European Union's Horizon 2020 research and innovation program (grant agreement No.~850941), as well as by a postdoctoral fellowship from the Fondation Sciences Math\'ematiques de Paris (FSMP). C. Su is supported by the National Natural Science Foundation of China (No.~12522118).}

\begin{document}

\begin{abstract}
We establish low-regularity error estimates for a fully discrete, filtered Lie--Trotter splitting scheme for the Zakharov system on the $d$-dimensional torus. Our analysis relies on multilinear estimates in discrete Bourgain spaces. Specifically, for an exact solution $(E,z,z_t)$ with regularity $H^{s+r+1/2}\times H^{s+r}\times H^{s+r-1}$ on a fixed time interval, we prove an error bound of order $\mathcal{O}(\tau^{s/2}+N^{-s})$ in $H^{r+1/2}\times H^r\times H^{r-1}$ for $s\in(0,2]$, where $r=\max\{0,\tfrac d2-1\}$, $\tau$ is the time step, and $N$ is the number of spatial grid points in each coordinate direction. The analysis treats the coupled nonlinearities
directly, without introducing the time-differentiated auxiliary variables of an Ozawa-type transformation. Under the stated CFL condition, which renders the filter inactive on the Fourier collocation space, the scheme preserves discrete mass exactly in exact arithmetic. Numerical experiments illustrate the theoretical convergence rates and show small relative energy errors over the tested time intervals.
\end{abstract}


\maketitle

\section{Introduction}

We consider the Zakharov system (ZS) on a $d$-dimensional torus $\mathbb{T}^d$:
\begin{equation}\label{zs}
\left\{\begin{aligned}
&iE_t(t,x)=-\Delta E(t,x)+z(t,x)E(t,x),\\
&z_{tt}(t,x)=\Delta z(t,x)+\Delta(|E(t,x)|^2),\quad\quad\quad \ \ t>0,\\
&E(0,x)=E_0(x),\quad z(0,x)=z_0(x),\quad z_t(0,x)=z_1(x),
\end{aligned}\right.
\end{equation}
where $E(t,x): \mathbb{R}^+\times\mathbb{T}^d\rightarrow \mathbb{C}$ represents the highly oscillatory electric field envelope, and $z(t,x): \mathbb{R}^+\times \mathbb{T}^d\rightarrow \mathbb{R}$ denotes the ion density deviation from equilibrium. Introduced by Zakharov \cite{Zakspj} to model Langmuir wave propagation in plasmas, the ZS has emerged as a prototypical model for dispersive and nonlinear interactions with applications in nonlinear optics, fluid dynamics and plasma physics \cite{Davps,Degspj,Zakspj}. We assume that the initial data satisfy
$$
(E_0, z_0, z_1) \in H^s(\mathbb{T}^d)\times H^l(\mathbb{T}^d)\times H^{\,l-1}(\mathbb{T}^d).
$$
The well-posedness of the ZS with such data has been extensively studied in Bourgain spaces: in particular, it was shown in \cite{Kishijam} and \cite{Takdie} that the system is locally well-posed under the condition $0 \le s-l \le 1$ together with the following dimension-dependent constraints:
\begin{align*}
&1/2\leq l+1/2\leq 2s, &{\rm for} \quad d=1,\\
&1 \leq l+1 \leq 2s, &{\rm for} \quad d=2,\\
&d-1< l+d/2\leq 2s, &{\rm for} \quad d\geq 3.
\end{align*}
Moreover, global well-posedness for small initial data was established in one dimension \cite{Bourgduke}.
Integrating the wave equation in \eqref{zs} over $\mathbb{T}^d$ yields
$$
\widehat{(z_{tt})}_0=\int_{\mathbb{T}^d}z_{tt}(t,x)dx=0,
$$
thus, without loss of generality, we may assume
\begin{equation}\label{zeromean}
\int_{\mathbb{T}^d}z_t(t,x)dx\equiv\int_{\mathbb{T}^d}z_1(x)dx=0, \quad \forall \ t\ge 0.
\end{equation}

The ZS has a Hamiltonian structure and admits several conserved quantities. Specifically, mass is preserved:
\begin{equation}\label{massc}
\mathcal{M}(t)= \int_{\mathbb{T}^d} |E(t,x)|^2 dx \equiv \mathcal{M}(0),
\end{equation}
and, under the assumption \eqref{zeromean}, the system also conserves energy:
\begin{equation}\label{energyc}
\mathcal{H}(t) = \int_{\mathbb{T}^d} \big(|\nabla E(t,x)|^2 + z(t,x)|E(t,x)|^2 + \frac12||\nabla|^{-1}z_t(t,x)|^2 + \frac12|z(t,x)|^2\big) dx \equiv \mathcal{H}(0),
\end{equation}
where $|\nabla|=\sqrt{-\Delta}$.

The literature includes extensive numerical studies on the ZS. A variety of efficient and accurate numerical methods have been proposed and analyzed, including the finite difference method \cite{Bao3,Caimcom,Changjcp} the time-splitting spectral method \cite{Bao1,Bao2,Gaucknm} the exponential integrator spectral method \cite{Herrima,Limcom} and the discontinuous Galerkin method \cite{Xia} among others.
The principal difficulty in designing and analyzing numerical schemes for the ZS is the mismatch of regularity requirements between its Schr\"odinger and wave components. Writing the Duhamel formulas:
\begin{align*}
E(t)&=\fe^{it\Delta}E_0-i\int_0^t\fe^{i(t-\xi)\Delta}\big[z(\xi)E(\xi)\big]\, d\xi, \\
z(t)&=\cos(t|\nabla|)\, z_0+\frac{\sin(t|\nabla|)}{|\nabla|}\, z_1-\int_0^t|\nabla|\sin((t- \xi)|\nabla|)|E(\xi)|^2\, d\xi,
\end{align*}
makes this tension explicit: the Schr\"odinger equation suggests that $E$ and $z$ should have comparable regularity, whereas the wave equation imposes an extra derivative on $|E|^2$. Concretely, a direct calculation gives:
\begin{align*}
\|E(t)\|_{H^s} &\lesssim \|E_0\|_{H^s} + \int_0^t \|z(\xi)\|_{H^s} \|E(\xi)\|_{H^s} \, d\xi, \quad s > \tfrac{d}{2},\\
\|z(t)\|_{H^l} &\lesssim \|z_0\|_{H^l} + \|z_1\|_{H^{l-1}} + \int_0^t \|E(\xi)\|_{H^{l+1}}^2 \, d\xi, \quad l+1 > \tfrac{d}{2},
\end{align*}
which exhibits a genuine \emph{loss of derivatives}. This loss complicates both the construction and the analysis of numerical methods: standard discretizations and their error analyses can inherit or even amplify regularity degradation through the nonlinear coupling, while implicit or semi-implicit schemes---although effective in controlling derivative loss---are computationally costly. Hence there is significant interest in explicit schemes that avoid or mitigate the loss of derivatives.

Auxiliary-variable reformulations have been used to address the loss of derivatives in the ZS. Following the approach of Ozawa and Tsutsumi \cite{Ozawaprims}, Herr and Schratz \cite{Herrima} used the variables
$$
F = \partial_t E,\quad \mathcal{L}(t,x) = E_0(x) + \int_0^t F(s,x)\, ds,
$$
to express the system as
$$
\left\{\begin{aligned}
&i\partial_t F(t, x)+\Delta_{}  F(t, x)=z(t, x)F(t, x)+\partial_t z(t, x)\mathcal{L}(t, x),\\
&(\partial_{tt}-\Delta_{}) z(t, x)=\Delta_{} |E(t, x)|^2,\\
&E(t, x)=(-\Delta_{}+1)^{-1}\left[iF(t, x)-(z(t, x)-1)\mathcal{L}(t, x)\right],
\end{aligned}\right.
$$
with initial data
$$
F(0,x) = i(\Delta E_0 - z_0 E_0), \quad z(0,x) = z_0(x), \quad \partial_t z(0,x) = z_1(x), \quad E(0,x) = E_0(x).
$$
Based on this
reformulation, they constructed explicit trigonometric time
integrators. For $s>d/2$, their first-order convergence theorem
estimates the error in
$H^{s+2}\times H^{s+1}\times H^s$, assuming that the exact solution
$(E,z,\partial_t z)$ has the corresponding regularity in
$H^{s+4}\times H^{s+3}\times H^{s+2}$ on the time interval under
consideration. They also presented a second-order extension and
stated its convergence in the same error norm under the stronger
regularity assumption
$H^{s+6}\times H^{s+5}\times H^{s+4}$ \cite{Herrima}. More recently, Li and Su \cite{Limcom} developed first- and
second-order low-regularity exponential-type integrators whose
construction also employs this auxiliary-variable reformulation.
Their convergence results require one and two additional spatial
derivatives, respectively, relative to the Sobolev norms in which
the errors are estimated, thereby reducing the regularity
requirements compared with the classical trigonometric integrators.

Concerning on the splitting methods, Gauckler \cite{Gaucknm} analyzed an existing fully discrete Lie--Trotter splitting method combining Fourier collocation in space with splitting in time. The method is formulated in the original variables of the ZS, whereas a discrete analogue of the Ozawa--Tsutsumi transformation is introduced in the error analysis to control the loss of derivatives. For $s>d/2$ and $\sigma\geq 2$, he established an error bound of order $\mathcal{O}(\tau+K^{-\sigma})$ in $H^{s+2}\times H^{s+1}\times H^s$, assuming that the exact solution remains uniformly bounded in $H^{s+2+\sigma}\times H^{s+1+\sigma}\times H^{s+\sigma}$. This result holds for sufficiently small $\tau$ and sufficiently large $K$ under the CFL-type restriction $d\tau K^2\leq c<2\pi$, where $c$ is a fixed positive constant, $\tau$ is the time step, and $K$ is the Fourier discretization parameter corresponding to $2K$ collocation points in each
spatial direction.

Despite their success under restrictive smoothness assumptions, the aforementioned methods or the corresponding error analysis rely on auxiliary-variable reformulations. Recently, Li and Wang \cite{LiWang} developed a framework based on a different reformulation of the ZS, yielding a time-symmetric exponential wave integrator with two discrete conservation laws. They established second-order convergence under a time-step restriction independent of the spatial discretization.

Obtaining rigorous error estimates for rough solutions remains an important objective. Although low-regularity exponential-type integrators can reduce the number of additional spatial derivatives required for convergence, their regularity requirements must also be assessed relative to the norms in which the errors are measured. For example, the convergence results of Li and Su \cite{Limcom} discussed above still assume a base Sobolev index strictly greater than $d/2$. A standard tool in Sobolev-based error analyses is the multiplication estimate (see, e.g., \cite{CaoJ,LWnumermath,WuandYao2022,YABYBGK,MaLiNS,YAB,Limcom,Baiima,KWZ})
\begin{equation}\label{cbe}
\|fg\|_{H^r}\leq C_{d,s,r}\|f\|_{H^r}\|g\|_{H^s},\qquad 0\le r\le s,\quad s>d/2,
\end{equation}
where $H^r=H^r(\mathbb{T}^d)$. Without additional assumptions, this estimate does not extend in the same form to $s\leq d/2$. To address rougher solutions, it is therefore natural to seek space--time estimates that exploit the dispersive structure of the equations, rather than relying solely on spatial multiplication estimates.

Discrete Strichartz estimates \cite{Lifocm,Ostfocm} provide one such approach to the error analysis of numerical methods for dispersive PDEs. For periodic problems, a complementary approach
is provided by the discrete Bourgain-space framework \cite{Ostjems}, which adapts the continuous framework introduced by Bourgain \cite{Bourgduke} to discrete time. These spaces combine spatial regularity with temporal regularity, described by an additional index $b$ and measured relative to the linear dispersive flow. This allows the analysis to exploit the space--time structure of nonlinear interactions. The discrete framework has been extended to higher dimensions \cite{Jisiam,Jiima,Jiicos} and to other dispersive models \cite{CLWY,Jimcom,JiZhao,Roupaa}. Building on this framework, we establish low-regularity error estimates for a fully discrete filtered splitting method for the ZS, formulated in the original variables. The analysis uses multilinear estimates in discrete Bourgain spaces to control the nonlinear coupling between the Schr\"odinger and wave components.

\subsection*{Contributions and novelty.}

We develop a discrete Bourgain-space analysis of a fully discrete, filtered Lie--Trotter splitting method for the Zakharov system. The scheme uses the equivalent first-order wave variable \eqref{ztou}, while the error analysis avoids the Ozawa--Tsutsumi reformulation involving $(F,\mathcal L)$. The main results and technical contributions are as follows:
\begin{itemize}
\item[(i)] We establish discrete multilinear estimates for the nonlinear Schr\"odinger--wave coupling (Theorem~\ref{theomult}). These estimates account for the distinct dispersion relations of the two components, the discrete temporal frequency structure, and the imposed frequency cutoff. They provide the nonlinear bounds needed for the low-regularity convergence analysis.

\item[(ii)] Leveraging the multilinear estimates, we establish convergence of the Lie-Trotter splitting scheme for the original Zakharov system under some regularity assumptions that align with the known continuous well-posedness theory. This result breaks the classical barrier $s > d/2$ and successfully handles regimes of extremely low regularity ($s \leq d/2$).

\item[(iii)] We bridge the gap between critically low regularity and full discretisation at the $H^s$ level for methods with exponential terms and non-integer Sobolev exponents $s$. While complete error analyses have previously been carried out for methods with finite spectral support---for the KdV equation \cite{CLWY} and for the nonlinear Schr\"odinger equation at the $L^2$ level \cite{Jisiam}---our work extends these results to a substantially more general setting.

\item[(iv)] Under the stated CFL condition, for which the filter is inactive on the collocation space, the scheme preserves discrete mass in exact arithmetic. Numerical experiments additionally show small relative energy errors over the tested time interval, together with a visible long-time drift. These observations motivate a separate analysis of the long-time energy behavior; such an analysis is beyond the scope of the present paper.
\end{itemize}

\subsection*{Outline of the paper.}
The remainder of this paper is organized as follows. Section~\ref{sectionmethod} introduces the filtered Lie splitting method and states the main convergence result (Theorem~\ref{mainthm}). In Section~\ref{sectionbourgain} we present the Bourgain space framework for the ZS. The local and global error analyses are given in Sections~\ref{sectionlocal} and \ref{sectionglobal}, respectively. Finally, numerical results that confirm the analysis are shown in Section~\ref{sectionnumerexp}, and the conclusion is given in Section~\ref{conclusion}.

\subsection*{Notations.}
\begin{itemize}
\item We denote the time step size by $\tau$.

\item For $y\in \mathbb{R}^d$, the Japanese bracket is $\langle y\rangle=(1+|y|^2)^{\frac12}$.

\item For $a, b\ge 0$, we use the notation $a \lesssim b$ whenever there exists a generic constant $C$ independent of $\tau\in(0,1)$ such that $a\le Cb$. The notation $\lesssim_\gamma$ emphasizes that the constant $C$ depends in particular on $\gamma$. Moreover, $a\sim b$ means that $a\lesssim b\lesssim a$.

\item For a function $f(x)=\sum\limits_{k\in\mathbb{Z}^d}c_k\fe^{i\langle k,x\rangle}$, we define the fractional derivative by the corresponding Fourier multiplier:
$$
|\nabla|^\alpha f=\sum\limits_{k\neq 0}|k|^\alpha c_k\fe^{i\langle k, x\rangle},\quad \alpha\in\mathbb{R},
$$
where $\langle\cdot,\cdot\rangle$ denotes the Euclidean inner product on $\mathbb R^d$. Note that if $\alpha<0$, we require that the zero Fourier mode vanish, i.e., $\widehat{f}_0=0$.

\item Parameters are fixed as
\begin{equation}\label{para}
\begin{aligned}
s_0=\max\{0,\tfrac d2-1\},&\quad s_1>0,\quad s_2=s_0+s_1,\\
b_0\in\big(\tfrac12,\min\{\tfrac58,\tfrac12+\tfrac12s_1\}\big), &\quad b_1\in\big(\!\max\{\tfrac38,\tfrac12-\tfrac12s_1\},1-b_0\big).
\end{aligned}
\end{equation}
\end{itemize}

\section{The filtered Lie splitting method for ZS}\label{sectionmethod}

Owing to the zero-mean assumption \eqref{zeromean}, we set
\begin{equation}\label{ztou}
u=z-i|\nabla|^{-1}z_t.
\end{equation}
System \eqref{zs} can then be rewritten as an equivalent first-order system (cf. \cite{Limcom}):
\begin{equation}\label{zseqv}
\left\{
\begin{aligned}
&iE_t=-\Delta E+\frac12(u+\overline{u})E,\\
&iu_t=-|\nabla| u-|\nabla|(E\overline{E}),\\
&E(0,x)=E_0(x),\quad u(0,x)=u_0(x)=z_0(x)-i|\nabla|^{-1}z_1(x).
\end{aligned}
\right.
\end{equation}
The system \eqref{zseqv} admits a natural splitting into two subsystems:
\begin{equation}\label{zssub1}
\left\{\begin{aligned}
&F_t=i\Delta F,\\
&v_t=i|\nabla|v,
\end{aligned}\right.
\end{equation}
and
\begin{equation}\label{zssub2}
\left\{\begin{aligned}
&G_t=-\frac i2(w+\overline{w})G,\\
&w_t=i|\nabla|(G\overline{G}),
\end{aligned}\right.
\end{equation}
which corresponds to the linear and nonlinear parts, respectively. Note that this splitting is different from that given in \cite{Gaucknm}. The linear subsystem~\eqref{zssub1} can be solved exactly: $F(t)=\fe^{it\Delta}F_0$, $v(t)=\fe^{it|\nabla|}v_0$. For the nonlinear subsystem, noticing that $\mathrm{Re}(w_t)=0$, i.e., $w+\overline{w}\equiv w_0+\overline{w_0}$, we obtain $G(t)=\fe^{-it\mathrm{Re}(w_0)}G_0$; in particular, $|G|$ is time-independent and $w(t)=w_0+it|\nabla|(G_0\overline{G_0})$.
Combining these exact subflows via the splitting technique \cite{Mcacta,Lubich08} yields the Lie splitting scheme:
\begin{equation}\label{lie}
\left\{\begin{aligned}
&E_{n+1}=\fe^{i\tau\Delta}\fe^{-\frac i2\tau(u_n+\overline{u_n})}E_n,\\
&u_{n+1}=\fe^{i\tau|\nabla|}(u_n+i\tau|\nabla|(E_n\overline{E_n})).
\end{aligned}\right.
\end{equation}
When combined with a Fourier collocation in space to form a fully discrete scheme, this method is subject to the CFL restriction  (see also Section~\ref{sectionnumerexp}):
\begin{equation}\label{CFL}
dN^2\tau\leq c<2\pi.
\end{equation}
Here we introduce a filter operator, commonly used in discrete Bourgain settings \cite{Jisiam,Jiima,Ostjems,Roufocm}:
\begin{equation}\label{proj}
\Pi_\tau =\overline{\Pi}_\tau=\chi\left(\dfrac{-i\nabla}{\tau^{-1/2}}\right),
\end{equation}
where $\chi$ denotes the characteristic function of the cube $[-\frac12c^{1/2}d^{-1/2},\frac12c^{1/2}d^{-1/2}]^d$ with $0<c<2\pi$ and therefore enforces a frequency cutoff. For the enlarged cutoff $\Pi_{\theta/4}$ used in the proof of Theorem~\ref{theomult}, the condition $c<2\pi$ provides a uniform separation from the discrete modulation period $2\pi/\tau$, excluding wrap-around resonances.

In the spirit of \cite{Jisiam} we also analyze the full discretization error. For this purpose we introduce the standard Fourier interpolation operator $T_N$ associated with the discrete Fourier transform (DFT):

\begin{definition}\label{tndefine}
For every even $N>0$, let $T_N$ act on continuous functions $u: \mathbb{T}^d \to \mathbb C$ by
\begin{equation}\label{tndef}
T_N(u)(x)=\frac1{N^d}\sum\limits_{k_1=-N/2}^{N/2-1}\cdots\sum\limits_{k_d=-N/2}^{N/2-1}\mathcal{F}_N(u)(k)\fe^{i\langle k,x\rangle},
\end{equation}
where $\mathcal{F}_N(u)$ denotes the $d$-dimensional discrete Fourier transform:
\begin{equation}\label{dft}
\mathcal{F}_N(u)(k)=\sum\limits_{l_1=-N/2}^{N/2-1}\cdots\sum\limits_{l_d=-N/2}^{N/2-1}u(\tfrac{2\pi l}{N})\fe^{-\tfrac{2i\pi\langle l,k\rangle}{N}}, \quad 0 \leq k_1,\ldots,k_d\leq N-1.
\end{equation}
\end{definition}

Denote
\begin{equation}\label{theta}
\theta=\max\{\tau,cd^{-1}N^{-2}\},
\end{equation}
where $c\in(0, 2\pi)$ is a fixed number for the cutoff estimates below. Now we define the fully discrete filtered Lie splitting method
\begin{equation}\label{filtlie}
\begin{aligned}
&E_{n+1}=\Phi_N^\tau(E_n,u_n)=\fe^{i\tau\Delta}\Pi_\theta T_N(\fe^{-\frac i2\tau(\Pi_\theta u_n+\Pi_\theta\overline{u_n})}\Pi_\theta E_n),\\
&u_{n+1}=\Psi_N^\tau(E_n,u_n)=\fe^{i\tau|\nabla|}(i\tau|\nabla|\Pi_\theta T_N(\Pi_\theta E_n\Pi_\theta\overline{E_n})+\Pi_\theta u_n),
\end{aligned}
\end{equation}
with $E_0=\Pi_\theta E_0$, and $u_0=\Pi_\theta u_0$.

To establish the error estimate, similar to \cite{Jisiam}, we introduce the filtered ZS:
\begin{equation}\label{filtzs}
\left\{\begin{aligned}
&iE^\theta_t=-\Delta E^\theta+\frac12\Pi_\theta((\Pi_\theta u^\theta+\Pi_\theta\overline{u^\theta})\Pi_\theta E^\theta),\\
&iu^\theta_t=-|\nabla| u^\theta-|\nabla|\Pi_\theta(\Pi_\theta E^\theta\Pi_\theta\overline{E^\theta}),\\
&E^\theta(0,x)=\Pi_\theta E_0(x),\quad u^\theta(0,x)=\Pi_\theta u_0(x)=\Pi_\theta(z_0(x)-i|\nabla|^{-1}z_1(x)).
\end{aligned}\right.
\end{equation}
Note that \eqref{filtlie} is obtained by incorporating the Fourier interpolation operator $T_N$ for spatial discretization and the filters $\Pi_\theta$ into the original Lie splitting scheme \eqref{lie}, rather than by applying exact Lie splitting to the filtered system \eqref{filtzs}.

Under the CFL condition \eqref{CFL}, the filter is inactive in the collocation space. This means that the nonlinear phase step preserves the modulus at each grid point, while the linear step preserves the discrete $L^2$ norm. Consequently, the scheme preserves the discrete mass exactly in exact arithmetic. In floating-point computations, small deviations may arise from accumulated rounding errors in the FFTs,
phase evaluations, and other arithmetic operations; see \cite{Einkcma} for a discussion of round-off propagation in FFT-based methods.

Consequently, the numerical solutions $(z_n,\dot{z}_n)$ for approximating $(z(\cdot, t_n)$, $z_t(\cdot, t_n))$ can be recovered from $(E_n, u_n)$ via \eqref{ztou}:
\begin{equation}\label{zurelation}
z_n=\tfrac12(u_n+\overline{u_n}), \quad \dot{z}_n=\tfrac{i|\nabla|}{2}(u_n-\overline{u_n}).
\end{equation}

We conclude this section with the main result of the paper.
\begin{theorem}\label{mainthm}
Let $s_0=\max\{0,d/2-1\}$ be as defined in \eqref{para}, and consider initial data $(E_0,z_0,z_1)\in H^{s_2+\frac12}(\mathbb{T}^d)\times H^{s_2}(\mathbb{T}^d)\times H^{s_2-1}(\mathbb{T}^d)$ with $s_2>s_0$. Denote by $(E,z,z_t)$ the exact solution of \eqref{zs} on $[0,T]$, and by $(E_n, z_n, \dot{z}_n)$ the numerical solution obtained from the scheme \eqref{filtlie} together with \eqref{zurelation}. Then there exist $N_0\in\mathbb{N}$, $\tau_0>0$ and a constant $C_T>0$ such that for all $\tau\in(0,\tau_0]$ and $N>N_0$, the following estimate holds:
$$
\|E_n-E(t_n)\|_{H^{s_0+\frac12}(\mathbb{T}^d)}+ \|z_n-z(t_n)\|_{H^{s_0}(\mathbb{T}^d)}+\|\dot{z}_n-z_t(t_n)\|_{H^{s_0-1}(\mathbb{T}^d)}\leq C_T\theta^{s/2},
$$
where $\theta$ is given by \eqref{theta}, $s=\min\{s_2-s_0,2\}$, and $0\leq n\tau\leq T$. The constants $\tau_0$, $N_0$ and $C_T$ depend on $T$ but are independent of $n$ and $\tau$.
\end{theorem}

\section{A Bourgain framework}\label{sectionbourgain}

In this section, we introduce the Bourgain space framework which is essential in our low-regularity analysis; further background can be found in \cite{Jiima,Jimcom} and \cite{Ostjems}.

We define two Bourgain spaces associated with the Schr\"odinger and wave components of the ZS:
\begin{align} \label{bourg1}
\|v\|_{X_1^{s,b}}&=\|\fe^{-it\Delta}v\|_{H_t^bH_x^s}=\|\langle k\rangle^s\langle\sigma+|k|^2\rangle^b\widetilde{v}(\sigma,k)\|_{L_\sigma^2 l_k^2},\\
\label{bourg2}\|v\|_{X_2^{s,b}}&=\|\fe^{-it|\nabla|}v\|_{H_t^bH_x^s}=\|\langle k\rangle^s\langle\sigma-|k|\rangle^b\widetilde{v}(\sigma,k)\|_{L_\sigma^2 l_k^2},
\end{align}
where $\widetilde{v}(\sigma,k)$ denotes the space-time Fourier transform of $v$:
$$
\widetilde{v}(\sigma,k)=\int_{\mathbb{R}\times\mathbb{T}^d}\fe^{-i\sigma t-i\langle k,x\rangle}v(t,x)\,dxdt, \quad\sigma\in\mathbb{R},\quad k\in\mathbb{Z}^d.
$$
Next we collect several standard properties of the projector $\Pi_\theta$ and of the Bourgain spaces. Note that, if we write $X_{2,-}^{s,b}=X_2^{s,b}$ and define $X_{2,+}^{s,b}$ by replacing $\langle\sigma-|k|\rangle^b$ in \eqref{bourg2} with $\langle\sigma+|k|\rangle^b$, then $\|\overline v\|_{X_{2,+}^{s,b}}=\|v\|_{X_{2,-}^{s,b}}$. Noticing that the linear parts of the equations for $u$ and $\overline{u}$ in \eqref{zseqv} involve opposite half-wave propagators. Accordingly, we measure $u$ in $X_2^{s,b}=X_{2,-}^{s,b}$ and $\overline{u}$ in $X_{2,+}^{s,b}$. We define the working space for real-valued pairs $(z, w)$, with $w$ of zero spatial mean, by
\begin{equation}\label{zworkspace}
\mathcal{X}_2^{s,b}=\big\{(z,w):z-i|\nabla|^{-1}w\in X_2^{s,b}\big\},
\qquad\|(z,w)\|_{\mathcal X_2^{s,b}}=\|z-i|\nabla|^{-1}w\|_{X_2^{s,b}}.
\end{equation}
The standard embedding for half-wave Bourgain spaces therefore yields
$$
\mathcal X_2^{s,b}\hookrightarrow C\bigl(\mathbb{R};H^s\times H^{s-1}\bigr),\qquad b>1/2.
$$
In particular, $(z,z_t)\in\mathcal X_2^{s,b}$ implies $z\in C(\mathbb{R};H^s)\cap C^1(\mathbb{R};H^{s-1})$.

\begin{lemma}\label{contiprop}
For $s,~s^\prime,~b\in\mathbb{R}$ with $s^\prime>s$, $\eta\in\mathcal{C}_c^\infty$, and for $v$ supported in $[-2T,2T]$, we have:
\begin{align}
\label{pitau1}&\|(I-\Pi_\theta)v\|_{X_j^{s,b}}+\|\Pi_\theta v\|_{X_j^{s,b}}\lesssim \|v\|_{X_j^{s,b}},\\
\label{pitau2}&\|(I-\Pi_\theta)v\|_{X_j^{s,b}}\lesssim \theta^{\frac{s^\prime-s}{2}}\|v\|_{X_j^{s^\prime,b}},\\
\label{pitau3}&\|\Pi_\theta v\|_{X_j^{s^\prime,b}}\lesssim\theta^{-\frac{s^\prime-s}2}
\|v\|_{X_j^{s,b}},\quad j=1,2.
\end{align}
Moreover, for $b\in(\frac12,1)$ and $b'\in(0,\frac12)$ with $b+b'\leq1$, one has
\begin{align}
\label{linf}&\|v\|_{L^\infty H^s}\lesssim\|v\|_{X_j^{s,b}},\\
\label{b-1}&\Big\|\eta(t)\int_0^t\fe^{i(t-t')A_j}v(t')dt'\Big\|_{X_j^{s,b}}\lesssim T^{1-b-b'}\|v\|_{X_j^{s,-b'}},\quad j=1,2,
\end{align}
where $A_1=\Delta$ and $A_2=|\nabla|$.
\end{lemma}

Estimates \eqref{pitau1}-\eqref{pitau3} are fundamental properties of $\Pi_\theta$ and can be proved directly from the definition. Proofs of \eqref{linf}--\eqref{b-1} can be found in \cite[Section~2.6]{Tao}. In particular, \eqref{b-1} follows from several estimates in \cite[Section~2.6]{Tao} (see also \cite[Equation 9]{Jiima}).

In the following lemma, we present the multilinear estimates in Bourgain spaces related to the ZS.
\begin{lemma}\label{contimult}
Let $s_0$ be defined in \eqref{para}. For any $s_2>s_0$ and $\max\{\frac38,\frac12-\frac12(s_2-s_0)\}<b<1/2$, and for functions $v$ and $w$ supported in $[-2T,2T]$, the following estimates hold:
\begin{align}
\label{eqmult1}&\|vw\|_{X_1^{s_2+\frac12,-b}}+\|\overline{v}w\|_{X_1^{s_2+\frac12,-b}}\lesssim\|v\|_{X_2^{s_2,b}}\|w\|_{X_1^{s_2+\frac12,b}},\\	
\label{eqmult2}&\|\overline{v}w\|_{X_2^{s_2+1,-b}}+\|v\overline{w}\|_{X_2^{s_2+1,-b}}\lesssim\|v\|_{X_1^{s_2+\frac12,b}}\|w\|_{X_1^{s_2+\frac12,b}}.
\end{align}
\end{lemma}

For the proof, the principle arguments for $d\geq2$ are contained in \cite{Kishijam}; the case $d=1$ follows by the same estimates adapted to one dimension. The parameters $s_2$ and $b_1$ introduced in \eqref{para} are admissible in this lemma. Following the approach in \cite{Jiima}, and to avoid the loss of derivatives when $s_2>s_0$, we also use the following estimates:
\begin{align}
&\label{eqmult3}\|vw\|_{X_1^{s_0+\frac12,-b}}+\|\overline{v}w\|_{X_1^{s_0+\frac12,-b}}\lesssim\|v\|_{X_2^{s_0,b}}\|w\|_{X_1^{s_2+\frac12,b}},\\
&\label{eqmult4}\|vw\|_{X_1^{s_0+\frac12,-b}}+\|\overline{v}w\|_{X_1^{s_0+\frac12,-b}}\lesssim\|v\|_{X_2^{s_2,b}}\|w\|_{X_1^{s_0+\frac12,b}},\\
&\label{eqmult5}\|\overline{v}w\|_{X_2^{s_0+1,-b}}+\|v\overline{w}\|_{X_2^{s_0+1,-b}}\lesssim\|v\|_{X_1^{s_0+\frac12,b}}\|w\|_{X_1^{s_2+\frac12,b}}.
\end{align}
We shall prove the discrete analogues of these estimates,  \eqref{discmult3}-\eqref{discmult5}, in Section~\ref{sectionproof}; the same techniques introduced in Section~\ref{sectionproof} yield the continuous versions \eqref{eqmult3}-\eqref{eqmult5}.

Subsequently, we present the well-posedness result.

\begin{theorem}\label{theoexist}
Let the parameters $s_0,~s_1,~s_2,~b_0$ be as in \eqref{para}. For initial data $(E_0,z_0,z_1)\in H^{s_2+\frac12}\times H^{s_2}\times H^{s_2-1}$, there exists $T>0$ and a unique solution $(E,z, z_t)\in X_1^{s_2+\frac12,b_0}\times \mathcal{X}_2^{s_2,b_0}$ of the ZS \eqref{zs} on $[0,T]$. Moreover, for the filtered system \eqref{filtzs}, there exists a unique solution $(E^\theta,z^\theta,z_t^\theta)$ in the same spaces and the following estimates hold:
\begin{align}	
\label{contbound}\|E\|_{X_1^{s_2+\frac12,b_0}}+\|(z,z_t)\|_{\mathcal{X}_2^{s_2,b_0}}&\leq C_T,\\	
\label{contboundtau}\|E^\theta\|_{X_1^{s_2+\frac12,b_0}}+\|(z^\theta,z_t^\theta)\|_{\mathcal{X}_2^{s_2,b_0}}&\leq C_T,\\	
\label{minusest}\|E-E^\theta\|_{X_1^{s_0+\frac12,b_0}}+\|(z-z^\theta,z_t-z_t^\theta)\|_{\mathcal{X}_2^{s_0,b_0}}&\leq C_T\theta^{\frac{s_1}2},
\end{align}
where the constant $C_T>0$ depends solely on $T,~s_1,~b_0$ and the norms $\|E_0\|_{H^{s_2+\frac12}}$, $\|z_0\|_{H^{s_2}}$, $\|z_1\|_{H^{s_2-1}}$.
\end{theorem}
\begin{proof}
It is noteworthy that all functions are defined globally in time but solve the respective equations only on $[0,T]$ (cf. \cite[Remark~2.7]{Jiima}).

Existence of $(E,z,z_t)$ and the estimate in \eqref{contbound} can be found in \cite{Kishijam}. Local well-posedness for the filtered system \eqref{filtzs} follows similarly from \eqref{pitau1}. The key issue here is which of the two equations admits the longer maximal existence interval. Indeed, for each fixed $\theta$, the filtered system \eqref{filtzs} is finite-dimensional. Since $\Pi_\theta$ is self-adjoint, $E^\theta=\Pi_\theta E^\theta$, and $\operatorname{Re}u^\theta$ is real-valued, the system conserves the Schr\"odinger mass:
$$
\frac{d}{dt}\|E^\theta(t)\|_{L^2}^2=2\re\big\langle-i\Pi_\theta\bigl((\operatorname{Re}u^\theta)E^\theta\bigr),E^\theta\big\rangle=0.
$$
Combining this conservation law with Bernstein's inequality on the retained frequencies and Duhamel's formula for the wave component yields
$$
\|u^\theta(t)\|_{L^2}\leq\|\Pi_\theta u_0\|_{L^2}+C_d t\,\theta^{-1/2-d/4}\|\Pi_\theta E_0\|_{L^2}^2,\qquad 0<\theta\leq1.
$$
These bounds rule out finite-time blow-up and therefore imply global well-posedness for each fixed $\theta$ (see also \cite{ggkz}). Although this bound depends on $\theta$, the comparison argument below establishes \eqref{contboundtau} uniformly in $\theta$ on $[0,T]$.

We now prove \eqref{minusest}. Recalling $\eqref{ztou}$, it suffices to demonstrate
$$
\|E-E^\theta\|_{X_1^{s_0+\frac12,b_0}}+\|u-u^\theta\|_{X_2^{s_0,b_0}}\leq C_T\theta^{\frac{s_1}2}.
$$
To accomplish this, we write Duhamel's formula for $E$ and $E^\theta$ with a time cutoff $\eta$:
\begin{align}
\label{duhe}E(t)&=\eta(t)\fe^{it\Delta}E_0-\frac i2\eta(t)\int_0^t\fe^{i(t-s)\Delta}(u(s)+\overline{u}(s))E(s)ds,\\
\label{duhet}E^\theta(t)&=\eta(t)\fe^{it\Delta}\Pi_\theta E_0-\frac i2\eta(t)\int_0^t\fe^{i(t-s)\Delta}\Pi_\theta\big((\Pi_\theta u^\theta(s)+\Pi_\theta\overline{u^\theta}(s))\Pi_\theta E^\theta(s)\big)ds,
\end{align}

On $[0,T_1]$ with $T_1\leq1$ to be chosen below, we use a fixed smooth cutoff $\eta$ equal to one on $[-T_1,T_1]$ in the free initial-data term and take restriction norms on this interval, so that \eqref{b-1} yields the factor $T_1^{1-b_0-b_1}$ while the constant in \eqref{iniest} remains independent of $T_1$.

Note that we have the standard estimate
\begin{equation}\label{iniest}
\|\eta(t)\fe^{it\Delta}E_0\|_{X_1^{s,b_0}}\lesssim \|E_0\|_{H^s},
\end{equation}
which can be found in \cite[Section~2.6]{Tao}. By comparing \eqref{duhe} and \eqref{duhet}, and utilizing \eqref{b-1} and \eqref{iniest}, we obtain
\begin{align*}
\|E-E^\theta\|_{X_1^{s_0+\frac12,b_0}}&\lesssim T_1^{\varepsilon_0}\big(\|(I-\Pi_\theta)((u+\overline{u})E)\|_{X_1^{s_0+\frac12,-b_1}}\\
&\quad+\|\Pi_\theta((u+\overline{u})(E-E^\theta))\|_{X_1^{s_0+\frac12,-b_1}}+\|\Pi_\theta((u-u^\theta+\overline{u}-\overline{u^\theta})E)\|_{X_1^{s_0+\frac12,-b_1}}\\
&\quad+\|\Pi_\theta((u-u^\theta+\overline{u}-\overline{u^\theta})(E-E^\theta))\|_{X_1^{s_0+\frac12,-b_1}}\big)+\|(I-\Pi_\theta)E_0\|_{H^{s_0+\frac12}}.
\end{align*}
Utilizing the properties of $\Pi_\theta$ \eqref{pitau1}-\eqref{pitau3} (noting $\Pi_\theta u^\theta=u^\theta$ and $\Pi_\theta E^\theta=E^\theta$), we derive that
\begin{align*}
\|E-E^\theta\|_{X_1^{s_0+\frac12,b_0}}&\lesssim T_1^{\varepsilon_0}\big(\theta^{\frac{s_1}2}\|(u+\overline{u})E\|_{X_1^{s_2+\frac12,-b_1}}+\|(u+\overline{u})(E-E^\theta)\|_{X_1^{s_0+\frac12,-b_1}}\\
&\quad+\|(u-u^\theta+\overline{u}-\overline{u^\theta})E\|_{X_1^{s_0+\frac12,-b_1}}+\|((I-\Pi_\theta)(u+\overline{u}))(E-E^\theta)\|_{X_1^{s_0+\frac12,-b_1}}\\
&\quad+\|(\Pi_\theta u-u^\theta+\Pi_\theta\overline{u}-\overline{u^\theta})(E-E^\theta))\|_{X_1^{s_0+\frac12,-b_1}}\big)+\theta^{\frac{s_1}2}.
\end{align*}
Moreover, employing the nonlinear estimates \eqref{eqmult1}, \eqref{eqmult3} and \eqref{eqmult4}, and recalling \eqref{contbound}, we arrive at
\begin{equation}\label{e-et}
\begin{aligned}
\|E-E^\theta\|_{X_1^{s_0+\frac12,b_0}}&\lesssim\theta^{\frac{s_1}2}+T_1^{\varepsilon_0}\big(\|u-u^\theta\|_{X_2^{s_0,b_0}}+\|E-E^\theta\|_{X_1^{s_0+\frac12,b_0}}\\
&\quad+\|\Pi_\theta(u-u^\theta)\|_{X_2^{s_0+s_3,b_0}}\|E-E^\theta\|_{X_1^{s_0+\frac12,b_0}}\big)\\
&\lesssim\theta^{\frac{s_1}2}+T_1^{\varepsilon_0}\big(\|u-u^\theta\|_{X_2^{s_0,b_0}}+\|E-E^\theta\|_{X_1^{s_0+\frac12,b_0}}\\
&\quad+\theta^{-\frac{s_3}2}\|\Pi_\theta(u-u^\theta)\|_{X_2^{s_0,b_0}}\|E-E^\theta\|_{X_1^{s_0+\frac12,b_0}}\big),
\end{aligned}
\end{equation}
where $s_3\in(1-2b_1,s_1)$ and $\varepsilon_0=1-b_0-b_1>0$.

Similarly, using Duhamel's formula for $u$ and $u^\theta$, and applying \eqref{pitau1} and \eqref{pitau3}, we get
\begin{equation}\label{u-ut}
\|u-u^\theta\|_{X_2^{s_0,b_0}}\lesssim\theta^{\frac{s_1}2}+T_1^{\varepsilon_0}\bigg(\|E-E^\theta\|_{X_1^{s_0+\frac12,b_0}}+\theta^{-\frac{s_3}2}\|E-E^\theta\|^2_{X_1^{s_0+\frac12,b_0}}\bigg).\\
\end{equation}

Note that the estimates \eqref{e-et} and \eqref{u-ut} hold uniformly on $[0,T]$. Therefore, by choosing $T_1$ sufficiently small, we obtain the desired estimate on $[0,T_1]$ by the a priori assumption
$$
\|E-E^\theta\|_{X_1^{s_0+\frac12,b_0}}+\|u-u^\theta\|_{X_2^{s_0,b_0}}\lesssim\theta^{\frac{s_3}2},
$$
which is related to the difference of initial data
$$
\|(I-\Pi_\theta)E_0\|_{H^{s_0+\frac12}}+\|(I-\Pi_\theta)u_0\|_{H^{s_0}}\lesssim\theta^{\frac{s_1}2}.
$$

Since the choice of $T_1$ is independent of $u^\theta$ and $E^\theta$, the estimate in \eqref{minusest} then follows by the standard iteration argument (see also~\cite[Proposition~3.5]{Jimcom}).

It remains to establish \eqref{contboundtau}. Indeed, we have
$$
\|E^\theta\|_{X_1^{s_2+\frac12,b_0}}\leq\|\Pi_\theta(E-E^\theta)\|_{X_1^{s_2+\frac12,b_0}}+\|\Pi_\theta E\|_{X_1^{s_2+\frac12,b_0}}\lesssim\theta^{-\frac{s_1}2}\|E-E^\theta\|_{X_1^{s_0+\frac12,b_0}}+\|E\|_{X_1^{s_2+\frac12,b_0}}\lesssim 1.
$$
The estimates for $z$ and $z_t$ follow from \eqref{zworkspace}. This completes the proof.
\end{proof}

We now introduce the discrete Bourgain spaces and present several of their fundamental properties; further details can be found in \cite{Jiima,Jimcom,Ostjems}. Let $\{v_n(x)\}_n$ be a sequence of functions on the torus $\mathbb{T}^d$. We define the time-space Fourier transform of this sequence by
\begin{equation}\label{deff}
\mathcal{F}(v_n)=\widetilde{v_n}(\sigma,k)=\tau\sum\limits_{m\in\mathbb{Z}}\widehat{v_m}(k)\fe^{-im\tau\sigma},\qquad \widehat{v_m}(k)=\int_{\mathbb{T}^d}v_m(x)\fe^{-i\langle k,x\rangle}dx.
\end{equation}
In this context, $\widetilde{v_n}$ is periodic in $\sigma$ with period $2\pi\tau^{-1}$.

The discrete Bourgain spaces $X_{1,\tau}^{s,b}$ and $X_{2,\tau}^{s,b}$ are endowed with the norms:
\begin{align}
\label{discdef1}\|v_n\|_{X_{1,\tau}^{s,b}}&=\|\langle D_\tau\rangle^b\fe^{-in\tau\Delta}v_n\|_{l_\tau^2H_x^s}\sim\|\langle k\rangle^s\langle d_\tau(\sigma+|k|^2)\rangle^b\widetilde{v_n}(\sigma,k)\|_{L^2_\sigma l^2_k((-\frac\pi\tau,\frac\pi\tau)\times\mathbb{Z}^d)},\\
\label{discdef2}\|v_n\|_{X_{2,\tau}^{s,b}}&=\|\langle D_\tau\rangle^b\fe^{-in\tau|\nabla|}v_n\|_{l_\tau^2H_x^s}\sim\|\langle k\rangle^s\langle d_\tau(\sigma-|k|)\rangle^b\widetilde{v_n}(\sigma,k)\|_{L^2_\sigma l^2_k((-\frac\pi\tau,\frac\pi\tau)\times\mathbb{Z}^d)},
\end{align}
where $\big(D_\tau(u_n)\big)_n=\big(\frac{u_n-u_{n-1}}{\tau}\big)_n$, $d_\tau(\sigma)=\frac{\fe^{i\tau\sigma}-1}{\tau}$. The equivalence of the two norms on each space was proved in \cite{Jiima} and \cite{Ostjems}. Moreover, for $j=1,2$, we have the continuous embeddings
\begin{equation}\label{monotone}
X_{j,\tau}^{s,b}\subset X_{j,\tau}^{s^\prime,b^\prime},\quad \mathrm{for}\quad s\geq s^\prime,~b\geq b^\prime.
\end{equation}
Regarding the discrete Bourgain spaces, we present the following properties.
\begin{lemma}
For every sequence $\{v_n\}_{n\in\mathbb{Z}}$ supported on $n\tau\in[-2T,2T]$, we have the following estimates:
\begin{align}
\label{shift}\sup\limits_{\delta\in[-3,3]}\|\fe^{i\tau\delta A_j} v_n\|_{X_{j,\tau}^{s,b}}&\lesssim\|v_n\|_{X_{j,\tau}^{s,b}},\quad s,b\in\mathbb{R},\\
\label{s-sp}\|\Pi_\theta v_n\|_{X_{j,\tau}^{s,b}}&\lesssim\theta^{\frac{s^\prime-s}{2}}\|\Pi_\theta v_n\|_{X_{j,\tau}^{s^\prime,b}},\quad s\geq s^\prime,\\
\label{s-sp2}\|(I-\Pi_\theta)v_n\|_{X_{j,\tau}^{s^\prime,b}}&\lesssim\theta^{\frac{s-s^\prime}{2}}\|v_n\|_{X_{j,\tau}^{s,b}},\quad s\geq s^\prime,\\
\label{b-bp}\|v_n\|_{X_{j,\tau}^{s,b}}&\lesssim\tau^{b^\prime-b}\|v_n\|_{X_{j,\tau}^{s,b^\prime}},\quad b\geq b^\prime,\\
\label{disclinf}\|v_n\|_{l_\tau^\infty H^s}&\lesssim\|v_n\|_{X_{j,\tau}^{s,b}},\quad b>\tfrac12,\\
\label{discb-1}\Big\|\tau\sum\limits_{m=0}^n\fe^{i(n-m)\tau A_j}v_m(x)\Big\|_{X_{j,\tau}^{s,b}}&\lesssim T^{1-b-b^\prime}\|v_n\|_{X_{j,\tau}^{s,-b^\prime}},\quad b\in(\tfrac12,1),~b^\prime\in(0, 1-b),
\end{align}
where $\theta\gtrsim\tau$, $j=1,2$ and $A_1=\Delta$, $A_2=|\nabla|$.
\end{lemma}

The proof for the case $j=1$ appears in \cite{Ostjems}; the case $j=2$ can be obtained by the same argument and is omitted here for brevity.

At the end of this section, we provide the multilinear estimate for the ZS in discrete Bourgain spaces.

\begin{theorem}\label{theomult}
Let $s_0=\max\{0,\frac d2-1\}$. For any $s_2>s_0$ and  $\max\{\frac38,\frac12-\frac12(s_2-s_0)\}<b<\frac12$, and for any $\theta\geq\tau$, we have
\begin{align}
\label{discmult1}\|\Pi_\theta v_n\Pi_\theta w_n\|_{X_{1,\tau}^{s_2+\frac12,-b}}+\|\Pi_\theta\overline{v_n}\Pi_\theta w_n\|_{X_{1,\tau}^{s_2+\frac12,-b}}&\lesssim\|v_n\|_{X_{2,\tau}^{s_2,b}}\|w_n\|_{X_{1,\tau}^{s_2+\frac12,b}},\\
\label{discmult3}\|\Pi_\theta v_n\Pi_\theta w_n\|_{X_{1,\tau}^{s_0+\frac12,-b}}+\|\Pi_\theta\overline{v_n}\Pi_\theta w_n\|_{X_{1,\tau}^{s_0+\frac12,-b}}&\lesssim\|v_n\|_{X_{2,\tau}^{s_0,b}}\|w_n\|_{X_{1,\tau}^{s_2+\frac12,b}},\\
\label{discmult4}\|\Pi_\theta v_n\Pi_\theta w_n\|_{X_{1,\tau}^{s_0+\frac12,-b}}+\|\Pi_\theta\overline{v_n}\Pi_\theta w_n\|_{X_{1,\tau}^{s_0+\frac12,-b}}&\lesssim\|v_n\|_{X_{2,\tau}^{s_2,b}}\|w_n\|_{X_{1,\tau}^{s_0+\frac12,b}},\\
\label{discmult2}\|\Pi_\theta v_n\Pi_\theta\overline{w_n}\|_{X_{2,\tau}^{s_2+1,-b}}+\|\Pi_\theta\overline{v_n}\Pi_\theta w_n\|_{X_{2,\tau}^{s_2+1,-b}}&\lesssim\|v_n\|_{X_{1,\tau}^{s_2+\frac12,b}}\|w_n\|_{X_{1,\tau}^{s_2+\frac12,b}},\\
\label{discmult5}\|\Pi_\theta v_n\Pi_\theta\overline{w_n}\|_{X_{2,\tau}^{s_0+1,-b}}+\|\Pi_\theta
\overline{v_n}\Pi_\theta w_n\|_{X_{2,\tau}^{s_0+1,-b}}&\lesssim\|v_n\|_{X_{1,\tau}^{s_2+\frac12,b}}\|w_n\|_{X_{1,\tau}^{s_0+\frac12,b}},
\end{align}
where $\{v_n\}$ and $\{w_n\}$ are arbitrary sequences belonging to the appropriate spaces.
\end{theorem}
The proof is lengthy and would interrupt the subsequent error analysis, we postpone it to Appendix.

\section{Local error analysis}\label{sectionlocal}

In this section we analyze the local error of the fully discretized filtered Lie splitting method \eqref{filtlie} when applied to the system~\eqref{zseqv}. We begin by demonstrating the boundedness of the exact solution in the discrete Bourgain space. In the next lemma, we improve the result presented in \cite{Jiima} and \cite{Roupaa}, which leads to the boundedness straightforwardly.

\begin{lemma}
For any $s\geq0$ and $b>\frac12$, let $\{v_n(x)\}_{n\in\mathbb{Z}}$ be the time-sampled sequence $v_n(x)=v(n\tau,x)$. Then it holds that
\begin{equation}\label{boundutau}
\|v_n\|_{X_{j,\tau}^{s,b}}\lesssim\|v\|_{X_j^{s,b}},\quad j=1,2.
\end{equation}
\end{lemma}
\begin{proof}
Following the argument in \cite[Lemma~3.4]{Roupaa}, we set $f_k(t)=\langle k\rangle^s\fe^{-ita_j(k)}\widehat v(t,k)$, where $a_1(k)=-|k|^2$ and $a_2(k)=|k|$. By Poisson's summation formula and the periodicity of $d_\tau$ we obtain
$$
\langle d_\tau(\sigma)\rangle^b\widetilde{f_k(n\tau)}(\sigma)=\sum_{m\in\mathbb{Z}}\Big\langle d_\tau\Big(\sigma+\frac{2\pi m}\tau\Big)\Big\rangle^b\widehat f_k\Big(\sigma+\frac{2\pi m}\tau\Big).
$$
Separating the terms $m=0$ and $m\neq0$, and using $|d_\tau(\sigma)|\lesssim|\sigma|$ and $|d_\tau(\sigma+\frac{2\pi m}\tau)|\lesssim\tau^{-1}$, we obtain by Cauchy--Schwarz inequality,
$$
\big|\big\langle d_\tau(\sigma)\big\rangle^b\widetilde{f_k(n\tau)}(\sigma)\big|^2\lesssim\big|\langle\sigma\rangle^b\widetilde{f_k}(\sigma)\big|^2+\sum\limits_{m\neq0}\tau^{-2b}\Big\langle\sigma+\frac{2\pi m}\tau\Big\rangle^{-2b}\sum\limits_{m\neq0}\Big|\Big\langle\sigma+\frac{2\pi m}\tau\Big\rangle^b\widehat f_k\Big(\sigma+\frac{2\pi m}\tau\Big)\Big|^2.
$$
Since $|\tau\sigma+2\pi m|\geq\pi(2|m|-1)$ for $|\sigma|\leq\pi/\tau$, the first sum is uniformly bounded for $b>1/2$. Integrating over $\sigma\in[-\pi/\tau,\pi/\tau]$ and summing over $k$ yields the desired estimate.
\end{proof}

Combining estimates \eqref{contboundtau} and \eqref{boundutau} then yields the boundedness of the solution $(u^\theta,E^\theta)$ of~\eqref{filtzs}:
\begin{equation}\label{discbounde}
\sup_{\vartheta\in[-3\tau,3\tau]}\|E^\theta(t_n+\vartheta)\|_{X_{1,\tau}^{s_2+\frac12,b_0}}+ \sup_{\vartheta\in[-3\tau,3\tau]}\|u^\theta(t_n+\vartheta)\|_{X_{2,\tau}^{s_2,b_0}}\leq C_T,
\end{equation}
where $s_2$ and $b_0$ are defined in \eqref{para}.

With these preparatory bounds we proceed to compute the temporal local error. Applying \eqref{filtlie} and Duhamel's formula produces the following decomposition of the local error:
\begin{align}
\mathcal{E}_{eloc}(t_n)&=\Phi^\tau_N(E^\theta(t_n),u^\theta(t_n))-E^\theta(t_{n+1})\nonumber\\
&=\fe^{i\tau\Delta}\Pi_\theta(T_N-I)(\fe^{-\tfrac i2\tau\Pi_\theta (u^\theta(t_n)+\overline{u^\theta}(t_n))}\Pi_\theta E^\theta(t_n))\nonumber\\
&\quad+\fe^{i\tau\Delta}\Pi_\theta(\fe^{-\tfrac i2\tau\Pi_\theta (u^\theta(t_n)+\overline{u^\theta}(t_n))}\Pi_\theta E^\theta(t_n))-\fe^{i\tau\Delta}\Pi_\theta E^\theta(t_n)\nonumber\\
&\quad+\frac i2\int_0^\tau\fe^{i(\tau-\vartheta)\Delta}\Pi_\theta\big(\Pi_\theta( u^\theta(t_n+\vartheta)+\overline{u^\theta}(t_n+\vartheta))\Pi_\theta E^\theta(t_n+\vartheta)\big)d\vartheta\nonumber\\
&=\mathcal{E}_{esloc}(t_n)+\frac i2\fe^{i\tau\Delta}\Big(\int_0^\tau(\fe^{-i\vartheta\Delta}-1)\Pi_\theta\big(\Pi_\theta( u^\theta(t_n+\vartheta)+\overline{u^\theta}(t_n+\vartheta))\Pi_\theta E^\theta(t_n+\vartheta)\big)d\vartheta\nonumber\\
&\quad+\int_0^\tau\Pi_\theta\big(\Pi_\theta( u^\theta(t_n+\vartheta)+\overline{u^\theta}(t_n+\vartheta))\Pi_\theta (E^\theta(t_n+\vartheta)- E^\theta(t_n))\big)d\vartheta\nonumber\\
&\quad+\int_0^\tau\Pi_\theta\big(\Pi_\theta (u^\theta(t_n+\vartheta)- u^\theta(t_n)+\overline{u^\theta}(t_n+\vartheta)-\overline{u^\theta}(t_n))\Pi_\theta E^\theta(t_n)\big)d\vartheta\nonumber\\
&\quad-2i\Pi_\theta\bigl(\big(\fe^{-\tfrac i2\tau\Pi_\theta (u^\theta(t_n)+\overline{u^\theta}(t_n))}-1+\tfrac i2\tau\Pi_\theta (u^\theta(t_n)+\overline{u^\theta}(t_n))\big)\Pi_\theta E^\theta(t_n)\bigr)\!\Bigr)\nonumber\\
&=\mathcal{E}_{esloc}(t_n)+\fe^{i\tau\Delta}\big(\mathcal{E}_1(t_n)+\mathcal{E}_2(t_n)+
\mathcal{E}_3(t_n)+\mathcal{E}_4(t_n)\big)\nonumber\\
&=\mathcal{E}_{esloc}(t_n)+\mathcal{E}_{etloc}(t_n),\label{etloc}
\end{align}
\begin{align}
\mathcal{E}_{uloc}(t_n)&=\Psi^\tau_N(E^\theta(t_n),u^\theta(t_n))-
u^\theta(t_{n+1})\nonumber\\
&=i\tau|\nabla|\fe^{i\tau|\nabla|}\Pi_\theta(T_N-I)(\Pi_\theta E^\theta(t_n)\Pi_\theta\overline{E^\theta}(t_n))\nonumber\\
&\quad+\fe^{i\tau|\nabla|}(i\tau|\nabla|\Pi_\theta(\Pi_\theta E^\theta(t_n)\Pi_\theta\overline{E^\theta}(t_n))+\Pi_\theta u^\theta(t_n))-\fe^{i\tau|\nabla|}\Pi_\theta u^\theta(t_n)\nonumber\\
&\quad-i|\nabla|\int_0^\tau\fe^{i(\tau-\vartheta)|\nabla|}\Pi_\theta\big(\Pi_\theta E^\theta(t_n+\vartheta)\Pi_\theta\overline{E^\theta}(t_n+\vartheta)\big)d\vartheta\nonumber\\
&=\mathcal{E}_{usloc}(t_n)-i|\nabla|\fe^{i\tau|\nabla|}\int_0^\tau(\fe^{-i\vartheta|\nabla|}-1)\Pi_\theta\big(\Pi_\theta E^\theta(t_n+\vartheta)\Pi_\theta\overline{E^\theta}(t_n+\vartheta)\big)d\vartheta\nonumber\\
&\quad-i|\nabla|\fe^{i\tau|\nabla|}\int_0^\tau\Pi_\theta\big(\Pi_\theta(E^\theta(t_n+\vartheta)-E^\theta(t_n))\Pi_\theta\overline{E^\theta}(t_n+\vartheta)\big)d\vartheta\nonumber\\
&\quad-i|\nabla|\fe^{i\tau|\nabla|}\int_0^\tau\Pi_\theta\big(\Pi_\theta E^\theta(t_n)\Pi_\theta(\overline{E^\theta}(t_n+\vartheta)-\overline{E^\theta}(t_n))\big)d\vartheta\nonumber\\
&=\mathcal{E}_{usloc}(t_n)+\fe^{i\tau|\nabla|}\big(\mathcal{E}_5(t_n)+\mathcal{E}_6(t_n)+
\mathcal{E}_7(t_n)\big)\nonumber\\
&=\mathcal{E}_{usloc}(t_n)+\mathcal{E}_{utloc}(t_n).\label{utloc}
\end{align}
Specifically, the solution-increment term $\mathcal{E}_2$ and the
exponential-remainder term $\mathcal{E}_4$ in \eqref{etloc} together account for the difference between the projected nonlinear evolution and the projected phase step. The former retains the nested projection arising from the filtered equation, whereas the
latter retains the remainder of the pointwise exponential.
Similarly, $\mathcal{E}_6$ and $\mathcal{E}_7$ in \eqref{utloc} account for the variation
of the quadratic wave source that is frozen in the numerical
update. These contributions are included in the temporal local
defect treated in Proposition~1. The additional terms involving
$T_N-I$ are identified separately as spatial interpolation defects.

In order to establish the local error, we present a lemma pertaining to the well-known $\varphi$-functions \cite{Hochacta}:

\begin{lemma}\label{lemmaphik}
For any $k\in\mathbb{N}$, we define $\varphi_k$ as
$$
\varphi_0(\alpha)=\fe^\alpha,\quad \varphi_{k+1}(\alpha)=\frac{\varphi_k(\alpha)-1/k!}{\alpha},\quad k\ge 0.
$$
Suppose $\theta\geq\tau$, and $s,p\geq0$ with $s+p>\frac d2$.  Then the following estimates hold:
\begin{align}
\label{phik1}\|\varphi_k(\tau\Pi_\theta u)\|_{H^s}&\lesssim\varphi_k(\tau^{1-\frac p2}\|\Pi_\theta u\|_{H^s}),\\
\label{phik2}\|\varphi_k(\tau\Pi_\theta u+\tau\Pi_\theta v)\|_{H^s}&\lesssim\varphi_k(2\tau^{1-\frac p2}\|\Pi_\theta u\|_{H^s})+\varphi_k(2\tau^{1-\frac p2}\|\Pi_\theta v\|_{H^s}).
\end{align}
\end{lemma}
\begin{proof}
Since $\Pi_\theta u$ and $\Pi_\theta v$ possess only finitely many Fourier modes they are smooth. Taylor expansion together with the classical bilinear estimate \eqref{cbe} yields
\begin{align*}
\|\varphi_k(\tau\Pi_\theta u)\|_{H^s}&\leq\|\varphi_k(\tau\Pi_\theta u)\|_{H^{s+p}}\leq\sum\limits_{j=k}^\infty\frac{\tau^{j-k}}{j!}\|(\Pi_\theta u)^{j-k}\|_{H^{s+p}}\leq\sum\limits_{j=k}^\infty\frac{\tau^{j-k}}{j!}\|\Pi_\theta u\|_{H^{s+p}}^{j-k}.
\end{align*}
Noticing that $\theta\geq\tau$, utilizing \eqref{s-sp} gives the asserted estimate \eqref{phik1}:
$$
\sum\limits_{j=k}^\infty\frac{\tau^{j-k}}{j!}\|\Pi_\theta u\|_{H^{s+p}}^{j-k}\lesssim\sum\limits_{j=k}^\infty\frac{\tau^{\frac{(j-k)(2-p)}2}}{j!}\|\Pi_\theta u\|_{H^s}^{j-k}=\varphi_k(\tau^{1-\frac p2}\|\Pi_\theta u\|_{H^s}).
$$
For \eqref{phik2}, applying \eqref{phik1}, it therefore suffices to show
$$
\varphi_k(a+b)\lesssim\varphi_k(2a)+\varphi_k(2b),\quad a,b\geq0. $$
Indeed, by employing Taylor's expansion and the power-mean inequality, we find
\begin{align*} \varphi_k(a+b)&=\sum\limits_{j=k}^\infty\frac{(a+b)^{j-k}}{j!}\leq\sum\limits_{j=k}^\infty\frac{2^{j-k-1}(a^{j-k}+b^{j-k})}{j!}\\
&=\frac12\sum\limits_{i=k}^\infty\frac{(2a)^{j-k}}{j!}+
\frac12\sum\limits_{j=k}^\infty\frac{(2b)^{j-k}}{j!}=
\frac12\varphi_k(2a)+\frac12\varphi_k(2b),
\end{align*}
which completes the proof.
\end{proof}

We now estimate the temporal component of the local error, followed by the spatial component.
\begin{proposition}\label{proptlocal}
Let $s_0,~s_1,~b_0$ be as in \eqref{para}, let $(E^\theta,u^\theta)$ be the solution from  Theorem~\ref{theoexist}, and set $\theta$ defined in \eqref{theta} sufficiently small. Then the temporal local error satisfies:
\begin{equation}\label{tlocal}		
\|\mathcal{E}_{etloc}(t_n)\|_{X_{1,\tau}^{s_0+\frac12,b_0-1}}+\|\mathcal{E}_{utloc}(t_n)\|_{X_{2,\tau}^{s_0,b_0-1}}\leq C_T\tau\theta^{\frac{s_c}2},
\end{equation}
where $s_c=\min\{2,s_1\}$.
\end{proposition}
\begin{proof}
We assume without loss of generality that $s_1\leq2$, thus $s_c=s_1$. Utilizing \eqref{shift}, \eqref{utloc} and \eqref{etloc}, it suffices to demonstrate
\begin{align}
&\|\mathcal{E}_j(t_n)\|_{X_{1,\tau}^{s_0+\frac12,b_0-1}}\leq C_T\tau\theta^{\frac{s_1}2},\quad j=1,2,3,4;\label{j1234}\\
&\|\mathcal{E}_j(t_n)\|_{X_{2,\tau}^{s_0,b_0-1}}\leq C_T\tau\theta^{\frac{s_1}2},\quad j=5,6,7.\label{j567}
\end{align}
We begin by estimating $\mathcal{E}_1$. By \eqref{shift}, \eqref{discmult1}, \eqref{discbounde}, we obtain
\begin{equation}\label{e1}
\begin{aligned}
&\quad\|\mathcal{E}_1(t_n)\|_{X_{1,\tau}^{s_0+\frac12,b_0-1}}\\
&\lesssim\tau\sup_{\vartheta\in[0,\tau]}\big\|(\tau\Delta)^{\frac{s_1}2}\Pi_\theta\big(\Pi_\theta(u^\theta(t_n+\vartheta)+\Pi_\theta\overline{u^\theta}(t_n+\vartheta))\Pi_\theta E^\theta(t_n+\vartheta)\big)\big\|_{X_{1,\tau}^{s_0+\frac12,b_0-1}}\\
&\lesssim\tau^{1+\frac{s_1}2}\sup_{\vartheta\in[0,\tau]}\|u^\theta(t_n+\vartheta)\|_{X_{2,\tau}^{s_2,b_0}}\|E^\theta(t_n+\vartheta)\|_{X_{1,\tau}^{s_2+\frac12,b_0}}\lesssim\tau\theta^{\frac{s_1}2}.
\end{aligned}
\end{equation}
Similar arguments yield \eqref{j567} for $j=5$.

Interpolating \eqref{disclinf} with the trivial result $X_{j,\tau}^{s,0}=l_\tau^2H^s$ $(j=1,2)$ gives, for any $s\in\mathbb{R}$, $b>\frac14$, $b^\prime>\frac13$ and $\{v_n\}_n\in X_{j,\tau}^{\gamma,b_2}$, the auxiliary bound
\begin{equation}\label{discl4}
\|v_n\|_{l_\tau^4 H^s}\lesssim\|v_n\|_{X_{j,\tau}^{s,b}},  \qquad
\|v_n\|_{l_\tau^6 H^s}\lesssim\|v_n\|_{X_{j,\tau}^{s,b^\prime}},\quad j=1,2.
\end{equation}
Applying Duhamel's formula
\begin{align*}
E^\theta(t_n+\theta)-E^\theta(t_n)&=(\fe^{i\theta\Delta}-1)E^\theta(t_n)\\
&\quad-\frac i2\int_0^\theta\fe^{i(\theta-\xi)\Delta}\Pi_\theta\big(\Pi_\theta(u^\theta(t_n+\xi)+\overline{u^\theta}(t_n+\xi))\Pi_\theta E^\theta(t_n+\xi)\big)d\xi,
\end{align*}
and then estimates \eqref{cbe}, \eqref{shift}, \eqref{discmult4} and \eqref{s-sp} yields
\begin{align*}
&\quad\|\mathcal{E}_2(t_n)\|_{X_{1,\tau}^{s_0+\frac12,b_0-1}}\\
&\lesssim\tau\sup_{\vartheta\in[0,\tau]}\big\|\Pi_\theta\big(\Pi_\theta(u^\theta(t_n+\vartheta)+\overline{u^\theta}(t_n+\vartheta))\Pi_\theta(E^\theta(t_n+\vartheta)-E^\theta(t_n))\big)\big\|_{X_{1,\tau}^{s_0+\frac12,b_0-1}}\\
&\lesssim\tau\sup_{\vartheta\in[0,\tau]}\|u^\theta(t_n+\vartheta)\|_{X_{2,\tau}^{s_2,b_0}}\|(\fe^{i\theta\Delta}-1)E^\theta(t_n)\|_{X_{1,\tau}^{s_0+\frac12,b_0}}\\
&\quad+\tau\theta\sup_{\vartheta\in[0,\tau]}\|\Pi_\theta u^\theta(t_n+\vartheta)\|_{l_\tau^\infty H^{s_0+1+\varepsilon}}\|\Pi_\theta u^\theta(t_n+\vartheta)\Pi_\theta E^\theta(t_n+\vartheta)\|_{l_\tau^2 H^{s_0+\frac12}}.
\end{align*}
Combining these with H\"older's inequality, \eqref{disclinf}, \eqref{discbounde} and \eqref{discl4}, note that $X_{1,\tau}^{s,0}=l_\tau^2H^s$, we obtain
\begin{align*}	\|\mathcal{E}_2(t_n)\|_{X_{1,\tau}^{s_0+\frac12,b_0-1}}&\lesssim\tau\sup_{\vartheta\in[0,\tau]}\|u^\theta(t_n+\vartheta)\|_{X_{2,\tau}^{s_2,b_0}}\|(\theta\Delta)^{\frac{s_1}2}E^\theta(t_n)\|_{X_{1,\tau}^{s_0+\frac12,b_0}}\\
&\quad+\tau\theta\sup_{\vartheta\in[0,\tau]}\|\Pi_\theta u^\theta(t_n+\vartheta)\|^2_{X_{2,\tau}^{s_0+1+\varepsilon,b_0}}\|\Pi_\theta E^\theta(t_n+\vartheta)\|_{{X_{1,\tau}^{s_2+\frac12,b_0}}}\\	&\lesssim\tau\theta^{\frac{s_1}2}+\tau\theta^s\sup_{\vartheta\in[0,\tau]}\|\Pi_\theta u^\theta(t_n+\vartheta)\|^2_{X_{2,\tau}^{s_2,b_0}}\|\Pi_\theta E^\theta(t_n+\vartheta)\|_{{X_{1,\tau}^{s_2+\frac12,b_0}}},
\end{align*}
where $\varepsilon>0$ can be chosen arbitrarily small, and $s=\min\{s_1-\varepsilon,1\}$. By taking $\varepsilon<\frac{s_1}2$ for $s_1\leq1$ and $\varepsilon<s_1-1$ for $s_1>1$, we get \eqref{j1234} for $j=2$. Similarly, \eqref{j567} for $j=6, 7$ can be yielded by similar derivations.

Next we estimate $\mathcal{E}_3$. Duhamel's formula gives
$$
u^\theta(t_n+\vartheta)-u^\theta(t_n)=(\fe^{i\vartheta|\nabla|}-1)u^\theta(t_n)+i|\nabla|\int_0^\vartheta\fe^{i(\vartheta-\xi)|\nabla|}\Pi_\theta\big(\Pi_\theta E^\theta(t_n+\xi)\Pi_\theta\overline{E^\theta}(t_n+\xi)\big)d\xi.
$$
We treat the cases $s_1\leq\frac12$ and $s_1>\frac12$ separately. For $s_1\leq\frac12$, by applying \eqref{cbe}, \eqref{shift} and  \eqref{discmult3}, we obtain
\begin{align*}
&\quad\|\mathcal{E}_3(t_n)\|_{X_{1,\tau}^{s_0+\frac12,b_0-1}}\\
&\lesssim\tau\!\sup_{\vartheta\in[0,\tau]}\!\big\|\Pi_\theta\big(\Pi_\theta(u^\theta(t_n+\vartheta)+\overline{u^\theta}(t_n+\vartheta)-u^\theta(t_n)-\overline{u^\theta}(t_n))\Pi_\theta E^\theta(t_n)\big)\big\|_{X_{1,\tau}^{s_0+\frac12,-b_1}}\\
&\lesssim\tau\sup_{\vartheta\in[0,\tau]}\|(\fe^{i\theta|\nabla|}-1)u^\theta(t_n)\|_{X_{2,\tau}^{s_0,b_0}}\|E^\theta(t_n)\|_{X_{1,\tau}^{s_2+\frac12,b_0}}\\
&\quad+\tau\theta\sup_{\vartheta\in[0,\tau]}\|\Pi_\theta E^\theta(t_n)\|_{l_\tau^\infty H^{s_2+1}}\||\nabla|\Pi_\theta(\Pi_\theta E^\theta(t_n+\vartheta)\Pi_\theta \overline{E^\theta}(t_n+\vartheta))\|_{l_\tau^2 H^{s_0+\frac12}},
\end{align*}
which together with H\"older's inequality, \eqref{disclinf},  \eqref{discbounde} and \eqref{discl4} yields
\begin{equation}\label{e31}
\begin{aligned}
\|\mathcal{E}_3(t_n)\|_{X_{1,\tau}^{s_0+\frac12,b_0-1}}&\lesssim\tau\sup_{\vartheta\in[0,\tau]}\|(\theta|\nabla|)^{\frac{s_1}2}u^\theta(t_n+\vartheta)\|_{X_{2,\tau}^{s_0,b_0}}\|E^\theta(t_n)\|_{X_{1,\tau}^{s_2+\frac12,b_0}}\nonumber\\
&\quad+\tau\theta^{\frac12+\frac{s_1}2}\sup_{\vartheta\in[0,\tau]}\|E^\theta(t_n)\|_{X_{1,\tau}^{s_2+\frac12,b_0}}\|\Pi_\theta E^\theta(t_n+\vartheta)\|^2_{{X_{1,\tau}^{s_2+1,b_0}}}\nonumber\\
&\lesssim\tau\theta^{\frac{s_1}2}+\tau\theta^{\frac{s_1}2}\sup_{\vartheta\in[0,\tau]}\|\Pi_\theta E^\theta(t_n+\vartheta)\|^2_{X_{1,\tau}^{s_2+\frac12,b_0}}\lesssim\tau\theta^{\frac{s_1}2}.
\end{aligned}
\end{equation}
Similarly, for the case $s_1>\frac12$, we get
\begin{align*}
&\quad\|\mathcal{E}_3(t_n)\|_{X_{1,\tau}^{s_0+\frac12,b_0-1}}\nonumber\\	&\lesssim\tau\theta^{\frac{s_1}2}+\tau\theta\sup_{\vartheta\in[0,\tau]}\|\Pi_\theta E^\theta(t_n)\|_{l_\tau^\infty H^{s_2+\frac12}}\|\Pi_\theta(\Pi_\theta E^\theta(t_n+\vartheta)\Pi_\theta \overline{E^\theta}(t_n+\vartheta))\|_{l_\tau^2 H^{s_0+\frac32}}\\ &\lesssim\tau\theta^{\frac{s_1}2}+\tau\theta^{\frac12+\frac{s_1}2}\sup_{\vartheta\in[0,\tau]}\|E^\theta(t_n)\|_{X_{1,\tau}^{s_2+\frac12,b_0}}\|\Pi_\theta E^\theta(t_n+\vartheta)\|^2_{{X_{1,\tau}^{s_2+\frac12,b_0}}}\\
&\lesssim\tau\theta^{\frac{s_1}2}+\tau\theta^{\frac{s_1}2}\sup_{\vartheta\in[0,\tau]}\|\Pi_\theta E^\theta(t_n+\vartheta)\|^2_{X_{1,\tau}^{s_2+\frac12,b_0}}\lesssim\tau\theta^{\frac{s_1}2},
\end{align*}
which completes the proof of \eqref{j1234} for $j=3$.

It remains to estimate $\mathcal{E}_4$. Write $F_n=u^\theta(t_n)+\overline{u^\theta}(t_n)$. By \eqref{cbe} and H\"older's inequality, we are led to
\begin{equation}\label{e4st}
\begin{aligned}
\|\mathcal{E}_4(t_n)\|_{X_{1,\tau}^{s_0+\frac12,b_0-1}}&\lesssim\|\fe^{-\frac i2\tau\Pi_\theta F_n}-1+\tfrac i2\tau\Pi_\theta F_n\|_{l_\tau^2H^{s_0+1+\varepsilon}}\|E^\theta(t_n)\|_{l_\tau^\infty H^{s_0+\frac12}}\\
&\lesssim\|\fe^{-\frac i2\tau\Pi_\theta F_n}-1+\tfrac i2\tau\Pi_\theta F_n\|_{l_\tau^2H^{s_0+1+\varepsilon}}\|E^\theta(t_n)\|_{X_{1,\tau}^{s_0+\frac12,b_0}},
\end{aligned}
\end{equation}
where $\varepsilon>0$ can be taken arbitrarily small. Note that for any $\alpha$, we have
$$
\fe^\alpha-1-\alpha=\alpha^2\varphi_2(\alpha)=\tfrac12\alpha^2+\alpha^3\varphi_3(\alpha),
$$
where $\varphi_2$ and $\varphi_3$ are defined in Lemma~\ref{lemmaphik}. Using $s_2+1>\frac d2$, H\"older's inequality, \eqref{cbe} and \eqref{s-sp} yields
\begin{align}
&\quad\|\fe^{-\frac i2\tau\Pi_\theta F_n}-1+\tfrac i2\tau\Pi_\theta F_n\|_{l_\tau^2 H^{s_0+1+\varepsilon}}\nonumber\\
&\lesssim\tau^2\|(\Pi_\theta F_n)^2\|_{l_\tau^2 H^{s_0+1+\varepsilon}}+\tau^3\|(\Pi_\theta F_n)^3\varphi_3(-\tfrac i2\tau\Pi_\theta F_n)\|_{l_\tau^2 H^{s_2+1}}\label{e4c1}\\
&\lesssim\tau^{1+s_1-\varepsilon}\|\Pi_\theta u^\theta(t_n)\|^2_{l_\tau^4 H^{s_2}}+\tau^{\frac32}\|\Pi_\theta u^\theta(t_n)\|^3_{l_\tau^6 H^{s_2}}\|\varphi_3(-\tfrac i2\tau\Pi_\theta F_n)\|_{l_\tau^\infty H^{s_2+1}},\nonumber
\end{align}
for any $0<s_1\leq1$. Consequently, by using \eqref{disclinf}, \eqref{phik1} and \eqref{discl4} we obtain
\begin{equation}\label{e4c2}
\begin{aligned}
&\quad\|\fe^{-\frac i2\tau\Pi_\theta F_n}-1+\tfrac i2\tau\Pi_\theta F_n\|_{l_\tau^2 H^{s_0+1+\varepsilon}}\\
&\lesssim\tau\theta^{\frac{s_1}2}\|\Pi_\theta u^\theta(t_n)\|^2_{X_{2,\tau}^{s_2,b_0}}+\tau\theta^{\frac12}\|\Pi_\theta u^\theta(t_n)\|^3_{X_{2,\tau}^{s_2,b_0}}\varphi_3\big(\tau^\frac13\|\Pi_\theta F_n\|_{X_{2,\tau}^{s_2,b_0}}\big)\\
&\lesssim\tau\theta^{\frac{s_1}2},\quad\quad  s_1\in (0,1].
\end{aligned}
\end{equation}
For $s_1>1$, we choose $\varepsilon>0$ such that $s_2>s_0+1+\varepsilon>\frac d2$. It follows from \eqref{cbe} that we always have by definition that $s_0+1\ge\frac d2$.
\begin{equation}\label{e4c3}
\begin{aligned}
\|\mathcal{E}_4(t_n)\|_{X_{1,\tau}^{s_0+\frac12,b_0-1}}&\lesssim\|(\fe^{-\frac i2\tau\Pi_\theta F_n}-1+\tfrac i2\tau\Pi_\theta F_n)E^\theta(t_n)\|_{l_\tau^2 H^{s_2}}\\
&\lesssim\tau^2\|\Pi_\theta F_n\|^2_{l_\tau^6 H^{s_2}}\|\varphi_2(-\tfrac i2\tau\Pi_\theta F_n)\|_{l_\tau^\infty H^{s_2}}\|E^\theta(t_n)\|_{l_\tau^6 H^{s_2}}\\
&\lesssim\tau\theta\|\Pi_\theta u^\theta(t_n)\|^2_{X_{2,\tau}^{s_2,b_0}}\varphi_2\big(\tau^\frac12\|\Pi_\theta u^\theta(t_n)\|_{X_{2,\tau}^{s_2,b_0}}\big)\|E^\theta(t_n)\|_{X_{1,\tau}^{s_2+\frac12,b_0}},
\end{aligned}
\end{equation}
where for the last inequality we utilized \eqref{monotone}, \eqref{phik1} and \eqref{discl4}. Taking $\varepsilon<\frac{s_1}2$ and combining \eqref{e4st}--\eqref{e4c3} concludes \eqref{j1234} for $j=4$ and the proof is completed.
\end{proof}

Finally, we estimate the spatial contributions to the local error, denoted $\mathcal{E}_{esloc}(t_n)$ and $\mathcal{E}_{usloc}(t_n)$. Similar work has been carried out for the nonlinear Schr\"odinger equation in \cite{Jisiam}, where a spatial local error of order $\mathcal{O}(\tau\theta^{\frac s2})$ was proved for initial data $u_0\in H^s$, $s\in(0,2]$. The estimates were obtained in the discrete Bourgain space $X_{1,\tau}^{0,b-1}$ with proper $b>\frac12$. Here we extend this result, which allows us to estimate the spatial error in discrete Bourgain spaces with non-integer regularity.

We first introduce a technical lemma which involves the interpolation error.
\begin{lemma}
For any $p\in[1,\infty]$, $0\leq s'\leq s$, $b\in\mathbb{R}$ and $s>\frac d2$, $\alpha>0$, we have
\begin{align}
\label{discont}\|(I-T_N)v_n\|_{l_\tau^pH^{s'}}&\lesssim N^{s'-s}\|v_n\|_{l_\tau^pH^s},\\
\label{discontbourg}\|(I-T_N)\Pi_{\alpha\theta}v_n\|_{X_{j,\tau}^{s',b}}&\lesssim_\alpha \theta^{\frac{s-s'}2}\|v_n\|_{X_{j,\tau}^{s,b}},\quad j=1,2.
\end{align}
\end{lemma}
\begin{proof}
We first prove \eqref{discont}. It suffices to prove the pointwise-in-time estimate for any $v\in H^s$ with $s>\frac d2$, it holds
\begin{equation}\label{discontspace}
\|(I-T_N)v\|_{H^{s'}}\lesssim N^{s'-s}\|v\|_{H^s},
\end{equation}
then \eqref{discont} follows by summing over time.

The one-dimensional case is standard (see, for instance, \cite[Theorem~11.8]{Kressspr}). Note that, unlike the symmetric Nyquist-frequency convention in \cite{Kressspr}, our interpolation retains only the frequency $-N/2$ in each coordinate. Thus for $d\geq 2$, the Fourier interpolation operator satisfies
$$
T_Nv(k)=\sum_{p\in\mathbb Z^d}\widehat v(k+Np), \qquad k\in K_N=\{-\tfrac N2,\cdots,\tfrac N2-1\}^d.
$$
By Cauchy--Schwarz inequality,
$$
\Big|\sum_{p\neq0}\widehat v(k+Np)\Big|^2\le\Big(\sum_{p\neq0}\langle k+Np\rangle^{-2s}\Big)\Big(\sum_{p\neq0}\langle k+Np\rangle^{2s}|\widehat v(k+Np)|^2\Big).
$$
Since $2s>d$, the first factor is bounded by $CN^{-2s}$ uniformly with respect to $k\in K_N$. Summing over $k$ and including the Fourier tail yields \eqref{discontspace}, and hence \eqref{discont}.

We now turn to \eqref{discontbourg}. Since $(I-T_N)\Pi_\theta v_n=0$ for all $v_n$, we may assume without loss of generality that $\alpha<1$. We follow the argument in \cite[Lemma~4.2]{Jisiam}.

Taking $m=\lceil\alpha^{-1/2}\rceil$, we have $T_{mN}\Pi_{\alpha\theta}v_n=\Pi_{\alpha\theta}v_n$.
We decompose $\Pi_{\alpha\theta}v_n=\sum_p v_{n,p}$ according to the disjoint frequency sets $G_p=(K_N+Np)\cap K_{mN}$. At most $(m+1)^d$ such sets are nonempty, and no two distinct frequencies within the same set are aliased under $T_N$. Hence, for any $s\ge0$ and $p\in\mathbb{Z}^d$ with $G_p\neq\emptyset$,
$$
\|T_Nv_{n,p}\|_{H^s}^2=\sum_{k\in K_N}\langle k\rangle^{2s}\big|\widehat{v_{n,p}}\big(k+Np\big)\big|^2\le\sum_{k\in K_N}\big\langle k+Np\big\rangle^{2s}\big|\widehat{v_{n,p}}\big(k+Np\big)\big|^2\le\|v_{n,p}\|_{H^s}^2.
$$
Summing over $p$ and using the disjointness of their Fourier supports, we obtain
$$
\big\|T_N\Pi_{\alpha\theta}v_n\big\|_{H^s}\le\sum_p\|T_Nv_{n,p}\|_{H^s}\le\sum_p\|v_{n,p}\|_{H^s}\lesssim_\alpha\big\|\Pi_{\alpha\theta}v_n\big\|_{H^s}.
$$
Following the discrete Parseval argument in \cite[Lemma~4.2]{Jisiam}, we obtain
$$
\big\|\langle D_\tau\rangle^b\fe^{-in\tau A}T_N\Pi_{\alpha\theta}v_n\big\|_{H^s}\lesssim_\alpha\big\|\langle D_\tau\rangle^b\fe^{-in\tau A}\Pi_{\alpha\theta}v_n\big\|_{H^s},
$$
where $A=\Delta$ or $A=|\nabla|$, as appropriate. Taking the $l_\tau^2$ norm in time yields
\begin{equation}\label{tnkeep}
\big\|T_N\Pi_{\alpha\theta}v_n\big\|_{X_{j,\tau}^{s,b}}\lesssim_\alpha\big\|\Pi_{\alpha\theta}v_n\big\|_{X_{j,\tau}^{s,b}},\quad j=1,2.
\end{equation}
Similarly, we derive that
\begin{equation}\label{discobou}
\big\|T_N(\Pi_{\alpha\theta}-\Pi_\theta)v_n\big\|_{X_{j,\tau}^{s',b}}\lesssim_\alpha\big\|(\Pi_{\alpha\theta}-\Pi_\theta)v_n\big\|_{X_{j,\tau}^{s',b}},\quad j=1,2.
\end{equation}
Since $(I-T_N)\Pi_\theta v_n=0$ and $\Pi_{\alpha\theta}\Pi_\theta=\Pi_\theta$, we have
$$
(I-T_N)\Pi_{\alpha\theta}v_n=(\Pi_{\alpha\theta}-\Pi_\theta)v_n-T_N(\Pi_{\alpha\theta}-\Pi_\theta)v_n.
$$
The triangle inequality, together with \eqref{s-sp2} and \eqref{discobou}, yields \eqref{discontbourg}. This completes the proof.
\end{proof}

\begin{proposition}\label{propslocal}
Let $s_0,~s_1,~b_0$ be as in \eqref{para}, let $(E^\theta,u^\theta)$ denote the solution provided by Theorem~\ref{theoexist}, and let $\theta$ \eqref{theta} be sufficiently small. Then the spatial local errors satisfy the estimate
\begin{equation}\label{slocal}
\|\mathcal{E}_{esloc}(t_n)\|_{X_{1,\tau}^{s_0+\frac12,b_0-1}}+\|\mathcal{E}_{usloc}(t_n)\|_{X_{2,\tau}^{s_0,b_0-1}}\le C_T\tau\theta^{\frac{s_c}2},
\end{equation}
where $s_c=\min\{2,s_1\}$.
\end{proposition}
\begin{proof}
We assume without loss of generality that $s_1\leq2$. If $s_1>2$, we simply set $s_1=2$ since $H^{s_0+s_1}\subset H^{s_0+2}$. For the case $s_1\leq2$ so that $s_c=s_1$,  we write $\mathcal{E}_{esloc}(t_n)$ and $\mathcal{E}_{usloc}(t_n)$ as follows:
\begin{align*}
\mathcal{E}_{esloc}(t_n)&=\fe^{i\tau\Delta}\Pi_\theta(T_N-I)(\fe^{-\frac i2\tau(\Pi_\theta u^\theta(t_n)+\Pi_\theta\overline{u^\theta}(t_n))}\Pi_\theta E^\theta(t_n))\\
&=-\frac i2\tau\fe^{i\tau\Delta}\Pi_\theta(T_N-I)\big(\Pi_\theta(u^\theta(t_n)+\overline{u^\theta}(t_n))\Pi_\theta E^\theta(t_n)\big)\\
&\quad+\fe^{i\tau\Delta}\Pi_\theta(T_N-I)\Big(\big(\fe^{-\frac i2\tau\Pi_\theta( u^\theta(t_n)+\overline{u^\theta}(t_n))}-1+\frac i2\tau\Pi_\theta(u^\theta(t_n)+\overline{u^\theta}(t_n))\big)\Pi_\theta E^\theta(t_n)\Big)\\
&=\mathcal{E}_1'(t_n)+\mathcal{E}_2'(t_n);\\
\mathcal{E}_{usloc}(t_n)&=i\tau|\nabla|\fe^{i\tau|\nabla|}\Pi_\theta(T_N-I)(\Pi_\theta E^\theta(t_n)\Pi_\theta\overline{E^\theta}(t_n))=\mathcal{E}_3'(t_n).
\end{align*}

Next we estimate each term separately. It follows from \eqref{discontbourg} that
\begin{align*}
\|\mathcal{E}_1'(t_n)\|_{X_{1,\tau}^{s_0+\frac12,b_0-1}}&\lesssim\tau\|(I-T_N)\big(\Pi_\theta(u^\theta(t_n)+\overline{u^\theta}(t_n))\Pi_\theta E^\theta(t_n)\big)\|_{X_{1,\tau}^{s_0+\frac12,b_0-1}}\\
&\lesssim\tau\theta^{\frac{s_1}2+\frac14}\|\Pi_\theta(u^\theta(t_n)+\overline{u^\theta}(t_n))\Pi_\theta E^\theta(t_n)\|_{X_{1,\tau}^{s_2+1,b_0-1}}.
\end{align*}
Note that the frequency of $\Pi_\theta(u^\theta(t_n)+\overline{u^\theta}(t_n))\Pi_\theta E^\theta(t_n)$ is no greater than $2\theta^{-\frac12}$, thus by \eqref{s-sp}, \eqref{discmult1}, \eqref{discbounde}, we arrive at
$$
\|\mathcal{E}_1'(t_n)\|_{X_{1,\tau}^{s_0+\frac12,b_0-1}}\lesssim\tau\theta^{\frac{s_1}2}\|\Pi_\theta(u^\theta(t_n)+\overline{u^\theta}(t_n))\Pi_\theta E^\theta(t_n)\|_{X_{1,\tau}^{s_2+\frac12,b_0-1}}\lesssim\tau\theta^{\frac{s_1}2}.
$$

Setting $F_n=u^\theta(t_n)+\overline{u^\theta}(t_n)$, and using \eqref{shift}, \eqref{discont}, \eqref{cbe}, \eqref{s-sp}, \eqref{disclinf}, \eqref{discbounde}, \eqref{phik2} and H\"older's inequality, we are led to the desired estimate for $\mathcal{E}_2'(t_n)$:
\begin{align*} \|\mathcal{E}_2'(t_n)\|_{X_{1,\tau}^{s_0+\frac12,b_0-1}}&\lesssim\tau^2\big\|(I-T_N)\big((\Pi_\theta F_n)^2\varphi_2(\tfrac\tau2\Pi_\theta F_n)\Pi_\theta E^\theta(t_n)\big)\big\|_{l_\tau^2H^{s_0+\frac12}}\\
&\lesssim\tau^2\theta^{\frac14+\frac{s_1}2}\|(\Pi_\theta u^\theta(t_n))^2\varphi_2(\tfrac\tau2\Pi_\theta F_n)\Pi_\theta E^\theta(t_n)\|_{l_\tau^2H^{s_2+1}}\\
&\lesssim\tau^2\theta^{\frac14+\frac{s_1}2}\|u^\theta(t_n)\|^2_{l_\tau^\infty H^{s_2+1}}\|E^\theta(t_n)\|_{l_\tau^2 H^{s_2+1}}\|\varphi_2(\tfrac\tau2\Pi_\theta F_n)\|_{l_\tau^\infty H^{s_2+1}}\\
&\lesssim\tau\theta^{\frac{s_1}2}\|u^\theta(t_n)\|^2_{X_{2,\tau}^{s_2,b_0}} \|E^\theta(t_n)\|_{X_{1,\tau}^{s_2+\frac12,b_0}}\varphi_2(\tau^\frac13\|u^\theta(t_n)\|_{X_{2,\tau}^{s_2,b_0}})\lesssim\tau\theta^{\frac{s_1}2}.
\end{align*}

Finally, \eqref{shift} and \eqref{discontbourg} again reduce $\mathcal{E}_3'(t_n)$ to terms controlled by \eqref{discmult2}:
\begin{align*}
\|\mathcal{E}_3'(t_n)\|_{X_{2,\tau}^{s_0,b_0-1}}&\lesssim\tau\|(I-T_N)\big(\Pi_\theta E^\theta(t_n)\Pi_\theta\overline{E^\theta}(t_n)\big)\|_{X_{2,\tau}^{s_0+1,b_0-1}}\\
&\lesssim\tau\theta^{\frac{s_1}2}\|\Pi_\theta E^\theta(t_n)\Pi_\theta\overline{E^\theta}(t_n)\|_{X_{2,\tau}^{s_2+1,b_0-1}}\lesssim\tau\theta^{\frac{s_1}2}\|E^\theta(t_n)\|^2_{X_{1,\tau}^{s_2+\frac12,b_0}}\lesssim \tau\theta^{\frac{s_1}2}.
\end{align*}
Combining the three component estimates yields the asserted spatial-local-error bound and completes the proof.
\end{proof}

\section{Global error analysis}\label{sectionglobal}

In this section we estimate the global error and prove our main result, Theorem~\ref{mainthm}. We derive simultaneous bounds for the temporal and spatial components of the global error and first establish a global error estimate; Theorem~\ref{mainthm} then follows from this estimate. By \eqref{minusest}, it suffices to control the two error sequences $e_{1,n}=E^\theta(t_n)-E_n$ and $e_{2,n}=u^\theta(t_n)-u_n$. Using \eqref{filtlie}, the global error can be decomposed as follows:
\begin{align}
e_{1,n}&=E^\theta(t_n)-E_n\nonumber\\
&=E^\theta(t_{n-1})-E_{n-1}-\frac i2\tau\fe^{i\tau\Delta}\Big(\phi_N^\tau(E^\theta(t_{n-1}), u^\theta(t_{n-1}))-\phi_N^\tau(E_{n-1},u_{n-1})\Big)-\mathcal{E}_{eloc}(t_{n-1})\nonumber\\
&=-\frac i2\tau\sum_{k=0}^{n-1}\fe^{i(n-k)\tau\Delta}\Big(\phi_N^\tau
(E^\theta(t_k),u^\theta(t_k))-\phi_N^\tau(E_k,u_k)\Big)-\sum_{k=0}^{n-1}\fe^{i(n-k-1)\tau\Delta}\mathcal{E}_{eloc}(t_k),\label{phi1}
\end{align}
where
$$
\phi_N^\tau(E,u)=\Pi_\theta T_N\Big(\frac{\fe^{-\frac i2\tau(\Pi_\theta u+\Pi_\theta\overline{u})}-1}{-i\tau/2}\Pi_\theta E\Big),
$$
and
\begin{align}
e_{2,n}&=u^\theta(t_n)-u_n\nonumber\\ &=u^\theta(t_{n-1})-u_{n-1}+i\tau|\nabla|\fe^{i\tau|\nabla|}
\Big(\psi_N^\tau(E^\theta(t_{n-1}))-\psi_N^\tau(E_{n-1})\Big)-\mathcal{E}_{uloc}(t_{n-1})\nonumber\\
&=i\tau|\nabla|\sum_{k=0}^{n-1}\fe^{i(n-k)\tau|\nabla|}\Big(\psi_N^\tau(E^\theta(t_k))-\psi_N^\tau(E_k)\Big)-\sum_{k=0}^{n-1}\fe^{i(n-k-1)\tau|\nabla|}\mathcal{E}_{uloc}(t_k),\label{psi1}
\end{align}
where
$$
\psi_N^\tau(E)=\Pi_\theta T_N(\Pi_\theta E\Pi_\theta\overline{E}).
$$

\begin{theorem}\label{theoglobal}
Let the parameters $s_0,~s_1,~b_0$ be as in \eqref{para} and let $(E^\theta,u^\theta)$ be the solution presented in Theorem~\ref{theoexist}. Then there exists $\tau_0>0$ such that for all $\tau\le \tau_0$ and all $n$ with $n\tau \le T$, the following estimate holds:
\begin{equation}\label{global}
\|e_{1,n}\|_{X_{1,\tau}^{s_0+\frac12,b_0}}+\|e_{2,n}\|_{X_{2,\tau}^{s_0,b_0}}\leq C_T\theta^{\frac{s_c}2},
\end{equation}
where $s_c=\min\{2,s_1\}$, and the constant $C_T$ depends on $T$ but is independent of $n$ and $\tau$.
\end{theorem}
\begin{proof}
Without loss of generality assume $s_1\leq2$, so $s_c=s_1$. Combining \eqref{discb-1} with Proposition~\ref{proptlocal} and Proposition~\ref{propslocal} yields
$$
\Big\|\sum_{k=0}^{n-1}\fe^{i(n-k-1)\tau\Delta}\mathcal{E}_{eloc}(t_k)\Big\|_{X_{1,\tau}^{s_0+\frac12,b_0}}+\Big\|\sum_{k=0}^{n-1}\fe^{i(n-k-1)\tau|\nabla|}\mathcal{E}_{uloc}(t_k)\Big\|_{X_{2,\tau}^{s_0,b_0}}\leq C_T\theta^{\frac{s_1}2}.
$$
Therefore, using \eqref{discb-1} together with \eqref{phi1} and \eqref{psi1}, we obtain for $e_{1,n}$ and $e_{2,n}$ supported in $[-T_1,T_1]$,
\begin{equation}\label{glost}
\begin{aligned}
\|e_{1,n}\|_{X_{1,\tau}^{s_0+\frac12,b_0}}&+\|e_{2,n}\|_{X_{2,\tau}^{s_0,b_0}}\lesssim T_1^{\varepsilon_0}\big\|\psi_N^\tau(E^\theta(t_n))-\psi_N^\tau(E_n)\big\|_{X_{2,\tau}^{s_0+1,-b_1}}\\
&+T_1^{\varepsilon_0}\big\|\phi_N^\tau(E^\theta(t_n),u^\theta(t_n))-\phi_N^\tau(E_n, u_n)\big\|_{X_{1,\tau}^{s_0+\frac12,-b_1}}+\theta^{\frac{s_1}2},
\end{aligned}
\end{equation}
where $\varepsilon_0=1-b_0-b_1$ and $T_1$ is to be determined later.

We next estimate the contribution from $\psi_N^\tau$. Utilizing \eqref{s-sp}, \eqref{discmult5} and \eqref{tnkeep}, and recalling $N^{-2}\lesssim\theta$ and $E_n=E^\theta(t_n)-e_{1,n}$, we arrive at
\begin{equation}\label{psiest}
\begin{aligned}
&\quad\big\|\psi_N^\tau(E^\theta(t_n))-\psi_N^\tau(E_n)\big\|_{X_{2,\tau}^{s_0+1,-b_1}}\\
&\lesssim\|\Pi_\theta E^\theta(t_n)\Pi_\theta e_{1,n}\|_{X_{2,\tau}^{s_0+1,-b_1}}+\|\Pi_\theta e_{1,n}\Pi_\theta e_{1,n}\|_{X_{2,\tau}^{s_0+1,-b_1}}\\
&\lesssim\|E^\theta(t_n)\|_{X_{1,\tau}^{s_2+\frac12,b_0}}\|e_{1,n}\|_{X_{1,\tau}^{s_0+\frac12,b_0}}+\theta^{-\frac{s_3}2}\|e_{1,n}\|^2_{X_{1,\tau}^{s_0+\frac12,b_0}},
\end{aligned}
\end{equation}
where $\varepsilon,s_3\in(0,s_1)$ will be chosen later. The term involving $\phi_N^\tau$ requires a more detailed treatment. Writing $F(u)=u+\overline{u}$ and splitting the expression as
\begin{equation}\label{phisp}
\begin{aligned}
&\quad\big\|\phi_N^\tau(E^\theta(t_n),u^\theta(t_n))-\phi_N^\tau(E_n,u_n)\big\|_{X_{1,\tau}^{s_0+\frac12,-b_1}}\\
&\lesssim\Big\|\Pi_\theta\Big(\frac{\fe^{-\frac i2\tau\Pi_\theta F(u^\theta(t_n))}-1}{-i\tau/2}\Pi_\theta E^\theta(t_n)-\frac{\fe^{-\frac i2\tau\Pi_\theta F(u_n)}-1}{-i\tau/2}\Pi_\theta E_n\Big)\Big\|_{X_{1,\tau}^{s_0+\frac12,-b_1}}\\
&\quad+\Big\|\Pi_\theta(I-T_N)\Big(\frac{\fe^{-\frac i2\tau\Pi_\theta F(u^\theta(t_n))}-1}{-i\tau/2}\Pi_\theta E^\theta(t_n)\Big)\Big\|_{X_{1,\tau}^{s_0+\frac12,-b_1}}\\
&\quad+\Big\|\Pi_\theta(I-T_N)\Big(\frac{\fe^{-\frac i2\tau\Pi_\theta F(u_n)}-1}{-i\tau/2}\Pi_\theta E_n\Big)\Big\|_{X_{1,\tau}^{s_0+\frac12,-b_1}}=\Phi_1+\Phi_2+\Phi_3.
\end{aligned}
\end{equation}
We estimate the resulting pieces separately.

Observing
\begin{equation}\label{phi1phi2}
\frac{\fe^{-\frac i2\tau\Pi_\theta F(u)}-1}{-i\tau/2}=\Pi_\theta F(u)-\tfrac i2\tau\Pi_\theta F(u)\Pi_\theta F(u)\varphi_2\bigl(-\tfrac i2\tau(\Pi_\theta F(u)\bigr),
\end{equation}
and applying \eqref{discmult3} and \eqref{discmult4} yields
\begin{equation}\label{phi11}
\begin{aligned}
\Phi_1&\lesssim\|E^\theta(t_n)\|_{X_{1,\tau}^{s_2+\frac12,b_0}}\|e_{2,n}\|_{X_{2,\tau}^{s_0,b_0}}+\|e_{1,n}\|_{X_{1,\tau}^{s_0+\frac12,b_0}}\|u^\theta(t_n)\|_{X_{2,\tau}^{s_2,b_0}}\\	
&\quad+\theta^{-\frac{s_3}2}\|e_{1,n}\|_{X_{1,\tau}^{s_0+\frac12,b_0}}\|e_{2,n}\|_{X_{2,\tau}^{s_0,b_0}}+\tau\|G(E^\theta(t_n),u^\theta(t_n),e_{1,n},e_{2,n})\|_{l_\tau^2 H^{s_0+\frac12}},
\end{aligned}
\end{equation}
where
$$
G(E,u,x,y)=(E-x)(u-y)^2\varphi_2(\tau u-\tau y)-Eu^2\varphi_2(\tau u).
$$
To bound the auxiliary function $G$, we first apply H\"older's inequality and \eqref{cbe} to obtain
\begin{align*}
&\quad\tau\|G(E^\theta(t_n),u^\theta(t_n),e_{1,n},e_{2,n})\|_{l_\tau^2 H^{s_0+\frac12}}\\
&\lesssim\tau(\|E^\theta(t_n)\|_{l_\tau^6H^{s_0+\frac12}}+\|e_{1,n}\|_{l_\tau^6H^{s_0+\frac12}})(\|\Pi_\theta u^\theta(t_n)\|_{l_\tau^6H^{s_0+1+\varepsilon}}+\theta^{-\frac12-\frac12\varepsilon}\|e_{2,n}\|_{l_\tau^6H^{s_0}})^2\\
&\quad\times\big(\varphi_2(\tau\|\Pi_\theta u^\theta(t_n)-e_{2,n}\|_{l_\tau^\infty H^{s_0+1+\varepsilon}})+\varphi_2(\tau\|\Pi_\theta u^\theta(t_n)\|_{{l_\tau^\infty H^{s_2+1}}})\big),
\end{align*}
where $\varepsilon\in(0,\frac16)$ can be taken arbitrarily small. Thus by applying \eqref{s-sp}, \eqref{disclinf}, \eqref{phik1}, \eqref{phik2} and \eqref{discl4}, we obtain
\begin{equation}\label{phi12}
\begin{aligned}
&\quad\tau\|G(E^\theta(t_n),u^\theta(t_n),e_{1,n},e_{2,n})\|_{l_\tau^2 H^{s_0+\frac12}}\\
&\lesssim(\|E^\theta(t_n)\|_{X_{1,\tau}^{s_2+\frac12,b_0}}+\|e_{1,n}\|_{X_{1,\tau}^{s_0+\frac12,b_0}})(\theta^{\frac{s_1}4}\|u^\theta(t_n)\|_{X_{2,\tau}^{s_2,b_0}}+\theta^{-\frac\varepsilon2}\|e_{2,n}\|_{X_{2,\tau}^{s_0,b_0}})^2\\
&\quad\times\big(\varphi_2(2\tau^\frac13\|\Pi_\theta u^\theta(t_n)\|_{X_{2,\tau}^{s_2,b_0}})+\varphi_2(2\tau^\frac13\|e_{2,n}\|_{X_{2,\tau}^{s_0,b_0}})+\varphi_2(\tau^\frac13\|\Pi_\theta u^\theta(t_n)\|_{X_{2,\tau}^{s_2,b_0}})\big)\\
&\lesssim\big(\theta^{\frac{s_1}2}+\|e_{1,n}\|_{X_{1,\tau}^{s_0+\frac12,b_0}}+\|e_{2,n}\|_{X_{2,\tau}^{s_0,b_0}}+\|e_{1,n}\|_{X_{1,\tau}^{s_0+\frac12,b_0}}\|e_{2,n}\|_{X_{2,\tau}^{s_0,b_0}}\\
&\quad+\theta^{-\varepsilon}\|e_{2,n}\|^2_{X_{2,\tau}^{s_0,b_0}}+\theta^{-\varepsilon}\|e_{1,n}\|_{X_{1,\tau}^{s_0+\frac12,b_0}}\|e_{2,n}\|^2_{X_{2,\tau}^{s_0,b_0}}\big)\big(1+\varphi_2(2\tau^\frac13\|e_{2,n}\|_{X_{2,\tau}^{s_0,b_0}})\big),
\end{aligned}
\end{equation}
where $s=\min\{1,s_1\}$. The factor $\theta^{-s_3/2}$ in \eqref{phi11} follows from the $\Pi_\theta$ support of the error and \eqref{s-sp}:
$$
\|e_{2,n}\|_{X_{2,\tau}^{s_0+s_3,b_1}}\lesssim\theta^{-s_3/2}\|e_{2,n}\|_{X_{2,\tau}^{s_0,b_0}}.
$$
Collecting \eqref{phi11}, \eqref{phi12} together with \eqref{discbounde} gives
\begin{equation}\label{phiest1}
\begin{aligned}
\Phi_1&\lesssim\big(\theta^{\frac{s_1}2}+\|e_{1,n}\|_{X_{1,\tau}^{s_0+\frac12,b_0}}+\|e_{2,n}\|_{X_{2,\tau}^{s_0,b_0}}+\theta^{-\frac{s_3}2}\|e_{1,n}\|_{X_{1,\tau}^{s_0+\frac12,b_0}}\|e_{2,n}\|_{X_{2,\tau}^{s_0,b_0}}\\
&\quad+\theta^{-\varepsilon}\|e_{2,n}\|^2_{X_{2,\tau}^{s_0,b_0}}+\theta^{-\varepsilon}\|e_{1,n}\|_{X_{1,\tau}^{s_0+\frac12,b_0}}\|e_{2,n}\|^2_{X_{2,\tau}^{s_0,b_0}}\big)\big(1+\varphi_2(2\tau^\frac13\|e_{2,n}\|_{X_{2,\tau}^{s_0,b_0}})\big).
\end{aligned}
\end{equation}

We then estimate $\Phi_2$. Using \eqref{discont}, \eqref{discontbourg} and \eqref{phi1phi2}, we obtain
\begin{equation}\label{phi21}
\begin{aligned}
\Phi_2&\lesssim\big\|(I-T_N)\big(\Pi_\theta F(u^\theta(t_n))\Pi_\theta E^\theta(t_n)\big)\big\|_{X_{1,\tau}^{s_0+\frac12,-b_1}}\\
&\quad+\tau\Big\|(I-T_N)\Big(\big(\Pi_\theta F(u^\theta(t_n))\big)^2\varphi_2\big(\tfrac\tau2\Pi_\theta F(u^\theta(t_n))\big)\Pi_\theta E^\theta(t_n)\Big)\Big\|_{l_\tau^2H^{s_0+\frac12}}\\
&\lesssim\theta^{\frac{s_1}2+\frac14}\|\Pi_\theta F(u^\theta(t_n))\Pi_\theta E^\theta(t_n)\|_{X_{1,\tau}^{s_2+1,-b_1}}\\
&\quad+\tau\theta^{\frac14+\frac{s_1}2}\big\|\big(\Pi_\theta F(u^\theta(t_n))\big)^2\varphi_2\big(\tfrac\tau2\Pi_\theta F(u^\theta(t_n))\big)\Pi_\theta E^\theta(t_n)\big\|_{l_\tau^2H^{s_2+1}}.
\end{aligned}
\end{equation}
The same arguments used in Proposition~\ref{propslocal} to estimate $\mathcal{E}_1'$ and $\mathcal{E}_2'$ yield the required control for $\Phi_2$:
\begin{equation}\label{phiest2}
\Phi_2\lesssim\theta^{\frac{s_1}2}.
\end{equation}
Similarly, note that $E_n=E^\theta(t_n)-e_{1,n}$ and $u_n=u^\theta(t_n)-e_{2,n}$, we obtain an estimate for $\Phi_3$ that
\begin{equation}\label{phiest3}
\begin{aligned} \Phi_3&\lesssim\big(\theta^{\frac{s_1}2}+\|e_{1,n}\|_{X_{1,\tau}^{s_0+\frac12,b_0}}+\theta^{-\frac{s_3}2}\|e_{1,n}\|_{X_{1,\tau}^{s_0+\frac12,b_0}}\|e_{2,n}\|_{X_{2,\tau}^{s_0,b_0}}\\
&\quad\quad\quad+\theta^{-s_3}\|e_{1,n}\|_{X_{1,\tau}^{s_0+\frac12,b_0}}\|e_{2,n}\|^2_{X_{2,\tau}^{s_0,b_0}}\big)\times\big(1+\varphi_2(2\tau^\frac13\|e_{2,n}\|_{X_{2,\tau}^{s_0,b_0}})\big).
\end{aligned}
\end{equation}

Combining \eqref{glost}, \eqref{psiest}, \eqref{phiest1}, \eqref{phiest2} and \eqref{phiest3} leads to
\begin{equation}\label{glot1}
\begin{aligned}			
&\quad\|e_{1,n}\|_{X_{1,\tau}^{s_0+\frac12,b_0}}+\|e_{2,n}\|_{X_{2,\tau}^{s_0,b_0}}\\
&\leq C_T\theta^{\frac{s_1}2}+C_TT_1^{\varepsilon_0}\big(\|e_{1,n}\|_{X_{1,\tau}^{s_0+\frac12,b_0}}+\|e_{2,n}\|_{X_{2,\tau}^{s_0,b_0}}+\theta^{-\frac{s_3}2}\|e_{1,n}\|^2_{X_{1,\tau}^{s_0+\frac12,b_0}}\\
&\quad+\theta^{-\frac{s_3}2}\|e_{2,n}\|^2_{X_{2,\tau}^{s_0,b_0}}+\theta^{-\frac{s_3}2}\|e_{1,n}\|_{X_{1,\tau}^{s_0+\frac12,b_0}}\|e_{2,n}\|_{X_{2,\tau}^{s_0,b_0}}\\
&\quad+\theta^{-s_3}\|e_{1,n}\|_{X_{1,\tau}^{s_0+\frac12,b_0}}\|e_{2,n}\|^2_{X_{2,\tau}^{s_0,b_0}}\big)\big(1+\varphi_2(2\tau^\frac13\|e_{2,n}\|_{X_{2,\tau}^{s_0,b_0}})\big).
\end{aligned}
\end{equation}

We’ll now select the parameters. Choose $T_1$ such that $C_T T_1^{\varepsilon_0}\le\frac13$, pick $s_3\in(1-2b_1,s_1)$, and choose $\varepsilon<\min\{\frac16,\frac{s_3}2\}$. With these choices, using the a priori assumption
$$
\|e_{1,n}\|_{X_{1,\tau}^{s_0+\frac12,b_0}}+\|e_{2,n}\|_{X_{2,\tau}^{s_0,b_0}}\le C_T\theta^{\frac{s_3}2},
$$
the estimate \eqref{glot1} yields
$$
\|e_{1,n}\|_{X_{1,\tau}^{s_0+\frac12,b_0}}+\|e_{2,n}\|_{X_{2,\tau}^{s_0,b_0}}\le C_T\theta^{\frac{s_1}2},
$$
which is the desired error bound on $[0,T_1]$. Since $T_1$ depends only on $T$, iterating this argument over successive subintervals covering $[0,T]$ completes the proof.
\end{proof}

Now we are ready to prove our main result, Theorem~\ref{mainthm}.

\noindent\emph{\textbf{Proof of Theorem~\ref{mainthm}.}} By \eqref{zworkspace}, \eqref{linf} and \eqref{minusest}, we have
\begin{align*}
&\quad\|E(t_n)-E^\theta(t_n)\|_{H^{s_0+\frac12}}+\|z(t_n)-z^\theta(t_n)\|_{H^{s_0}}+\|z_t(t_n)-z_t^\theta(t_n)\|_{H^{s_0-1}}\\
&\leq\|E-E^\theta\|_{L^\infty H^{s_0+\frac12}}+\|z-z^\theta\|_{L^\infty H^{s_0}}+\|z_t-z_t^\theta\|_{L^\infty H^{s_0-1}}\\
&\leq\|E-E^\theta\|_{X_1^{s_0+\frac12,b_0}}+\|(z-z^\theta,z_t-z_t^\theta)\|_{\mathcal{X}_2^{s_0,b_0}}\leq C_T\theta^{\frac{s_1}2}.
\end{align*}
Therefore, by the triangle inequality and the definitions of $e_{1,n}$ and $e_{2,n}$, it suffices to prove
$$
\|e_{1,n}\|_{H^{s_0+\frac12}}+\tfrac12\|e_{2,n}+\overline{e_{2,n}}\|_{H^{s_0}}
+\tfrac12\|e_{2,n}-\overline{e_{2,n}}\|_{H^{s_0}}\leq C_T\theta^{\frac{s_1}2}.
$$
Applying the triangle inequality, one gets
$$
\tfrac12\|e_{2,n}+\overline{e_{2,n}}\|_{H^{s_0}}+\tfrac12\|e_{2,n}
-\overline{e_{2,n}}\|_{H^{s_0}}\leq\|e_{2,n}\|_{l_\tau^\infty H^{s_0}}+\|\overline{e_{2,n}}\|_{l_\tau^\infty H^{s_0}}=2\|e_{2,n}\|_{l_\tau^\infty H^{s_0}}.
$$
Invoking \eqref{disclinf} together with Theorem~\ref{theoglobal} yields
$$
\|e_{1,n}\|_{H^{s_0+\frac12}}+\tfrac12\|e_{2,n}+\overline{e_{2,n}}\|_{H^{s_0}}
+\tfrac12\|e_{2,n}-\overline{e_{2,n}}\|_{H^{s_0}} \leq\|e_{1,n}\|_{X_{1,\tau}^{s_0+\frac12,b_0}}+2\|e_{2,n}\|_{X_{2,\tau}^{s_0,b_0}}\leq C_T\theta^{\frac{s_1}2},
$$
which establishes Theorem~\ref{mainthm}.

\section{Numerical results}\label{sectionnumerexp}

In this section we present numerical validations of Theorem~\ref{mainthm} and investigate the conservation properties of the filtered Lie splitting scheme \eqref{filtlie}. We generate random initial data with varying regularity and measure the numerical error in the norm $H^{s_0+\frac12}\times H^{s_0}\times H^{s_0-1}$, where $s_0=\max\{0,\frac d2-1\}$. For each convergence experiment, we generate one initial data as follows:
\begin{align*}
E_0(x)&=\sum_{k_1,\ldots,k_d=-N_0/2}^{N_0/2-1}\langle k\rangle^{-s_2-\frac d2-\frac12-\varepsilon}\widetilde{f}_k\fe^{i\langle k,x\rangle},\\
z_0(x)&=\re\sum_{k_1,\ldots,k_d=-N_0/2}^{N_0/2-1}\langle k\rangle^{-s_2-\frac d2-\varepsilon}\widetilde{g}_k\fe^{i\langle k,x\rangle},\quad
z_1(x)=\re\sum_{k_1,\ldots,k_d=-N_0/2}^{N_0/2-1}\langle k\rangle^{-s_2-\frac d2+1-\varepsilon}\widetilde{h}_k\fe^{i\langle k,x\rangle},
\end{align*}
where $s_2=s_0+s_1$ with $s_1>0$, $k=(k_1,\ldots,k_d)$. Here $N_0$ denotes the number of grid points in each direction, and $\widetilde{f}_k,~\widetilde{g}_k,~\widetilde{h}_k$ are independent uniformly distributed random variables in $[-1,1]$ (with $\widetilde{h}_0=0$). For any $\varepsilon>0$, the $H^{s_2+\frac12}\times H^{s_2}\times H^{s_2-1}$ norm of the Fourier series is uniform with respect to $N_0$, whereas in the experiments we take the borderline case $\varepsilon=0$. The initial data are then normalized at $N_0$ so that
$$
\|E_0\|_{H^{s_2+\frac12}}=\|z_0\|_{H^{s_2}}=\|z_1\|_{H^{s_2-1}}=1.
$$
For each coarser grid, we use Fourier projections of the same normalized initial data, without regenerating the random coefficients or renormalizing the projected data. Computational timings were measured on a cluster equipped with an Intel Xeon Gold 6226R CPU (for $d=1$) and an NVIDIA A100 GPU ($d=2,3$). 

\subsection{Accuracy test}

In this part, we assess the numerical error of the method \eqref{filtlie} at the final time $T=1$. For $d=1$ and $d=2$, we test the values $s_1=s_2=0.25, 0.5, 1, 2$; for $d=3$, we consider $s_1=0.25, 0.5, 1$, corresponding to $s_2=0.75, 1, 1.5$, respectively. All convergence plots compare the filtered Lie method against a reference solution computed using an extremely small time step.

\begin{figure}[h!]
\begin{center}
\subfigure{\includegraphics[height=0.33\textwidth, width=0.48\textwidth]{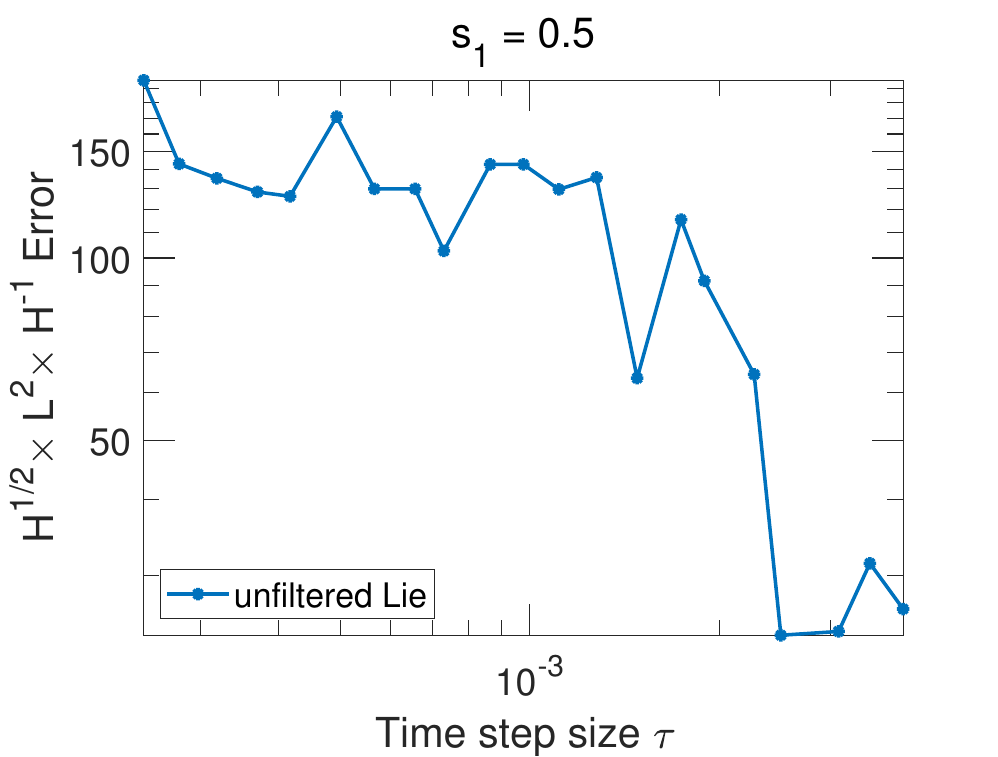}}
\subfigure{\includegraphics[height=0.33\textwidth, width=0.48\textwidth]{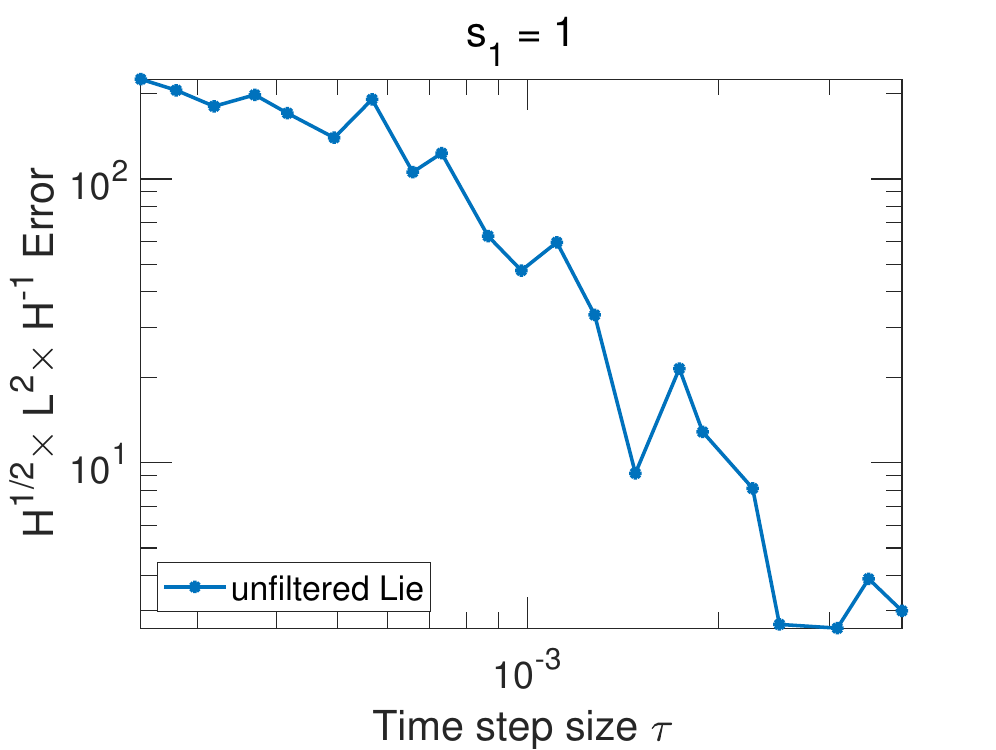}}
\subfigure{\includegraphics[height=0.33\textwidth, width=0.48\textwidth]{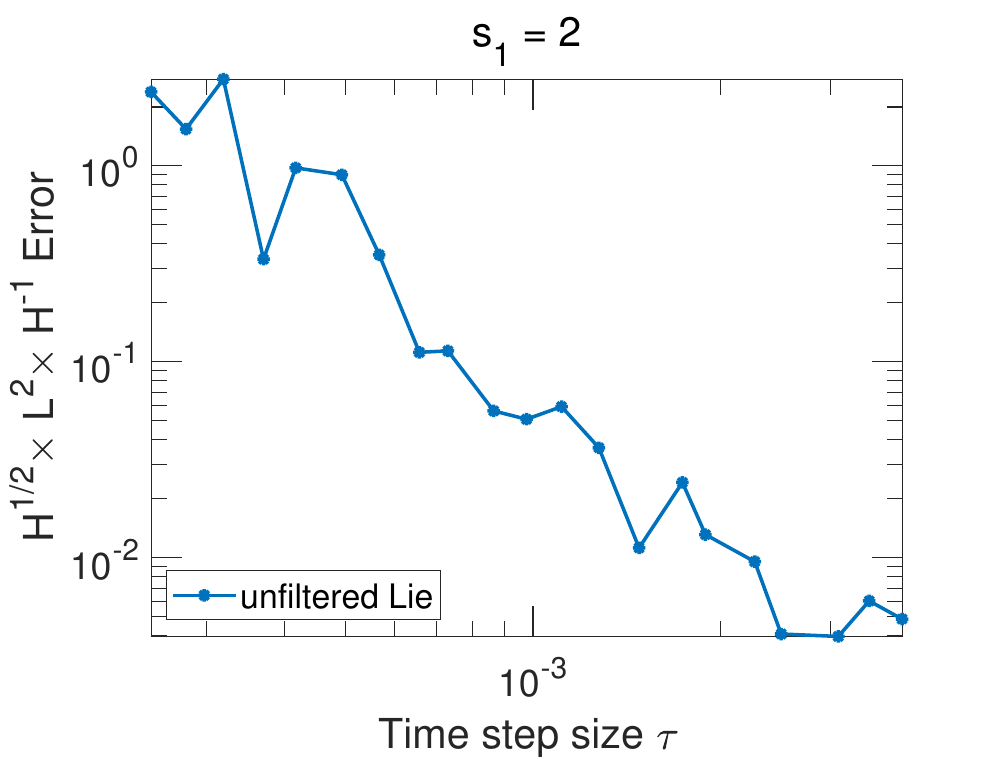}}
\subfigure{\includegraphics[height=0.33\textwidth, width=0.48\textwidth]{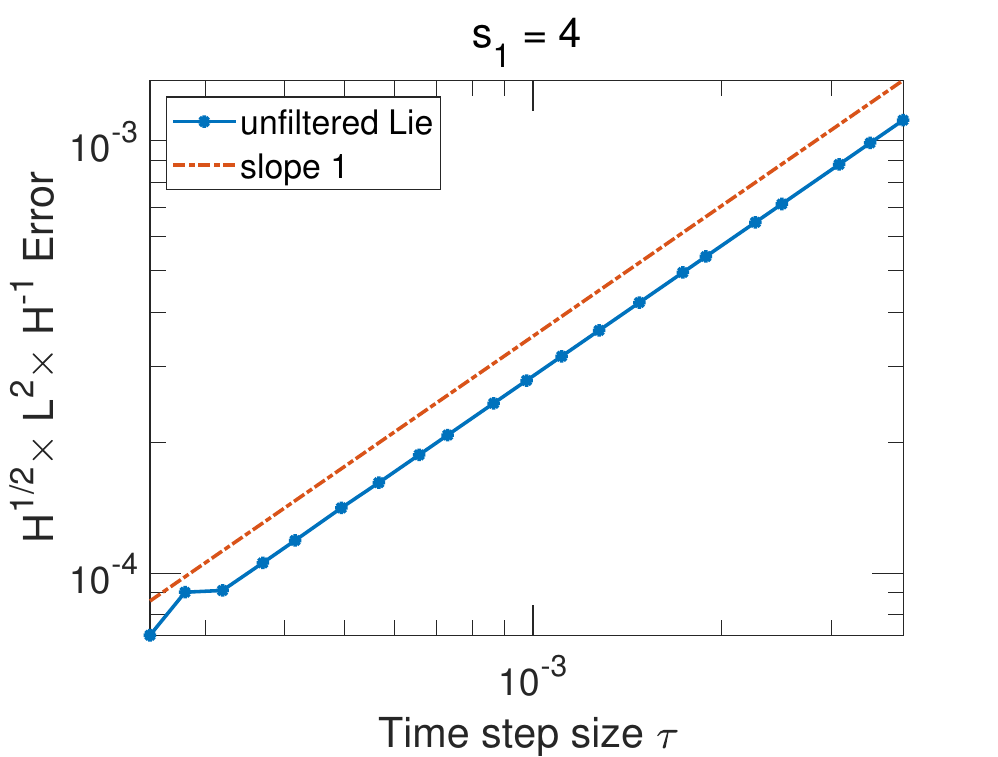}}
\end{center}
\vspace{-2mm}
\caption{The $H^{1/2}\times L^2\times H^{-1}$ error of the unfiltered Lie splitting scheme for the one-dimensional ZS with rough initial data $(E_0, z_0, z_1)\in H^{s_1+1/2}\times H^{s_1}\times H^{s_1-1}$.}\label{fig:error0}
\end{figure}

\begin{figure}
\begin{center}
\subfigure{\includegraphics[height=0.33\textwidth, width=0.48\textwidth]{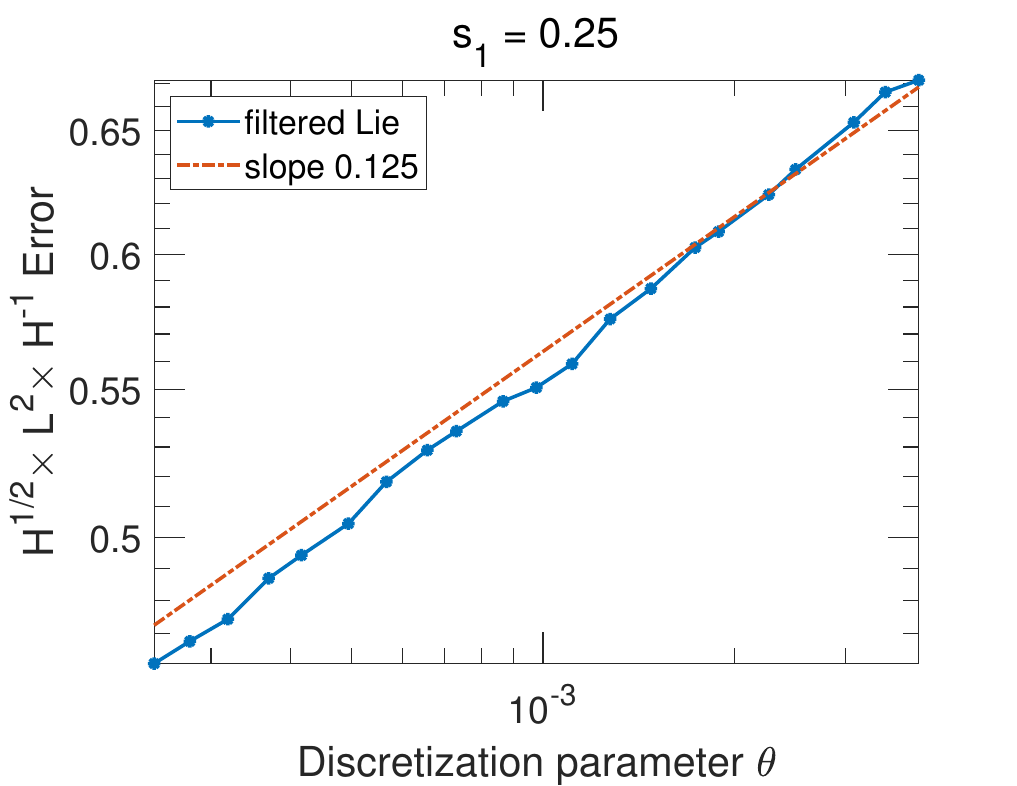}}
\subfigure{\includegraphics[height=0.33\textwidth, width=0.48\textwidth]{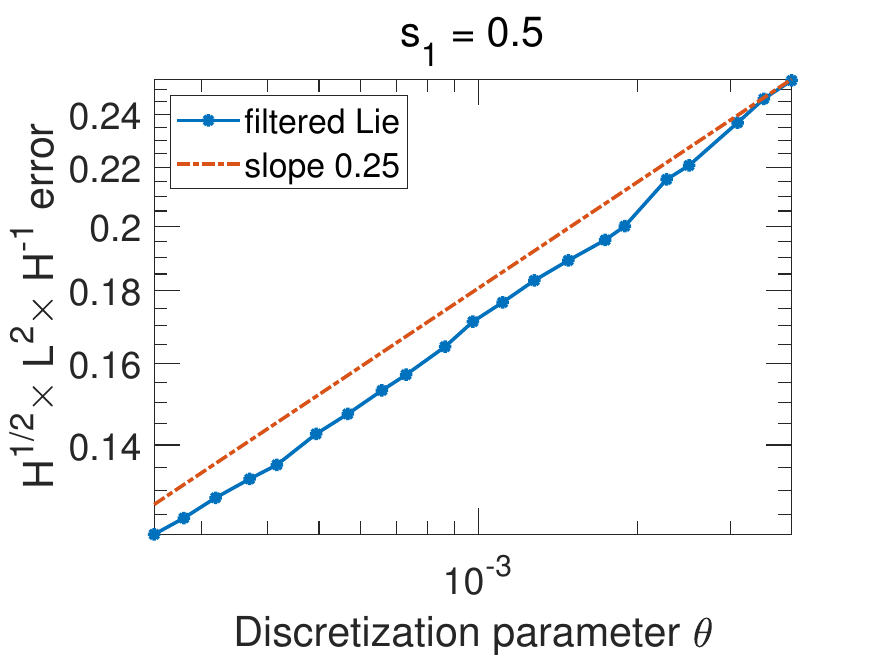}}
\subfigure{\includegraphics[height=0.33\textwidth, width=0.48\textwidth]{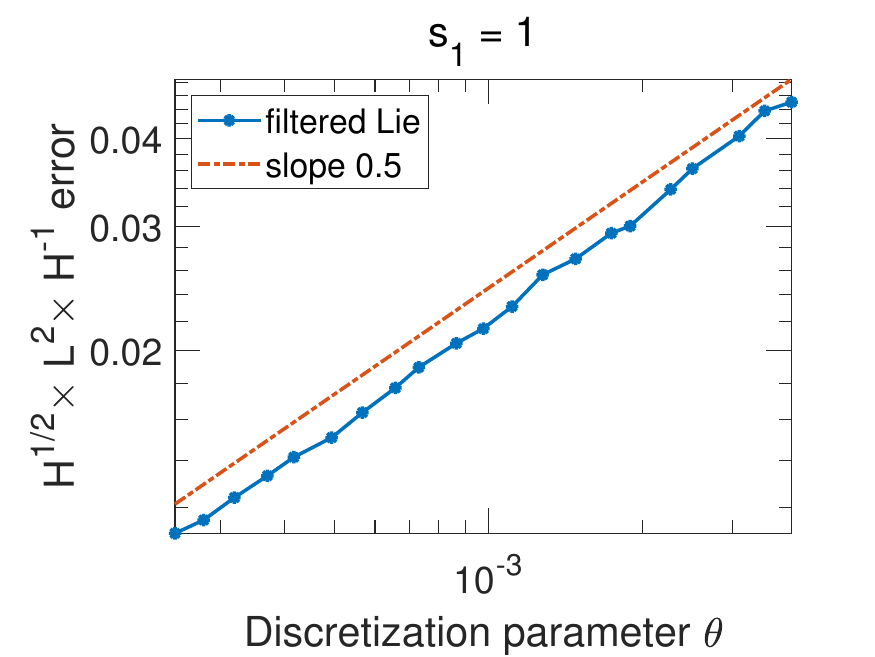}}
\subfigure{\includegraphics[height=0.33\textwidth, width=0.48\textwidth]{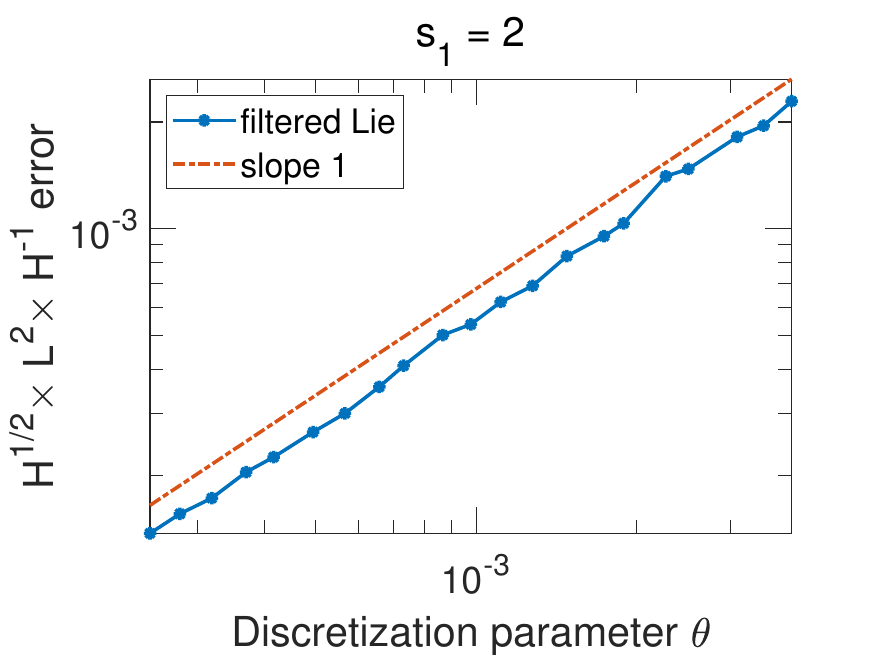}}
\end{center}
\vspace{-2mm}
\caption{The $H^{1/2}\times L^2\times H^{-1}$ error of the fully discretized filtered Lie splitting scheme for the one-dimensional ZS with rough initial data $(E_0, z_0, z_1)\in H^{s_1+1/2}\times H^{s_1}\times H^{s_1-1}$.}\label{fig:error1}
\end{figure}

\begin{figure}[h!]
\begin{center}
\subfigure{\includegraphics[height=0.33\textwidth, width=0.48\textwidth]{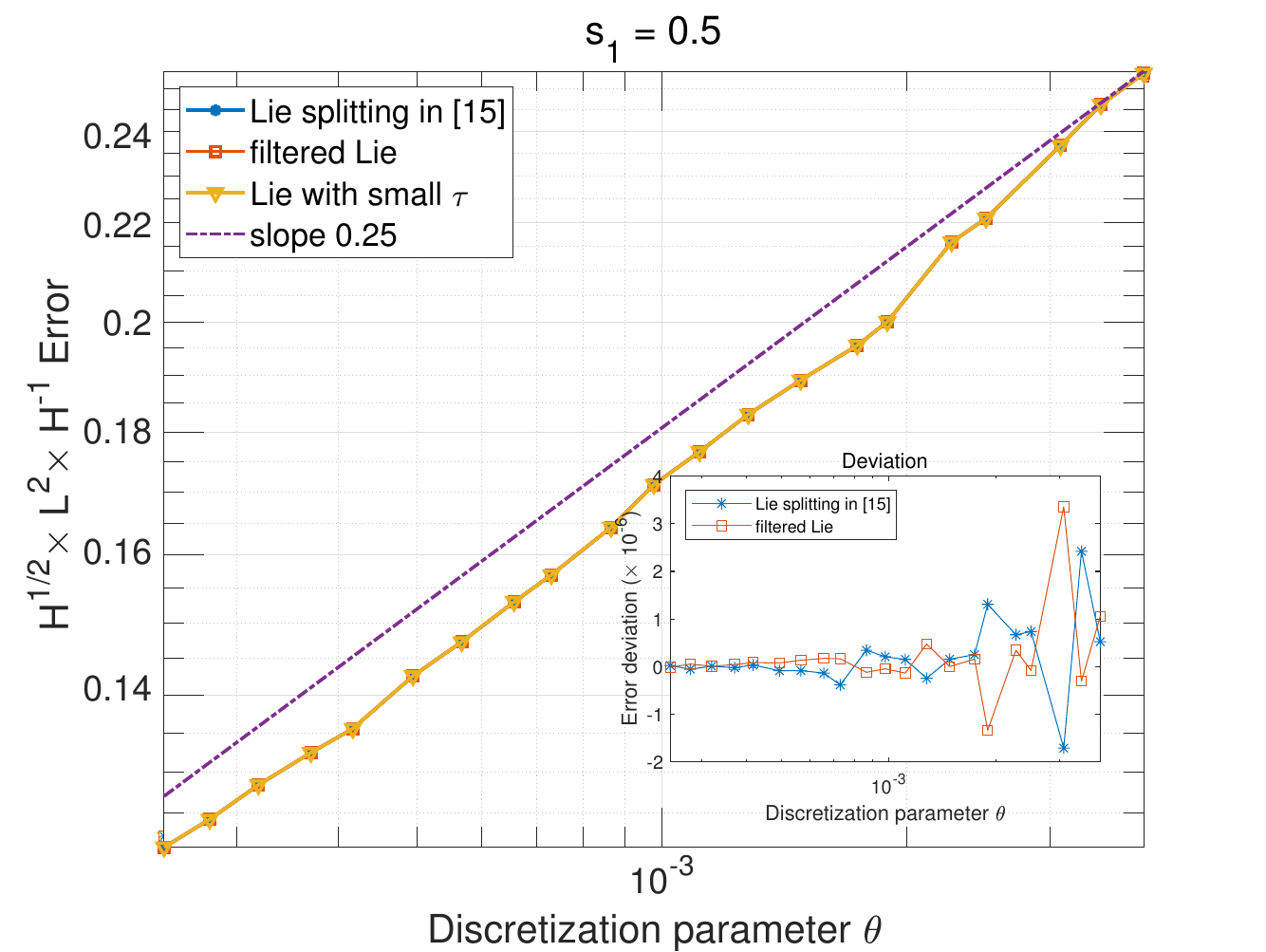}}
\subfigure{\includegraphics[height=0.33\textwidth, width=0.48\textwidth]{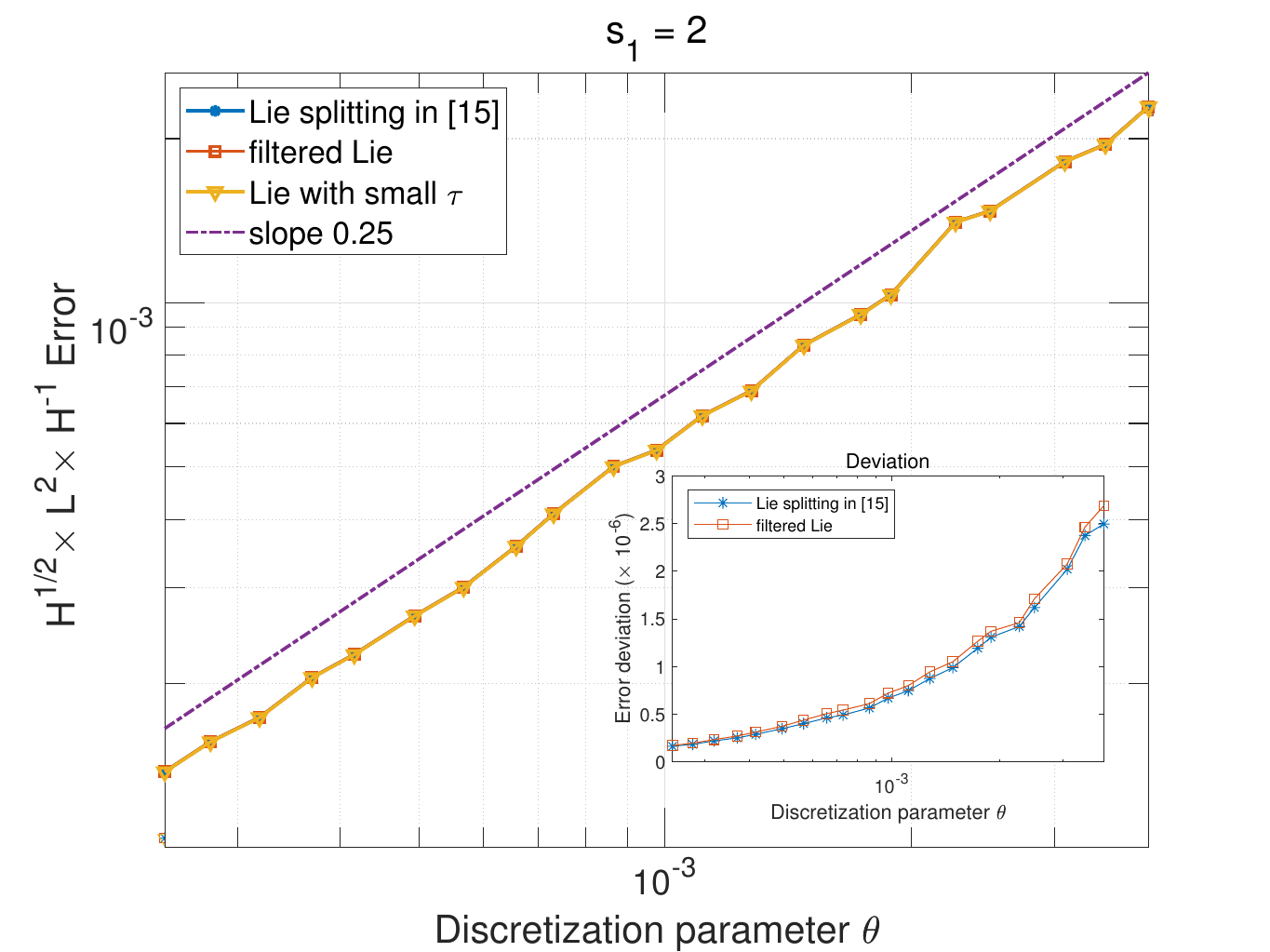}}
\end{center}
\vspace{-2mm}
\caption{The $H^{1/2}\times L^2\times H^{-1}$ error of the fully discretized filtered Lie splitting scheme compared with the same method with smaller temporal step size and the splitting in \cite{Gaucknm} for the one-dimensional ZS with initial data $(E_0, z_0, z_1)\in H^{s_1+1/2}\times H^{s_1}\times H^{s_1-1}$.}\label{fig:errorcompare}
\end{figure}

In the one-dimensional case, we compute the reference solution using the fully discretized filtered Lie splitting method with $N_0=2^{14}$ and $\tau=2^{-26}$. Figure~\ref{fig:error0} shows the numerical error of the unfiltered Lie splitting method with $N^2\tau=2^{16}$, which substantially exceeds the bound in \eqref{CFL}. For $s_1=0.5,1,2$, the results highlight the importance of filtering or imposing an appropriate CFL restriction. Indeed, the filter $\Pi_\tau$ is inactive on the collocation space when \eqref{CFL} holds. Nevertheless, first-order convergence is observed for the smoother case $s_1=4$. Figure~\ref{fig:error1} illustrates the numerical error of the filtered splitting method with $N=2\tau^{-1/2}$ and the time step varies over $\theta=\tau\in[2^{-12},2^{-8}]$. The results confirm that the filtered splitting method converges at $\mathcal{O}(\theta^{s_1/2})$ in the space $H^{1/2}\times L^2\times H^{-1}$ for solutions in $H^{s_1+1/2}\times H^{s_1}\times H^{s_1-1}$ for $s_1 = 0.25, 0.5, 1, 2$, in full agreement with Theorem~\ref{mainthm}.

We further compare our method \eqref{filtlie} with the Lie splitting method proposed in \cite{Gaucknm}, as shown in Figure~\ref{fig:errorcompare}. Both methods are applied with $\theta=\tau=4/N^2$.  To investigate the spatial discretization error, we also apply \eqref{filtlie} with a fixed time step $\tau=2^{-16}$ and varying spatial grid numbers $N\in[2^5,2^7]$. For both $s_1=0.5$ and $s_1=2$, the error curves of these methods almost coincide. The insets illustrate their small deviations, of order $\mathcal{O}(10^{-6})$, with the error of \eqref{filtlie} computed using a smaller time step taken as the reference (zero deviation).

\begin{figure}
\begin{center}
\subfigure{\includegraphics[height=0.33\textwidth,width=0.48\textwidth]{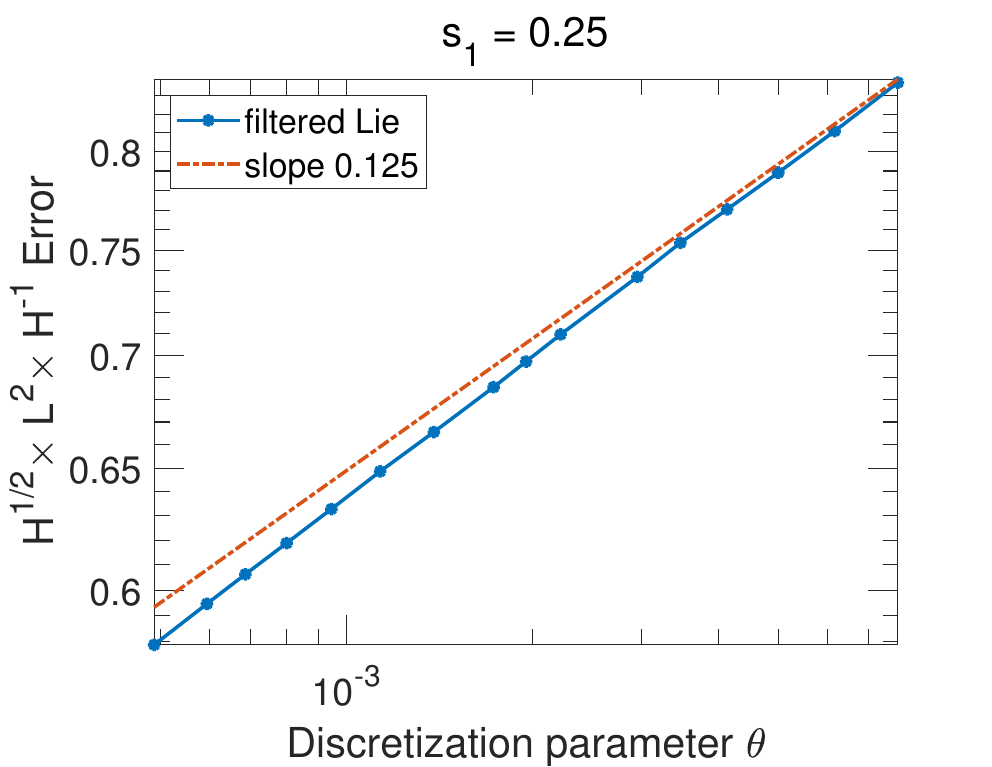}}
\subfigure{\includegraphics[height=0.33\textwidth,width=0.48\textwidth]{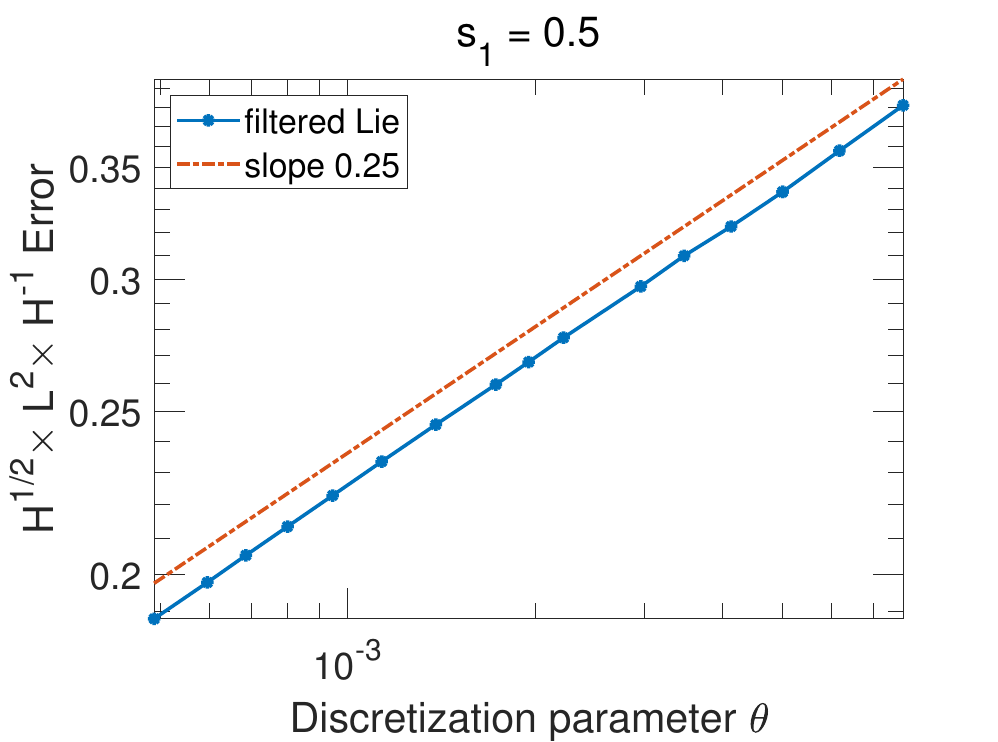}}
\subfigure{\includegraphics[height=0.33\textwidth,width=0.48\textwidth]{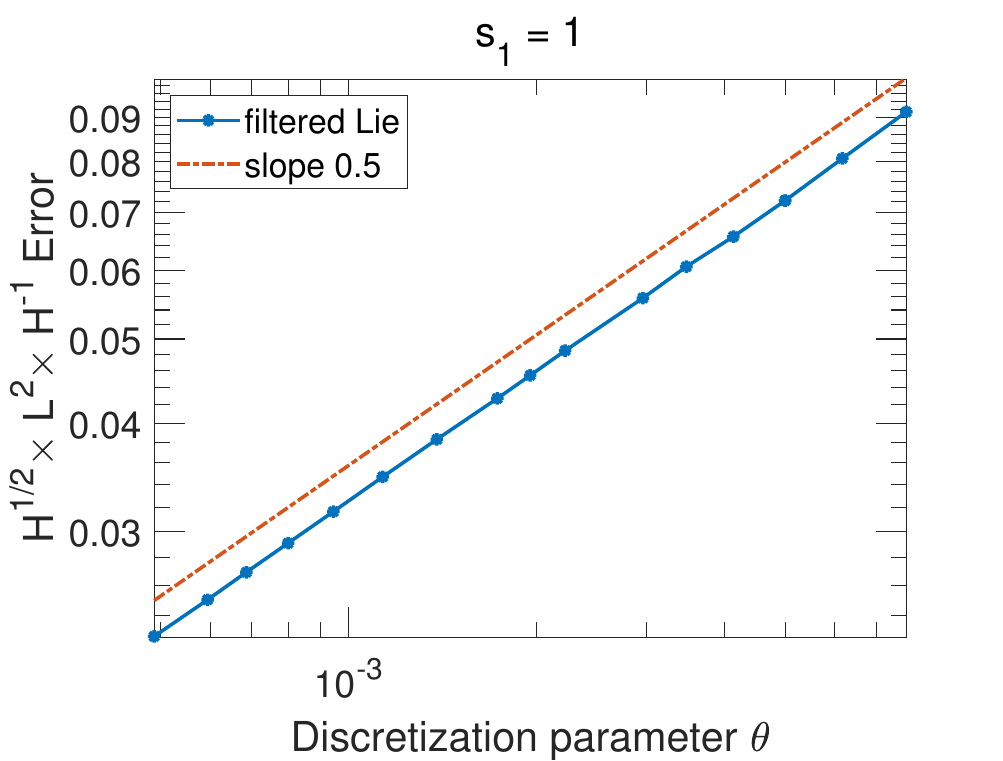}}
\subfigure{\includegraphics[height=0.33\textwidth,width=0.48\textwidth]{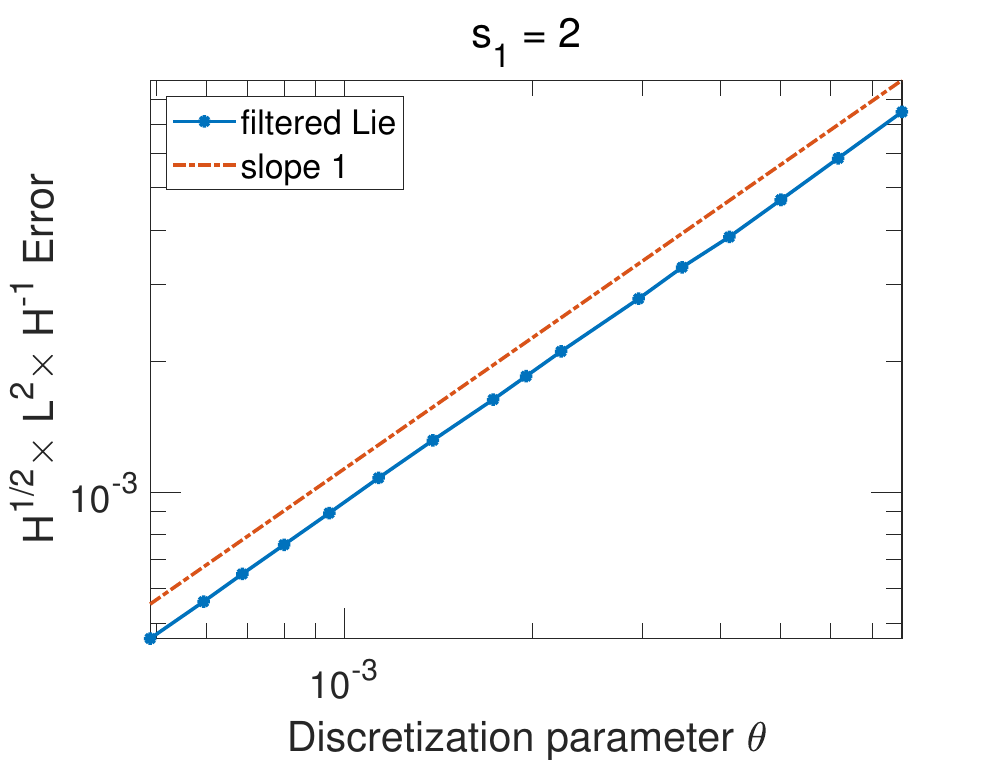}}
\end{center}
\vspace{-2mm}\caption{The $H^{1/2}\times L^2\times H^{-1}$ error of the filtered Lie splitting scheme for the two-dimensional ZS with rough initial data $(E_0, z_0, z_1)\in H^{s_1+1/2}\times H^{s_1}\times H^{s_1-1}$.}\label{fig:error2}
\end{figure}

In the case $d=2$, the reference solution is obtained using the filtered Lie splitting method with $N=2^{11}$ and $\tau=2^{-21}$. For $s_1=0.25$, we increased the resolution to $N_0=2^{13}$ and reduced the time step to $\tau=2^{-25}$ to ensure sufficient accuracy. The test runs employ the relation $N=\sqrt{2}\tau^{-1/2}$ with $\theta=\tau\in[2^{-11},2^{-7}]$. Figure~\ref{fig:error2} confirms the convergence rate $\mathcal{O}(\theta^{s_1/2})$ for $s_1=0.25,0.5,1,2$, in excellent agreement with our theoretical prediction.

\begin{figure}
\begin{center}
\subfigure{\includegraphics[height=0.33\textwidth,width=0.48\textwidth]{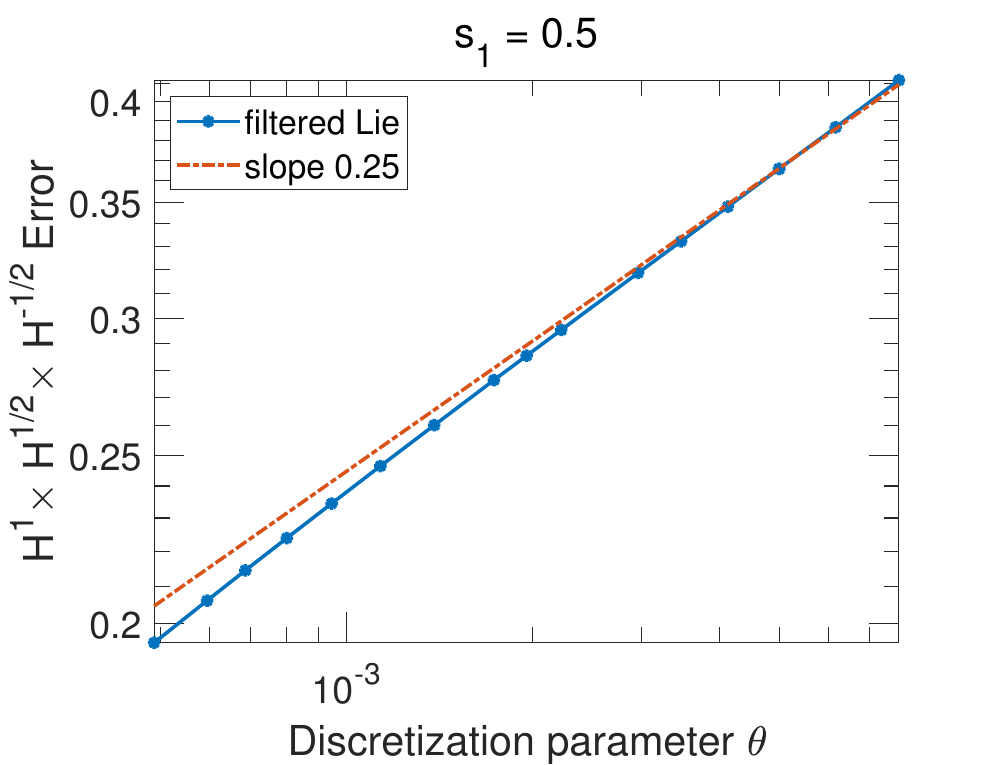}}
\subfigure{\includegraphics[height=0.33\textwidth,width=0.48\textwidth]{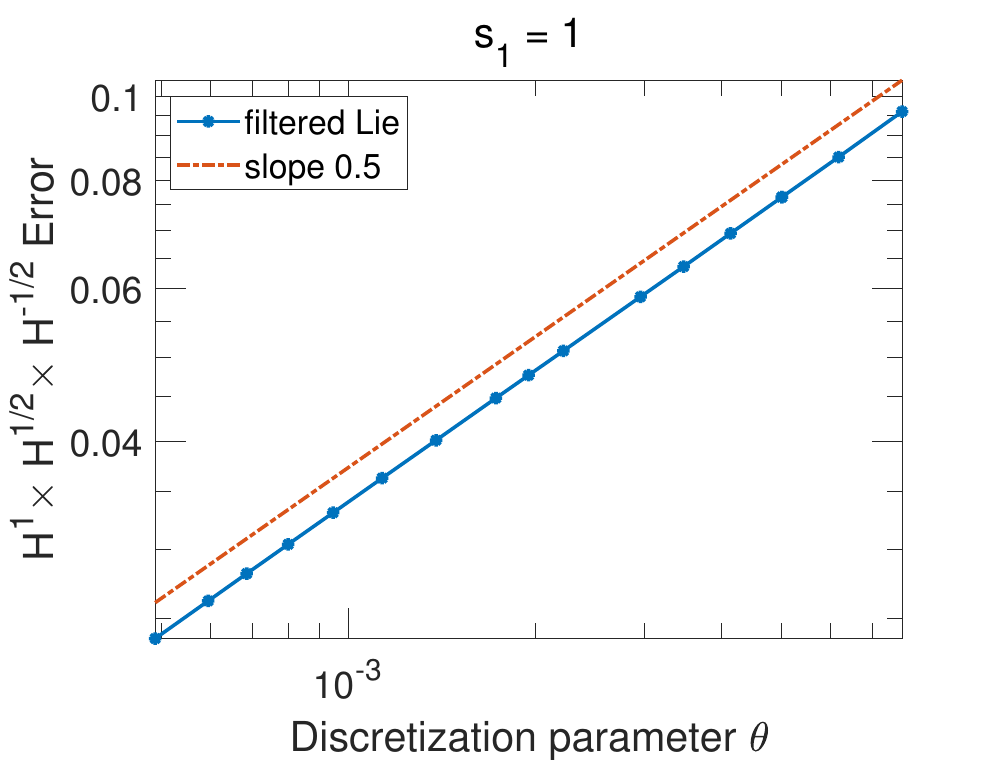}}
\end{center}
\vspace{-2mm}
\caption{The $H^1\times H^{1/2}\times H^{-1/2}$ error of the filtered Lie splitting scheme for the three-dimensional ZS with rough initial data $(E_0, z_0, z_1)\in H^{s_1+1}\times H^{s_1+1/2}\times H^{s_1-1/2}$.}\label{fig:error3}
\end{figure}

For the case $d=3$, reference computations were performed using the filtered Lie splitting method with $N_0=2^9$ and $\tau=2^{-17}$. The test runs employed $N=\sqrt{2}\tau^{-1/2}$ with $\theta=\tau\in[2^{-11},2^{-7}]$. Figure~\ref{fig:error3} illustrates the numerical errors for $s_1=0.5$ and $s_1=1$. The observed convergence rates align well with the theoretically predicted order $\mathcal{O}(\theta^{s_1/2})$, thereby confirming our analytical results.

\begin{figure}
\begin{center}
\subfigure{\includegraphics[height=0.33\textwidth,width=0.48\textwidth]{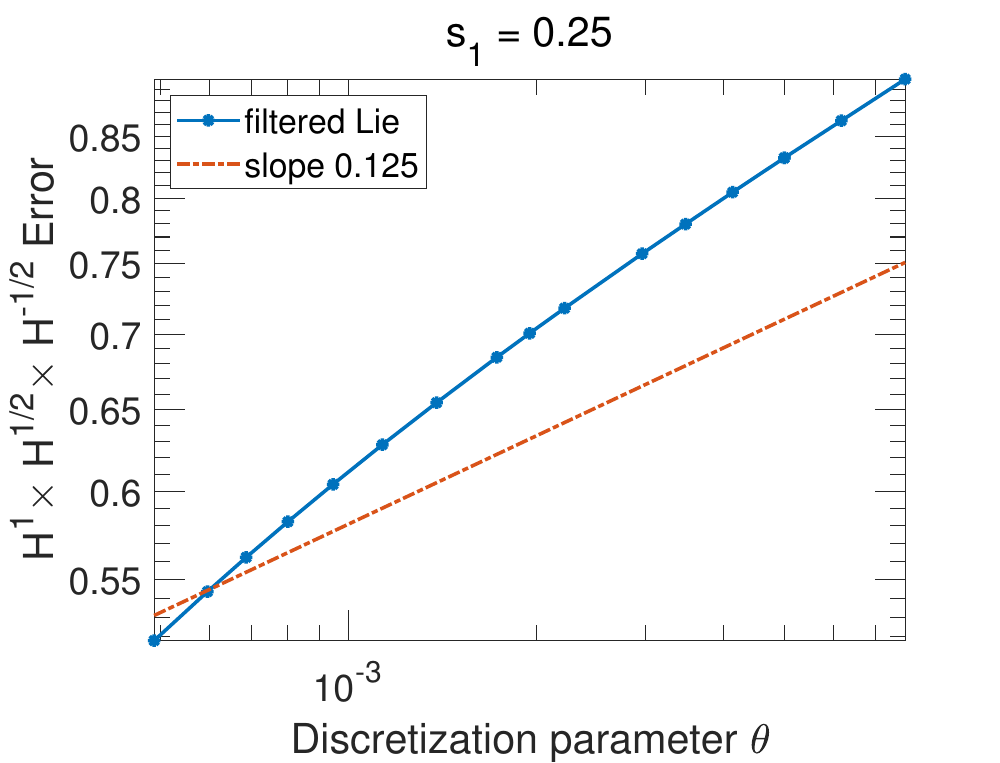}}
\subfigure{\includegraphics[height=0.33\textwidth,width=0.48\textwidth]{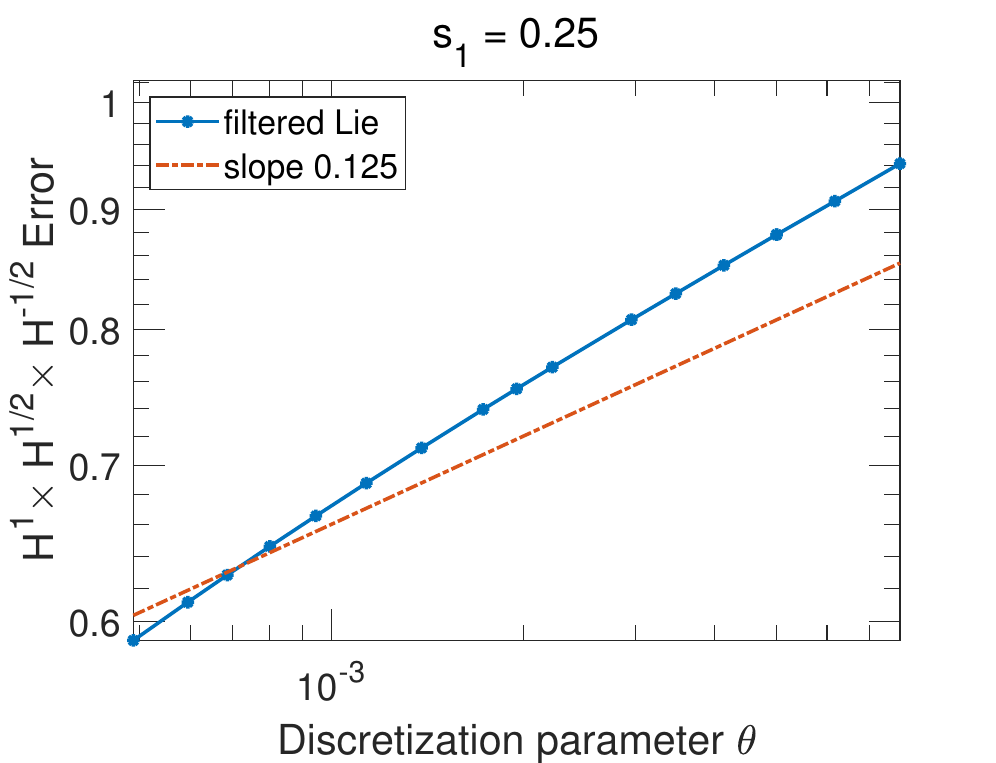}}
\end{center}
\vspace{-2mm}
\caption{The $H^1\times H^{1/2}\times H^{-1/2}$ error of the filtered Lie splitting scheme for the three-dimensional ZS with rough initial data $(E_0, z_0, z_1)\in H^{s_1+1}\times H^{s_1+1/2}\times H^{s_1-1/2}$ with $s_1=0.25$. Left: numerical errors computed against a reference solution with $N=2^8$ and $\tau=2^{-15}$. Right: numerical errors computed against a reference solution with $N=2^9$ and $\tau=2^{-17}$.}\label{fig:error4}
\end{figure}

However, for $s_1=0.25$, the numerical results are less conclusive. The slow convergence rate predicted by the theory makes its numerical verification particularly demanding, as the reference solution must be sufficiently accurate for its error to be negligible relative to the errors being measured. Computing such a reference solution at substantially finer spatial and temporal resolutions exceeds our available hardware resources. We therefore examine the sensitivity of the measured errors to the choice of reference solution, as shown in Figure~\ref{fig:error4}, using the same initial data with $N_0=2^9$ throughout. The results show better agreement with the predicted rate when the more finely resolved reference solution is used. The difference between the two reference solutions, measured in the same norm as the plotted errors, is approximately $0.28$, substantially smaller than the errors of the test solutions relative to their respective references. Although this comparison does not establish that the reference errors are negligible, it helps assess their influence on the observed convergence behavior. Similar difficulties and numerical comparisons are reported in \cite{Jisiam,Jiima}.

\subsection{Conservation properties}

In this part, we investigate the mass \eqref{massc} and energy \eqref{energyc} conservation properties of our numerical method \eqref{filtlie} for $s_2=0.5$, which corresponds to the minimal regularity required for the Hamiltonian to be finite. In comparison, we also tested the conservation behavior of the splitting method introduced in \cite{Gaucknm} for $d=1$. Note that, for $d=3$, the choice $s_2=s_0=0.5$ corresponds to the borderline case of our convergence analysis, where the theoretical convergence order vanishes. Thus, this experiment is intended to illustrate the conservation properties rather than verify the convergence results. Simulations were carried out to the final time $T = 1000$ with the following parameters:
\begin{enumerate}
\item [(i)] $d=1$: $N=2^9$, $\tau=2^{-16}$; CPU time is around $1.6$ hours for both methods;
\item [(ii)] $d=2$: $N=2^8$, $\tau=2^{-15}$; GPU time is around $5.5$ hours;
\item [(iii)] $d=3$: $N=2^7$, $\tau=2^{-13}$; GPU time is around $4.5$ hours.
\end{enumerate}

\begin{figure}
\begin{center}
\subfigure{\includegraphics[height=0.35\textwidth,width=0.48\textwidth]{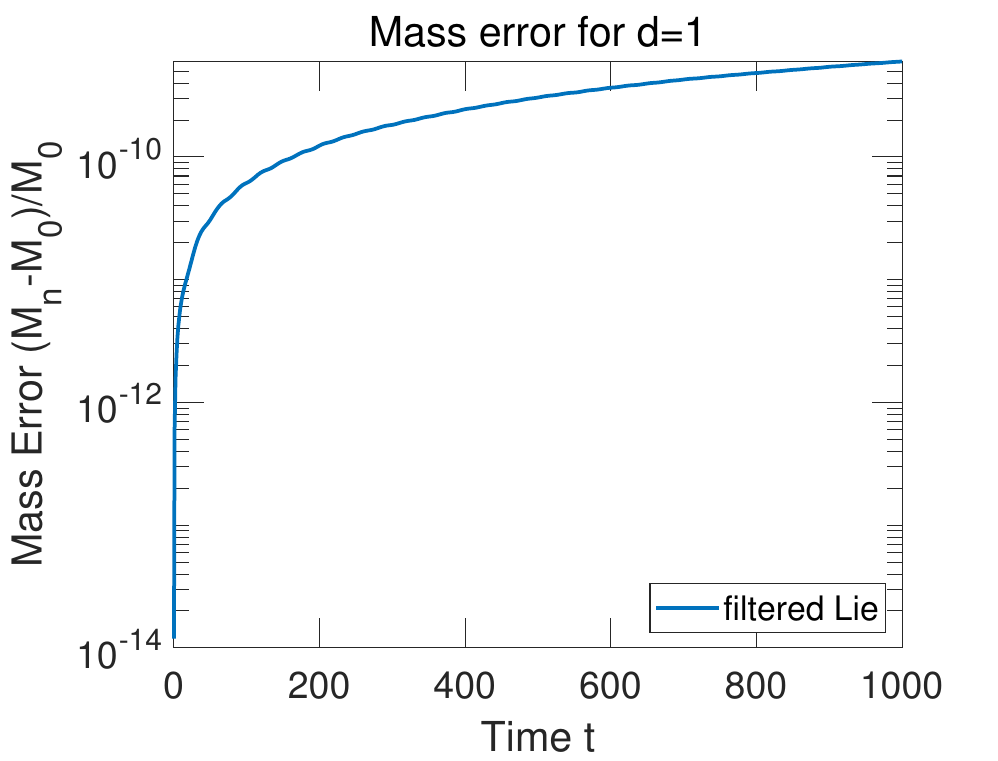}}
\subfigure{\includegraphics[height=0.35\textwidth,width=0.48\textwidth]{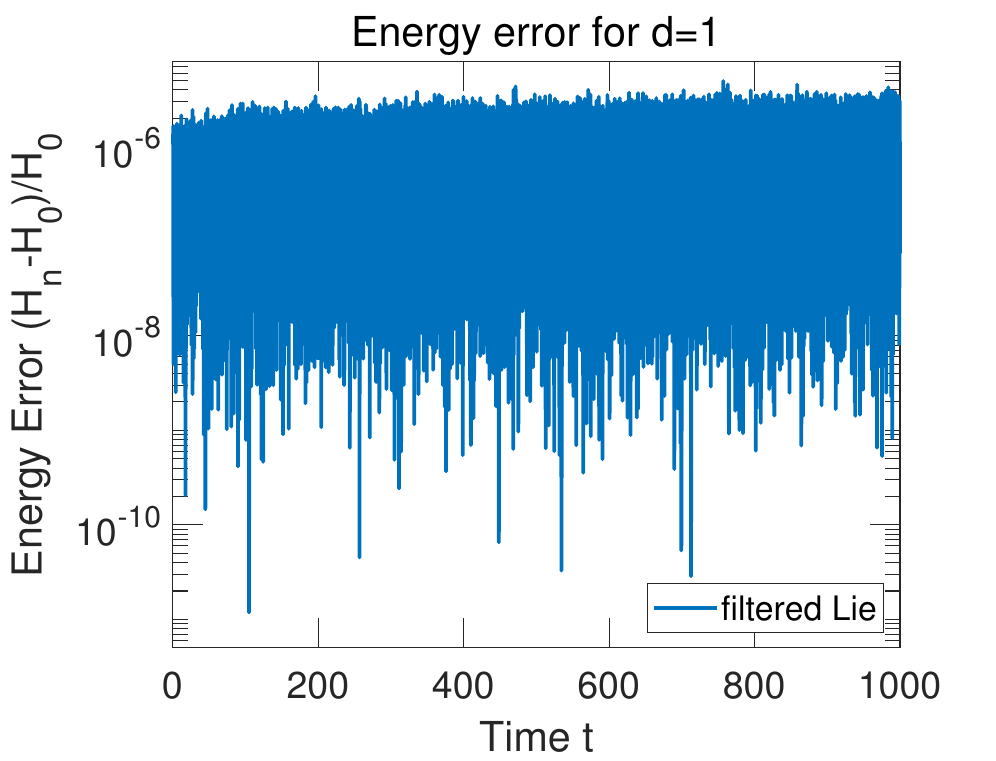}}
\subfigure{\includegraphics[height=0.35\textwidth,width=0.48\textwidth]{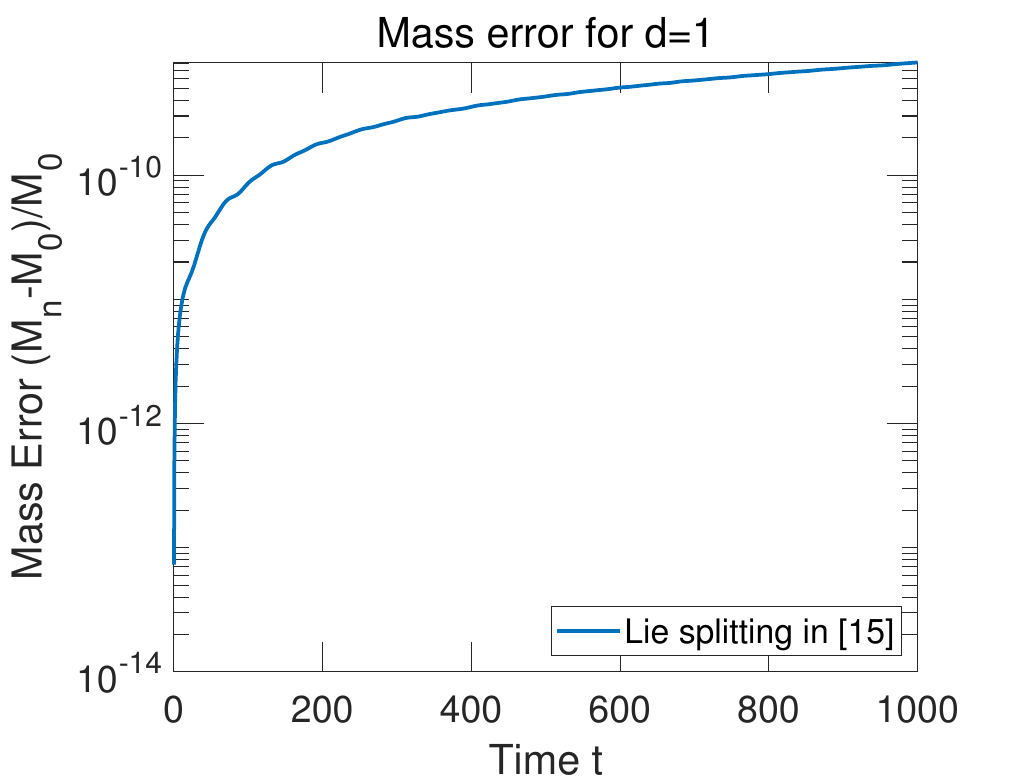}}
\subfigure{\includegraphics[height=0.35\textwidth,width=0.48\textwidth]{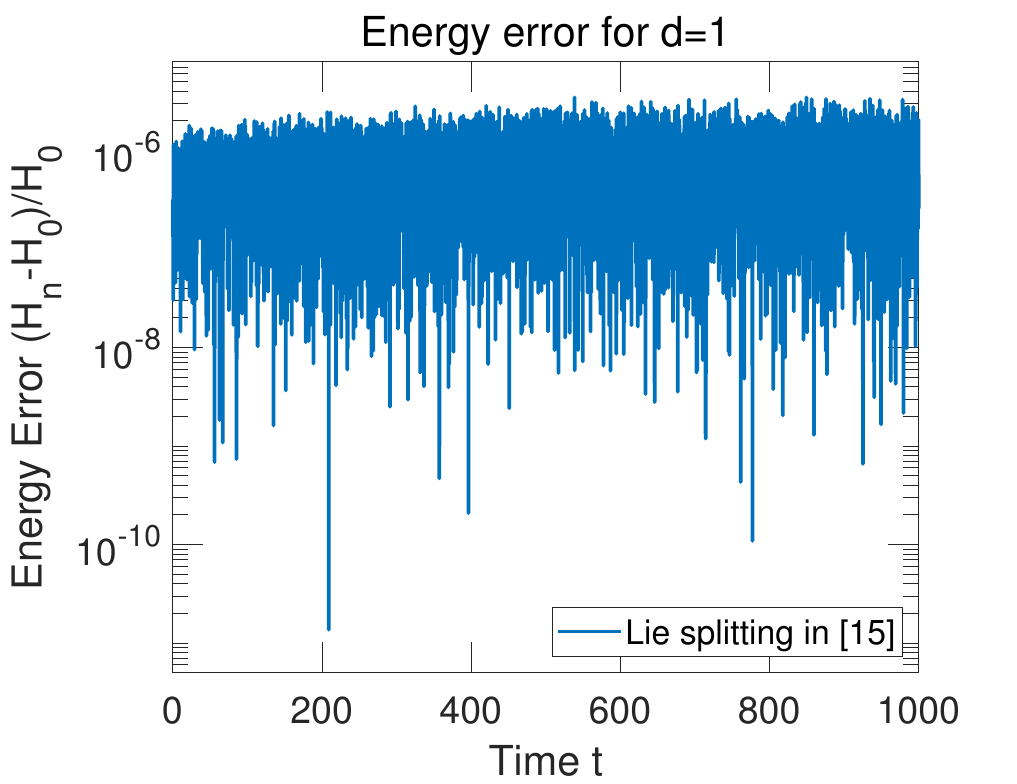}}
\end{center}
\vspace{-2mm}
\caption{The mass and energy error of the one-dimensional ZS. Top left: mass error of \eqref{filtlie}; top right: energy error of \eqref{filtlie}; bottom left: mass error of the method in \cite{Gaucknm}; top right: energy error of the method in \cite{Gaucknm}.}\label{fig:consv1}
\end{figure}

\begin{figure}
\begin{center}
\subfigure{\includegraphics[height=0.35\textwidth,width=0.48\textwidth]{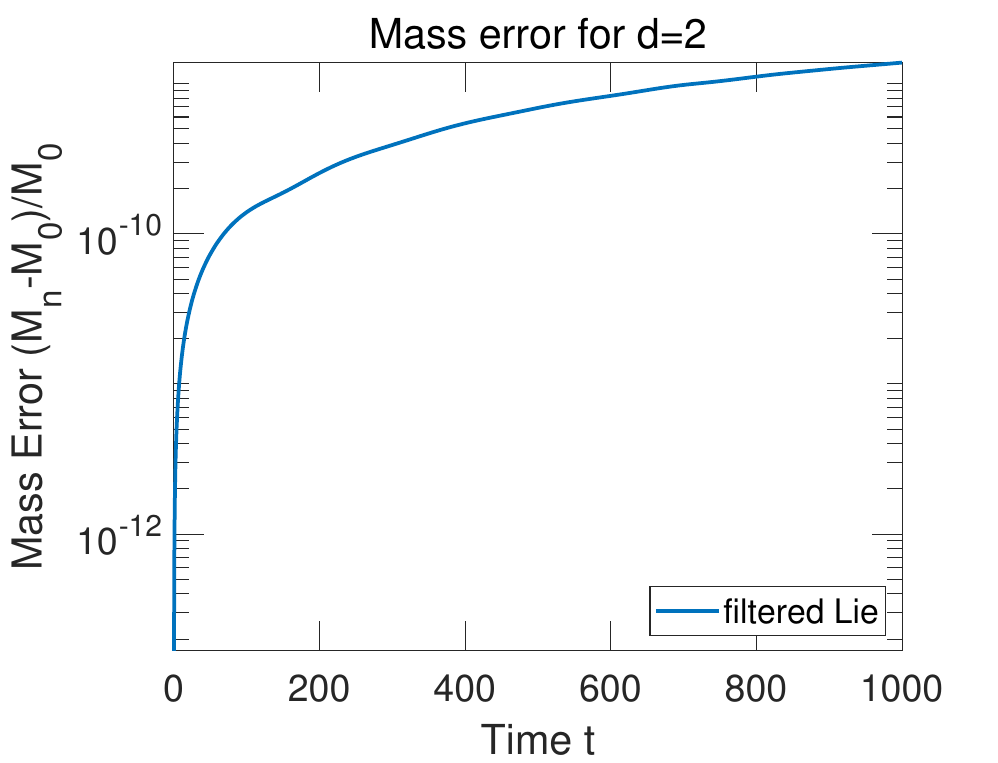}}
\subfigure{\includegraphics[height=0.35\textwidth,width=0.48\textwidth]{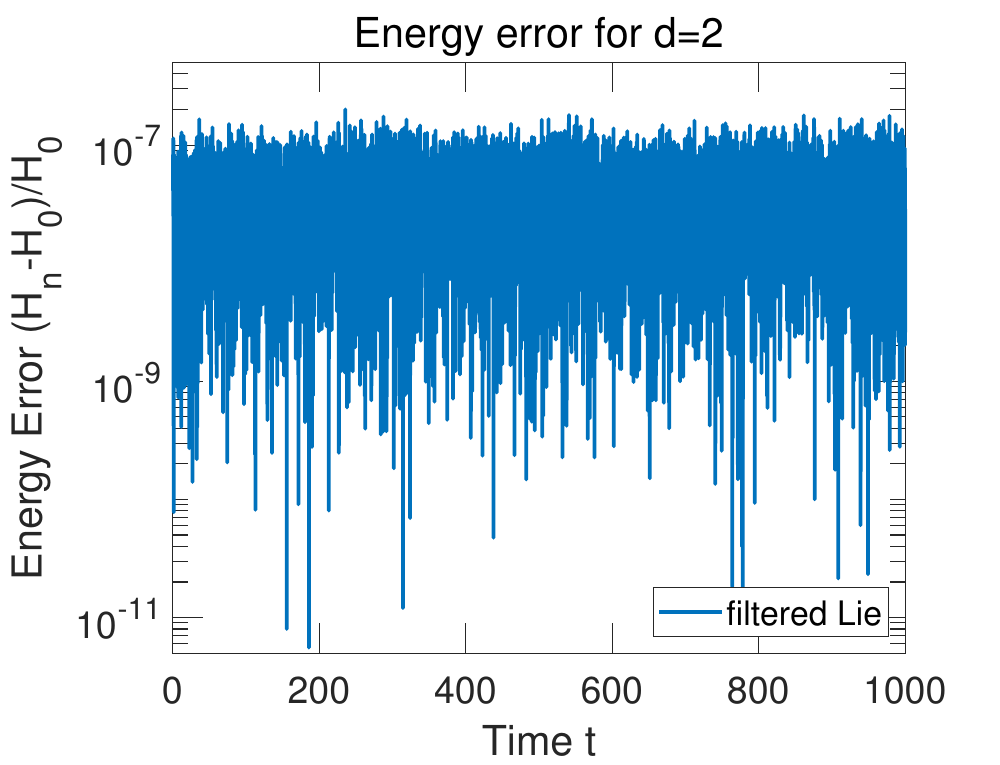}}
\end{center}
\vspace{-2mm}
\caption{The mass and energy error of the two-dimensional ZS. Left: mass error; right: energy error.}\label{fig:consv2}
\end{figure}

\begin{figure}
\begin{center}
\subfigure{\includegraphics[height=0.35\textwidth,width=0.48\textwidth]{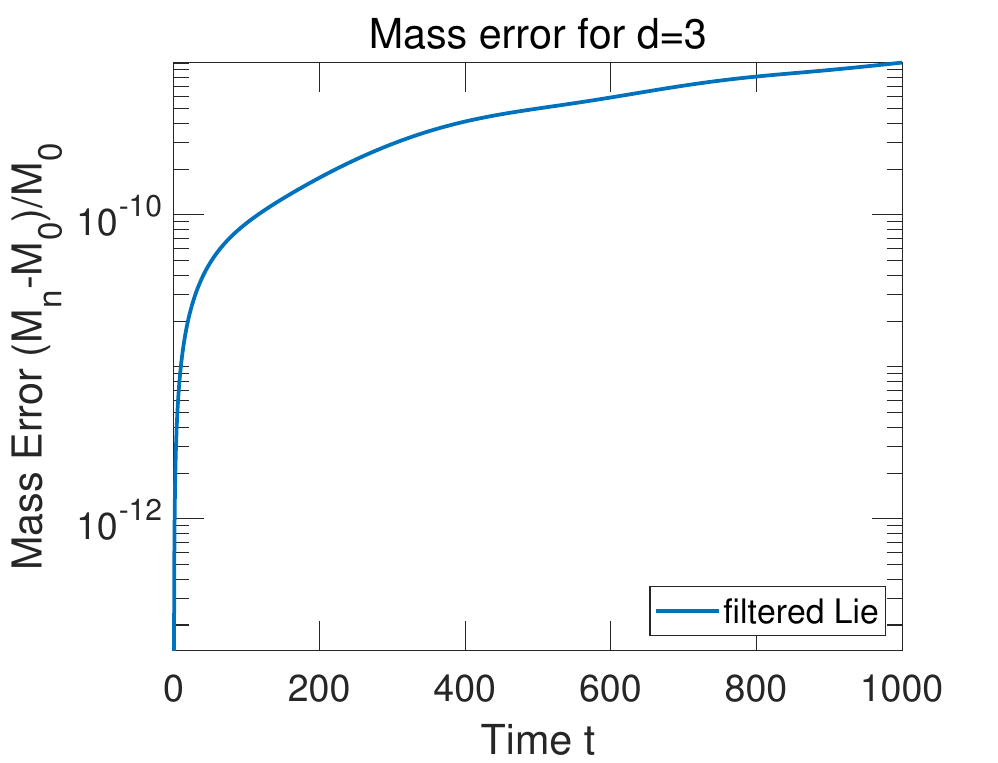}}
\subfigure{\includegraphics[height=0.35\textwidth,width=0.48\textwidth]{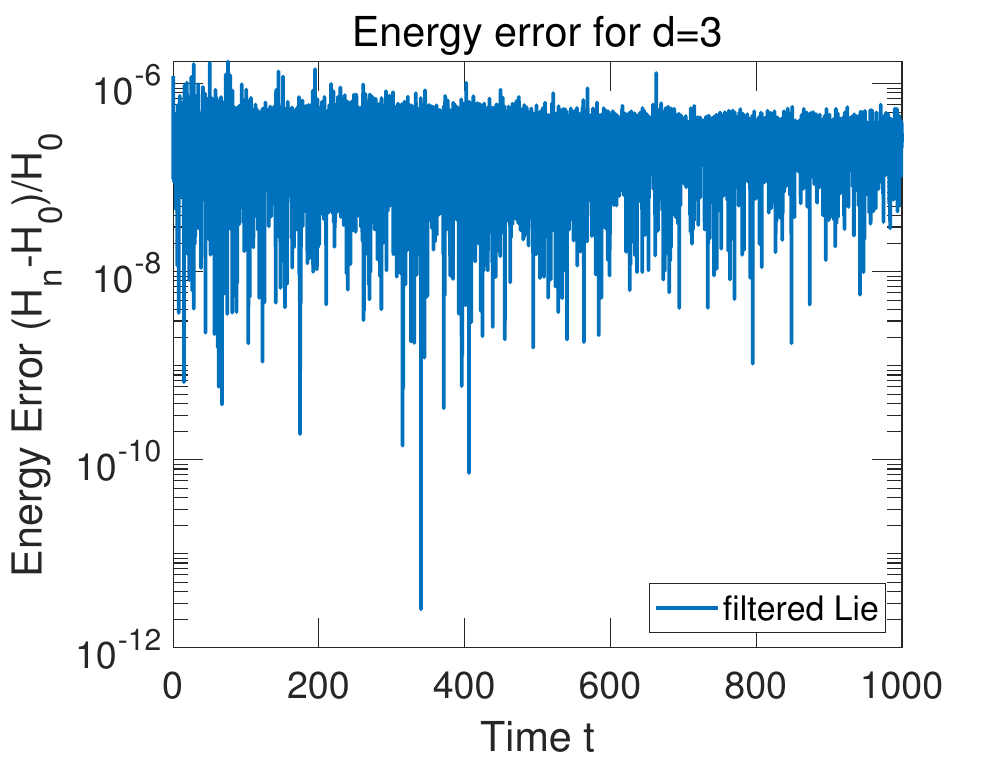}}
\end{center}
\vspace{-2mm}
\caption{The mass and energy error of the three-dimensional ZS. Left: mass error; right: energy error.}\label{fig:consv3}
\end{figure}

Figures~\ref{fig:consv1}--\ref{fig:consv3} illustrate the behavior of the mass and energy errors over the tested long time intervals. Under the CFL condition, our scheme conserves discrete mass in exact arithmetic. The small mass deviations observed in the computations are consistent with accumulated floating-point round-off errors, including those introduced by the FFT and inverse FFT, as also discussed in~\cite{Einkcma}. Figure~\ref{fig:consv1} further shows that our method \eqref{filtlie} and the Lie splitting method proposed in \cite{Gaucknm} exhibit nearly identical conservation behavior. Both methods have very small mass errors with similar accumulation over time. Their energy errors exhibit small oscillations on the scale of $10^{-6}$, together with a slight temporal drift. Thus, the experiments indicate approximate energy conservation over the tested time intervals, alongside discrete mass conservation up to floating-point errors. A rigorous analysis of the long-time energy behavior is a subject of future work.

\section{Conclusion}\label{conclusion}

We developed a fully discrete, filtered Lie splitting scheme for the Zakharov system and established the error estimates. By employing  multilinear estimates in discrete Bourgain spaces, we proved convergence under low regularity assumptions that match those of the continuous well-posedness theory. The resulting error is of order $\mathcal{O}(\tau^{s/2}+N^{-s})$, where $s$ denotes the regularity gap between the space in which the exact solution resides in and the space in which the error is measured. The numerical experiments show small relative energy errors with a visible envelope or drift over the displayed time interval. Numerical experiments corroborate the theoretical predictions and highlight the robustness of the approach for rough solutions.

\begin{appendices}
\section{Proof of Theorem \ref{theomult}}\label{sectionproof}

In this section, we establish Theorem~\ref{theomult}. We begin by presenting a sequence of inequalities \eqref{multdisc1}-\eqref{multdisc5}, from which \eqref{discmult1}-\eqref{discmult5} follows directly. We then prove a stronger one-dimensional estimate that immediately implies \eqref{multdisc1}-\eqref{multdisc4}. Next, we extend these estimates to higher dimensions $d\ge 2$. Finally, we show that both \eqref{multdisc2} and \eqref{multdisc5} are, in fact, corollaries of \eqref{multdisc1} in all dimensions.

To demonstrate Theorem~\ref{theomult}, it suffices to show the following:
\begin{align}
\label{multdisc1}	&\|\Pi_\nu(\Pi_\nu v_n\Pi_\nu w_n)\|_{X_{1,\tau}^{s_2+\frac12,-b}}+\|\Pi_\nu(\Pi_\nu\overline{v_n}\Pi_\nu w_n)\|_{X_{1,\tau}^{s_2+\frac12,-b}}\lesssim\|v_n\|_{X_{2,\tau}^{s_2,b}}\|w_n\|_{X_{1,\tau}^{s_2+\frac12,b}},\\
\label{multdisc3}	&\|\Pi_\nu(\Pi_\nu v_n\Pi_\nu w_n)\|_{X_{1,\tau}^{s_0+\frac12,-b}}+\|\Pi_\nu(\Pi_\nu\overline{v_n}\Pi_\nu w_n)\|_{X_{1,\tau}^{s_0+\frac12,-b}}\lesssim\|v_n\|_{X_{2,\tau}^{s_0,b}}\|w_n\|_{X_{1,\tau}^{s_2+\frac12,b}},\\
\label{multdisc4}	&\|\Pi_\nu(\Pi_\nu v_n\Pi_\nu w_n)\|_{X_{1,\tau}^{s_0+\frac12,-b}}+\|\Pi_\nu(\Pi_\nu\overline{v_n}\Pi_\nu w_n)\|_{X_{1,\tau}^{s_0+\frac12,-b}}\lesssim\|v_n\|_{X_{2,\tau}^{s_2,b}}\|w_n\|_{X_{1,\tau}^{s_0+\frac12,b}},\\
\label{multdisc2}	&\|\Pi_\nu(\Pi_\nu v_n\Pi_\nu\overline{w_n})\|_{X_{2,\tau}^{s_2+1,-b}}+
\|\Pi_\nu(\Pi_\nu\overline{v_n}\Pi_\nu w_n)\|_{X_{2,\tau}^{s_2+1,-b}}\lesssim\|v_n\|_{X_{1,\tau}^{s_2+\frac12,b}}\|w_n\|_{X_{1,\tau}^{s_2+\frac12,b}},\\
\label{multdisc5}&\|\Pi_\nu(\Pi_\nu v_n\Pi_\nu\overline{w_n})\|_{X_{2,\tau}^{s_0+1,-b}}+\|\Pi_\nu(\Pi_\nu\overline{v_n}\Pi_\nu w_n)\|_{X_{2,\tau}^{s_0+1,-b}}\lesssim\|v_n\|_{X_{1,\tau}^{s_2+\frac12,b}}\|w_n\|_{X_{1,\tau}^{s_0+\frac12,b}},
\end{align}
where $\nu=\frac\theta4\geq\frac\tau4$. Indeed, if we substitute $\Pi_\theta v_n$ and $\Pi_\theta w_n$ into the equations above, the resulting estimates coincide with those stated in Theorem~\ref{theomult}, since $\Pi_\nu\Pi_\theta u_n=\Pi_\theta u_n$ and $\Pi_\nu(\Pi_\theta v_n\Pi_\theta w_n)=\Pi_\theta v_n\Pi_\theta w_n$.

 Throughout this appendix, when estimating multilinear forms in Fourier space, we may replace the space-time Fourier transforms of the input functions by their absolute values. This leaves the corresponding Bourgain norms unchanged (see also \eqref{discdef1}, \eqref{discdef2}) and can only increase the convolution bounds.

\subsection{Proof of \eqref{multdisc1}--\eqref{multdisc4} in 1D}

\begin{proposition}\label{prop1d}
Let $d=1$. For $\nu\geq\frac\tau4$ and for arbitrary sequences $\{v_n\}$ and $\{w_n\}$ belonging to the appropriate discrete Bourgain spaces, we have the following estimate:
\begin{equation}\label{mult1d}
\|\Pi_\nu(\Pi_\nu v_n\Pi_\nu w_n)\|_{X_{1,\tau}^{\frac12,-\frac38}}+\|\Pi_\nu(\Pi_\nu\overline{v_n}\Pi_\nu w_n)\|_{X_{1,\tau}^{\frac12,-\frac38}}\lesssim\|v_n\|_{X_{2,\tau}^{0,\frac38}}
\|w_n\|_{X_{1,\tau}^{\frac12,\frac38}}.
\end{equation}
\end{proposition}
\begin{proof}
By duality, the inequality
$$
\|\Pi_\nu(\Pi_\nu v_{0,n}\Pi_\nu v_{2,n})\|_{X_{1,\tau}^{\frac12,-\frac38}}\lesssim\|v_{0,n}\|_{X_{2,\tau}^{0,\frac38}}\|v_{2,n}\|_{X_{1,\tau}^{\frac12,\frac38}}
$$
is equivalent to
\begin{equation}\label{multeqv}
S=\Big|\tau\sum_n\int_{\mathbb{T}}\Pi_\nu v_{0,n}\Pi_\nu\overline{v_{1,n}}\Pi_\nu v_{2,n}dx\Big|\lesssim\|v_{0,n}\|_{X_{2,\tau}^{0,\frac38}}\|v_{1,n}\|_{X_{1,\tau}^{-\frac12,\frac38}}\|v_{2,n}\|_{X_{1,\tau}^{\frac12,\frac38}}.
\end{equation}
Expanding the functions in Fourier space yields
$$
S=\Big|\sum_{k_1,k_2}\int_{\sigma_1,\sigma_2}\widetilde{\Pi_\nu v_{0,n}}(\sigma_0,k_0)\widetilde{\Pi_\nu\overline{v_{1,n}}}(\sigma_1,k_1)\widetilde{\Pi_\nu v_{2,n}}(\sigma_2,k_2)d\sigma_1d\sigma_2\Big|,
$$
where $k_0=k_1-k_2$, and $\sigma_0=\sigma_1-\sigma_2+2l\pi\tau^{-1}$ for some $l\in\mathbb{Z}$. For simplicity, we denote
$$
\mu_0=\sigma_0-|k_0|,\quad \mu_1=\sigma_1+|k_1|^2, \quad \mu_2=\sigma_2+|k_2|^2.
$$
We split the resulting sum $S$ into four contributions according to the frequency regions:
\begin{align*}
S_1&=\mathbbm{1}_{\langle k_1\rangle\leq 4\langle k_2\rangle}S,\qquad\qquad
&S_2=\mathbbm{1}_{\langle k_1\rangle>4\langle k_2\rangle,\langle k_1\rangle\lesssim\langle d_\tau(\mu_0)\rangle^{\frac12}}S,\\
S_3&=\mathbbm{1}_{\langle k_1\rangle>4\langle k_2\rangle,\langle k_1\rangle\lesssim\langle d_\tau(\mu_1)\rangle^{\frac12}}S,\qquad
&S_4=\mathbbm{1}_{\langle k_1\rangle>4\langle k_2\rangle,\langle k_1\rangle\lesssim\langle d_\tau(\mu_2)\rangle^{\frac12}}S,
\end{align*}
where $\mathbbm{1}_{P}$ denotes the characteristic function of the set $P$. We first show that
\begin{equation}\label{s1234}
S\leq S_1+S_2+S_3+S_4.
\end{equation}
It suffices to prove that when $\langle k_1\rangle>4\langle k_2\rangle$, it holds
$$
\langle k_1\rangle\lesssim\langle d_\tau(\mu_m)\rangle^{\frac12},
$$
where $\mu_m=\max\{|\mu_0|,|\mu_1|,|\mu_2|\}$. Since $d_\tau$ is a periodic function with period $2\pi\tau^{-1}$, we assume without loss of generality that $\mu_j\in[-\pi\tau^{-1},\pi\tau^{-1})$ for $j=0,1,2$; hence $\mu_m<a\tau^{-1}$ with $a>0$ to be decided, otherwise $\mu_m\gtrsim\tau^{-1}$ and the desired estimate follows from
$$
\langle k_1\rangle\lesssim\tau^{-\frac12}\lesssim\langle d_\tau(\mu_m)\rangle^{\frac12}.
$$
Under these assumptions we obtain
$$
\langle\mu_j\rangle\lesssim\langle d_\tau(\mu_j)\rangle,\quad j=0,1,2,
$$
and
\begin{align*}
3a\tau^{-1}&\geq\mu_0-\mu_1+\mu_2=\sigma_0-\sigma_1+\sigma_2-|k_1-k_2|-|k_1|^2+|k_2|^2\\
&=2l\pi\tau^{-1}-|k_1-k_2|-|k_1|^2+|k_2|^2.
\end{align*}
On the support of $\Pi_\nu$ with $\nu\geq\tau/4$, the bound $|k_j|\leq\sqrt{c/\tau}$ yields
$$
\big||k_1-k_2|+|k_1|^2-|k_2|^2\big|\leq2\sqrt{c/\tau}+c/\tau<(2\pi-3a)/\tau
$$
for $a=(2\pi-c)/4>0$, when $\tau<(2\pi-c)^2/(64c)$. Together with $|\mu_0-\mu_1+\mu_2|<3a/\tau$ whenever $|\mu_j|<a/\tau$ for $j=0,1,2$, this estimate forces $l=0$ in the identity above and hence $\sigma_0=\sigma_1-\sigma_2$. Consequently, when $\langle k_1\rangle>4\langle k_2\rangle$, we have
\begin{equation}\label{n1lm}
|k_1|^2\lesssim |k_1|^2-|k_2|^2+|k_1-k_2|=\sigma_1+|k_1|^2-\sigma_2-|k_2|^2-\sigma_0+|k_0|=\mu_1-\mu_2-\mu_0,
\end{equation}
and thus \eqref{s1234} is established.

We also require the auxiliary estimate
\begin{equation}\label{jems}
\|\Pi_\nu v_n\|_{l_\tau^4L^4}\lesssim\|v_n\|_{X_{1,\tau}^{0,\frac38}},
\end{equation}
for $\{v_n\}$ in $X_{1,\tau}^{0,\frac38}$; a proof of this is given in \cite{Ostjems}.

Now we return to the physical space and show
\begin{equation}\label{i1-4}
S_i\lesssim\|v_{0,n}\|_{X_{2,\tau}^{0,\frac38}}\|v_{1,n}\|_{X_{1,\tau}^{-\frac12,\frac38}}
\|v_{2,n}\|_{X_{1,\tau}^{\frac12,\frac38}},\quad i=1,2,3,4.
\end{equation}
Writing $\widetilde{w_{j,n}}(\sigma,k)=|\widetilde{v_{j,n}}(\sigma,k)|,~j=0,1,2$, we have $\|w_{j,n}\|_{X_{l,\tau}^{s,b}}=\|v_{j,n}\|_{X_{l,\tau}^{s,b}}$, for $l=1,2$, $j=0,1,2$ and any $s,b\in\mathbb{R}$.

For $S_1$, applying \eqref{discl4}, \eqref{jems} and H\"older's inequality yields
\begin{align*}
S_1&\lesssim\Big|\sum_{k_1,k_2}\int_{\sigma_1,\sigma_2}|\widetilde{\Pi_\nu v_{0,n}}(\sigma_0,k_0)||\langle k_1\rangle^{-\frac12}\widetilde{\Pi_\nu\overline{v_{1,n}}}(\sigma_1,k_1)||\langle k_2\rangle^{\frac12}\widetilde{\Pi_\nu v_{2,n}}(\sigma_2,k_2)|d\sigma_1d\sigma_2\Big|\\
&=\Big|\tau\sum_n\int_{\mathbb{T}}\Pi_\nu w_{0,n}\langle\partial_x\rangle^{-\frac12}\Pi_\nu\overline{w_{1,n}}\langle\partial_x\rangle^\frac12\Pi_\nu w_{2,n}dx\Big|\\
&\lesssim\|\Pi_\nu w_{0,n}\|_{l_\tau^2 L^2}\|\langle\partial_x\rangle^{-\frac12}\Pi_\nu\overline{w_{1,n}}\|_{l_\tau^4 L^4}\|\langle\partial_x\rangle^\frac12\Pi_\nu w_{2,n}\|_{l_\tau^4 L^4}\\
&\lesssim\|v_{0,n}\|_{X_{2,\tau}^{0,0}}\|v_{1,n}\|_{X_{1,\tau}^{-\frac12,\frac38}}\|v_{2,n}\|_{X_{1,\tau}^{\frac12,\frac38}}\leq\|v_{0,n}\|_{X_{2,\tau}^{0,\frac38}}\|v_{1,n}\|_{X_{1,\tau}^{-\frac12,\frac38}}\|v_{2,n}\|_{X_{1,\tau}^{\frac12,\frac38}}.
\end{align*}

For the remaining cases note that $\langle k_1\rangle>4\langle k_2\rangle$ and $k_0=k_1-k_2$ imply $k_0\sim k_1$. Let $D_\tau$ denote the time-difference $\big(D_\tau(u_n)\big)_n=\big(\frac{u_n-u_{n-1}}{\tau}\big)_n$ (see also \eqref{deff}).
Utilizing \eqref{discl4}, \eqref{jems}, the Sobolev embedding theorem and H\"older's inequality, we obtain
\begin{align*}
S_2&\lesssim\Big|\sum_{k_1,k_2}\int_{\sigma_1,\sigma_2}|\langle d_\tau(\mu_0)\rangle^{\frac14}\widetilde{\Pi_\nu v_{0,n}}(\sigma_0,k_0)||\langle k_1\rangle^{-\frac12}\widetilde{\Pi_\nu\overline{v_{1,n}}}(\sigma_1,k_1)||\widetilde{\Pi_\nu v_{2,n}}(\sigma_2,k_2)|d\sigma_1d\sigma_2\Big|\\
&\sim\Big|\tau\sum_n\int_{\mathbb{T}}\langle D_\tau\rangle^\frac14\Pi_\nu w_{0,n}\langle\partial_x\rangle^{-\frac12}\Pi_\nu\overline{w_{1,n}}\Pi_\nu w_{2,n}dx\Big|\\
&\lesssim\|\langle D_\tau\rangle^\frac14\Pi_\nu w_{0,n}\|_{l_\tau^2 L^2}\|\langle\partial_x\rangle^{-\frac12}\Pi_\nu\overline{w_{1,n}}\|_{l_\tau^4 L^4}\|\Pi_\nu w_{2,n}\|_{l_\tau^4 L^\infty}\\ &\lesssim\|v_{0,n}\|_{X_{2,\tau}^{0,\frac14}}\|v_{1,n}\|_{X_{1,\tau}^{-\frac12,\frac38}}\|v_{2,n}\|_{l_\tau^4 W^{\frac14+\varepsilon,4}}\leq\|v_{0,n}\|_{X_{2,\tau}^{0,\frac38}}\|v_{1,n}\|_{X_{1,\tau}^{-\frac12,\frac38}}\|v_{2,n}\|_{X_{1,\tau}^{\frac12,\frac38}},
\end{align*}
where $\varepsilon\in(0,\frac14)$. Similarly, for $\varepsilon\in(0,\frac18)$, we obtain
\begin{align*}	
S_3&\lesssim\Big|\sum_{k_1,k_2}\int_{\sigma_1,\sigma_2}|\widetilde{\Pi_\nu v_{0,n}}(\sigma_0,k_0)||\langle d_\tau(\mu_1)\rangle^\frac14\langle k_1\rangle^{-\frac12}\widetilde{\Pi_\nu\overline{v_{1,n}}}
(\sigma_1,k_1)||\widetilde{\Pi_\nu v_{2,n}}(\sigma_2,k_2)|d\sigma_1d\sigma_2\Big|\\
&\sim\Big|\tau\sum_n\int_{\mathbb{T}}\Pi_\nu w_{0,n}\langle D_\tau\rangle^\frac14\langle\partial_x\rangle^{-\frac12}\Pi_\nu\overline{w_{1,n}}\Pi_\nu w_{2,n}dx\Big|\\
&\lesssim\|\Pi_\nu w_{0,n}\|_{l_\tau^4 L^2}\|\langle D_\tau\rangle^\frac14\langle \partial_x\rangle^{-\frac12}\Pi_\nu w_{1,n}\|_{l_\tau^2 L^2}\|\Pi_\nu w_{2,n}\|_{l_\tau^4 L^\infty}\\
&\lesssim\|v_{0,n}\|_{X_{2,\tau}^{0,\frac14+\varepsilon}}\|v_{1,n}\|_{X_{1,\tau}^{-\frac12,\frac14}}\|v_{2,n}\|_{l_\tau^4 W^{\frac14+\varepsilon,4}}\leq\|v_{0,n}\|_{X_{2,\tau}^{0,\frac38}}\|v_{1,n}\|_{X_{1,\tau}^{-\frac12,\frac38}}\|v_{2,n}\|_{X_{1,\tau}^{\frac12,\frac38}}
\end{align*}
and
\begin{align*}
S_4&\lesssim\Big|\sum_{k_1,k_2}\int_{\sigma_1,\sigma_2}|\widetilde{\Pi_\nu v_{0,n}}(\sigma_0,k_0)||\langle k_1\rangle^{-\frac12}\widetilde{\Pi_\nu\overline{v_{1,n}}}
(\sigma_1,k_1)||\langle d_\tau(\mu_2)\rangle^{\frac14}\widetilde{\Pi_\nu v_{2,n}}(\sigma_2,k_2)|d\sigma_1d\sigma_2\Big|\\
&\sim\Big|\tau\sum_n\int_{\mathbb{T}}\Pi_\nu w_{0,n}\langle\partial_x\rangle^{-\frac12}\Pi_\nu\overline{w_{1,n}}\langle D_\tau\rangle^\frac14\Pi_\nu w_{2,n}dx\Big|\\
&\lesssim\|\Pi_\nu w_{0,n}\|_{l_\tau^4 L^2}\|\langle\partial_x\rangle^{-\frac12}\Pi_\nu w_{1,n}\|_{l_\tau^4 L^4}\|\langle D_\tau\rangle^\frac14\Pi_\nu w_{2,n}\|_{l_\tau^2 L^4}\\ &\lesssim\|v_{0,n}\|_{X_{2,\tau}^{0,\frac14+\varepsilon}}\|v_{1,n}\|_{X_{1,\tau}^{-\frac12,\frac38}}\|v_{2,n}\|_{X_{1,\tau}^{\frac14,\frac14}}\leq\|v_{0,n}\|_{X_{2,\tau}^{0,\frac38}}\|v_{1,n}\|_{X_{1,\tau}^{-\frac12,\frac38}}\|v_{2,n}\|_{X_{1,\tau}^{\frac12,\frac38}}.
\end{align*}
Thus we establish \eqref{i1-4}, which yields \eqref{multeqv} and correspondingly, the first part of \eqref{mult1d}.

Replacing $(\sigma_0,k_0)$ by $(-\sigma_0,-k_0)$ in the foregoing argument produces the second part of \eqref{mult1d}, i.e.,
$$
\|\Pi_\nu(\Pi_\nu\overline{v_n}\Pi_\nu w_n)\|_{X_{1,\tau}^{\frac12,-\frac38}}\lesssim\|v_n\|_{X_{2,\tau}^{0,\frac38}}\|w_n\|_{X_{1,\tau}^{\frac12,\frac38}},
$$
which completes the proof.
\end{proof}

\noindent\emph{Proof of \eqref{multdisc1}-\eqref{multdisc4} for $d=1$.} Estimate \eqref{multdisc1} follows from the generalized Leibniz rule \cite{Muscam} together with \eqref{mult1d}. Moreover, it is evident that \eqref{multdisc3} and \eqref{multdisc4} are direct corollaries of \eqref{mult1d}.

\subsection{Proof of \eqref{multdisc1}--\eqref{multdisc4} for $d\geq2$}\label{subsectprod2}

We begin by establishing estimate \eqref{multdisc1}. Without loss of generality we work with the component $v_{0,n}$; the same argument applies to its complex conjugate. Our approach follows the discrete Littlewood-Paley decomposition used in \cite{Baiima,Jiima} and \cite{Ostjems}, and adapts the corresponding continuous argument in \cite{Kishijam}.

Let us first recall the Littlewood-Paley decomposition. For $\sigma\in I_{\tau}=[-\frac\pi\tau,\frac\pi\tau)$ and  $m\in\mathbb{N}_0=\{0\}\cup\mathbb{N}^*$, we define
$$
\mathbbm{1}_m(\sigma)=\mathbbm{1}_{2^m\leq\langle\sigma\rangle<2^{m+1}\cap I_\tau},
$$
and denote by $\mathbbm{1}_{m}$ their $\frac{2\pi}\tau$-periodic extensions. For a sequence $\{v_n(x)\}$, define the operators by the Fourier multipliers:
$$
\widetilde{P_mv_n}(\sigma,k)=\mathbbm{1}_m(\sigma-|k|)
\widetilde{v_n}(\sigma,k),\quad
\widetilde{P'_mv_n}(\sigma,k)=\mathbbm{1}_m(\sigma+|k|^2)
\widetilde{v_n}(\sigma,k).
$$
Moreover, we set
$$
P_{\leq M}v_n=\sum_{0\leq m\leq M}P_mv_n,
$$
and the localizers
$$
\widetilde{Q_lv_n}(\sigma,k)=\mathbbm{1}_{2^l\leq\langle k\rangle<2^{l+1}} \widetilde{v_n}(\sigma,k),
$$
thus ensuring that
$$
v_n=\sum_{l,m\geq0}P_mQ_lv_n=\sum_{l,m\geq0}P'_mQ_lv_n.
$$

We set $N=(N_0,N_1,N_2)$, and $L=(L_0,L_1,L_2)$, where $N_j$ and $L_j$ range over dyadic values
$N_j=2^{n_j},L_j=2^{l_j} (j=0,1,2)$. We split
$$
v_{j,n}=\sum\limits_{l_j,n_j\in\mathbb{N}}v_{j,n}^{L_jN_j}, \qquad v_{0,n}^{L_jN_j}=P_{l_j}Q_{n_j}v_{0,n};\quad v_{j,n}^{L_jN_j}=P'_{l_j}Q_{n_j}v_{j,n},\,\,\,j=1,2,
$$
and denote
\begin{equation}\label{defSLN}
S(L,N)=\Big|\tau\sum_n\int_{\mathbb{T}^d}\Pi_\nu v_{0,n}^{L_0N_0}\Pi_\nu \overline{v_{1,n}}^{L_1N_1}\Pi_\nu v_{2,n}^{L_2N_2}dx\Big|.
\end{equation}
Clearly, we have
$$
S\leq\sum_{L,N}S(L,N),
$$
where $S$ is as defined in \eqref{multeqv}.

The subsequent technical lemmas control these dyadic contributions.

\begin{lemma}\label{lemma12}
Denote $N_{\min}=\min\{N_0,N_1,N_2\}$, $\check{L}_{ij}=\min\{L_i,L_j\}$, and $\hat{L}_{ij}=\max\{L_i,L_j\}$. Then we have the following estimates:
\begin{align}
\label{v1v2}\|P_{n_0}Q_{l_0}(\Pi_\nu v_{1,n}^{L_1N_1}\Pi_\nu\overline{v_{2,n}}^{L_2N_2})\|_{l_\tau^2L^2}&\lesssim N_{\min}^{\frac d2-1+\varepsilon}L_1^\frac12L_2^\frac12\|v_{1,n}^{L_1N_1}\|_{l_\tau^2L^2}\|v_{2,n}^{L_2N_2}\|_{l_\tau^2L^2},\\
\label{v0v1}\|P'_{n_2}Q_{l_2}(\Pi_\nu v_{0,n}^{L_0N_0}\Pi_\nu v_{1,n}^{L_1N_1})\|_{l_\tau^2L^2}&\lesssim N_{\min}^{\frac{d-1}2}\check{L}_{01}^\frac12(\frac{\hat{L}_{01}}{N_1}+1)^\frac12 \|v_{0,n}^{L_0N_0}\|_{l_\tau^2L^2}\|v_{1,n}^{L_1N_1}\|_{l_\tau^2L^2},
\end{align}
where $\nu\geq\frac\tau4$, and $\varepsilon>0$ can be taken to be arbitrarily small.
\end{lemma}
\begin{proof}
For notational convenience, we denote $$u_{j,n}=v_{j,n}^{L_jN_j},\quad j=0,1,2;\qquad E(B)=\cup_{m\in\mathbb{Z}}[\frac{2m\pi}{\tau}-B,\frac{2m\pi}{\tau}+B).$$ Since $|k_j|\lesssim\tau^{-\frac12}$ and $\sigma_j\in I_\tau$, the union above reduces to a uniformly bounded number of intervals, $E(B)=\cup_{|m|\leq M}[\frac{2m\pi}{\tau}-B,\frac{2m\pi}{\tau}+B)$ with $M=\mathcal{O}(1)$. Without loss of generality we alsonormalize $\|\Pi_\nu u_{j,n}\|_{l_\tau^2L^2}=1$ for $j=0,1,2$.

We first prove \eqref{v1v2}. The cases $N_{\min}=N_1$ or $N_{\min}=N_2$ have been treated in \cite[Lemma~8.3]{Jiima}. Note that the same argument extends to $d\geq3$ (see \cite{Jiicos}). Hence it suffices to consider $N_{\min}=N_0$. Here, we follow the lines of the proof in \cite{Jiima}. Since the estimate is symmetric in $u_{1,n}$ and $u_{2,n}$, we may assume, without loss of generality, that $L_1\leq L_2$.

Decompose $\widetilde{u_{j,n}},~j=1,2$  into $\mathcal{O}\big((\frac{N_j}{N_0})^d\big)$ pieces, each supported in a spatial-frequency cube of side length $N_0$. Let $R_a$ denote  the localization to the cube centered at $a\in\mathbb{R}^d$. By construction, we have
\begin{equation}\label{defra}
\sum\limits\|R_a u_{1,n}\|^2_{l_\tau^2L^2}=1.
\end{equation}
Note that this summation contains $\mathcal{O}\big((\frac{N_j}{N_0})^d\big)$ terms. By \cite[Lemma~8.2]{Jiima}, each localized piece satisfies
\begin{equation}\label{trans}
\|\Pi_\nu R_a u_{1,n}\|_{l_\tau^4L^4}\lesssim N_0^{\frac{d-2+\varepsilon}4}L_1^\frac12\|R_a u_{1,n}\|_{l_\tau^2L^2},
\end{equation}
and the analogous estimate holds for $u_{2,n}$.

Moreover, the left-hand side norm in \eqref{v1v2} does not vanish unless $a-a^\prime=\mathcal{O}(N_0)$; therefore the sum $\sum\limits_{a,a^\prime}\Pi_\nu(\Pi_\nu R_a u_{1,n}\Pi_\nu R_{a^\prime}\overline{u_{2,n}})$ exhibits quasi-orthogonality. Applying \eqref{defra}, \eqref{trans} and the Cauchy-Schwarz inequality yields
\begin{equation}\label{orth}
\begin{aligned}
\|P_{n_0}(\Pi_\nu u_{1,n}\Pi_\nu\overline{u_{2,n}})\|_{l_\tau^2L^2}&\lesssim\Big(\sum\limits_{a,a^\prime}\|\Pi_\nu(\Pi_\nu R_a u_{1,n}\Pi_\nu R_{a^\prime}\overline{u_{2,n}})\|^2_{l_\tau^2L^2}\Big)^\frac12\\
&\lesssim\Big(\sum\limits_{a,a^\prime}\|\Pi_\nu R_a u_{1,n}\|^2_{l_\tau^4L^4}\|\Pi_\nu R_{a^\prime}\overline{u_{2,n}})\|^2_{l_\tau^4L^4}\Big)^\frac12\\
&\lesssim N_0^{\frac{d-2+\varepsilon}2}L_1^\frac12L_2^\frac12\Big(\sum\limits_a\|R_a u_{1,n}\|^2_{l_\tau^2L^2}\Big)^\frac12\Big(\sum\limits_{a^\prime}\|R_{a^\prime}u_{2,n}\|^2_{l_\tau^2L^2}\Big)^\frac12\\
&\lesssim N_0^{\frac d2-1+\varepsilon}L_1^\frac12L_2^\frac12,
\end{aligned}
\end{equation}
which concludes \eqref{v1v2}.

We now turn to \eqref{v0v1}. By symmetry we assume $L_0\leq L_1$ and follow the argument of \cite[Lemma~2.5]{Kishijam}. Using $\sigma_1-\sigma_0=\sigma_2$ and $k_1-k_0=k_2$ and applying Fubini's theorem gives
$$
\Big\|\sum\limits_{k_0}\int|\widetilde{\Pi_\nu u_{0,n}}(\sigma_0,k_0)|^2|\widetilde{\Pi_\nu u_{1,n}}(\sigma_0+\sigma_2,k_0+k_2)|^2d\sigma_0\Big\|_{L^1l^1}=1.
$$
By Cauchy-Schwarz inequality, the desired norm is controlled by
$$
\Big\|\sum\limits_{k_0}\int\widetilde{\Pi_\nu u_{0,n}}(\sigma_0,k_0)\widetilde{\Pi_\nu u_{1,n}}(\sigma_0+\sigma_2,k_0+k_2)d\sigma_0\|_{L^2l^2}\lesssim\Big(\sum_{k_0}\int_\Omega1d\sigma_1\Big)^\frac12,
$$
where $\Omega=\{\sigma_0\mid\sigma_1-\sigma_0=\sigma_2,\sigma_0-|k_0|\in E(2L_0),\sigma_1+|k_1|^2\in E(2L_1)\}$.
Fix $\sigma_2$ and $k_2$. If it holds
$$
\sigma_2+|k_0|+|k_0+k_2|^2=\sigma_1+|k_1|^2-\sigma_0+|k_0|\in E(B),
$$
with $B=\mathcal{O}(L_1)$, then the integral is $\mathcal{O}(L_0)$; otherwise it vanishes. Since $|k_0|\lesssim |k_1|+|k_2|=|k_2|+\mathcal{O}(N_1)$, the desired norm is bounded by $L_0^\frac12|T(\sigma_2,k_2)|^\frac12$, where
$$
T(\sigma_2,k_2)=\{k_0\in\mathbb{Z}^d\mid |k_0|\sim N_0,~|k_0+k_2|\sim N_1,~|k_0+k_2|^2+\sigma_2+|k_2|\in E(B^\prime)\},
$$
with $B^\prime=\mathcal{O}(L_1+N_1)$.
If $N_{\min}\neq N_2$ and $L_1+N_1\geq N_1^2$, then $k_0$ is confined to a ball of radius $\mathcal{O}(N_{\min})$, leading to
$$
|T(\sigma_2,k_2)|\lesssim N_{\min}^d\lesssim N_1N_{\min}^{d-1}\lesssim\big(\frac{L_1}{N_1}+1\big)N_{\min}^{d-1},
$$
which implies the desired estimate. If $N_{\min}\neq N_2$ and $L_1+N_1\leq N_1^2$, then $k_0$ lies in the intersection of a ball with radius $\mathcal{O}(N_0)$ and an annulus centered at $-k_2$ with radius $\mathcal{O}(N_1)$ and thickness  $\mathcal{O}(\frac{L_1+N_1}{N_1})$. Each such intersection consists of $\mathcal{O}(1)$ connected components, and every component has thickness at least $\mathcal{O}(1)$. Consequently, the number of integer lattice points it contains is comparable to its volume. This yields the  lattice-point counting bound:
$$
|T(\sigma_2,k_2)|\lesssim\frac{L_1+N_1}{N_1}N_{\min}^{d-1}.
$$
Finally, when $N_{\min}=N_2$, one proceeds by applying the same localized decomposition and quasi-orthogonality argument as in \eqref{defra}-\eqref{orth} to establish \eqref{v0v1}.
\end{proof}

We next estimate $S(L,N)$.

\begin{lemma}\label{lemnonres}
Let $S(L,N)$ be as in \eqref{defSLN}. If $N_1\geq8N_2$, then it satisfies
\begin{equation}\label{slnnonres}
S(L,N)\lesssim N_0^{-\varepsilon}N_1^{-\frac12-\varepsilon}N_2^{\frac{d-1}2-\varepsilon}(L_0L_1L_2)^{\frac38+\varepsilon}\|v_{0,n}^{L_0N_0}\|_{l_\tau^2 L^2}\|v_{1,n}^{L_1N_1}\|_{l_\tau^2 L^2}\|v_{2,n}^{L_2N_2}\|_{l_\tau^2 L^2},
\end{equation}
where $\varepsilon>0$ can be taken arbitrarily small.
\end{lemma}
\begin{proof}
This lemma provides an estimate for $S(L,N)$ in the non-resonant regime. Since $k_0=k_1-k_2$ and $N_1\geq8N_2$, we have $N_0\sim N_1$, while the argument used in Proposition~\ref{prop1d} (cf. \eqref{n1lm}) gives $N_1^2\lesssim L_{\max}=\max\{L_0,L_1,L_2\}$. Writing $u_{j,n}=v_{j,n}^{L_jN_j}$ for $j=0, 1, 2$, we normalize $\|u_{j,n}\|_{l_\tau^2L^2}=1$ by homogeneity and take $0<\varepsilon\leq1/8$ throughout the proof.

For $L_{\max}=L_0$, H\"older's inequality together with \eqref{discl4} and the Sobolev embedding yields
$$
S(L,N)\lesssim\|\Pi_\nu u_{1,n}\|_{l_\tau^4L^2}\|\Pi_\nu u_{2,n}\|_{l_\tau^4L^\infty}\lesssim N_2^{\frac d4+\frac{\varepsilon}{2}}L_1^{\frac14+\varepsilon}\|\Pi_\nu u_{2,n}\|_{l_\tau^4L^4}.
$$
Applying the discrete Strichartz estimate \cite[Theorem~2.4]{Jiicos},
$$
\|\Pi_\nu u_n\|_{l_\tau^4L^4}\lesssim\|u_n\|_{X_{1,\tau}^{q,\frac12}},\qquad q>\frac d4-\frac12,
$$
and using $L_0^{1/2}\gtrsim N_1\sim N_0>N_2$, we obtain
\begin{align*}
S(L,N)\lesssim N_2^{\frac{d-1}{2}+\varepsilon}L_1^{\frac14+\varepsilon}L_2^{\frac12}\lesssim N_0^{-\varepsilon}N_1^{-\frac12-\varepsilon}N_2^{\frac{d-1}{2}-\varepsilon}L_0^{\frac14+2\varepsilon}L_1^{\frac14+\varepsilon}L_2^{\frac12}\lesssim N_0^{-\varepsilon}N_1^{-\frac12-\varepsilon}N_2^{\frac{d-1}2-\varepsilon}(L_0L_1L_2)^{\frac38+\varepsilon}.
\end{align*}
since $L_2\leq L_0$ and $\varepsilon\leq1/8$. When $L_{\max}=L_1$, placing $u_{1,n}$ in $l_\tau^2L^2$, $u_{0,n}$ in $l_\tau^4L^2$ and $u_{2,n}$ in $l_\tau^4L^\infty$ gives the same estimate by \eqref{discl4} and the preceding argument with $L_0$ and $L_1$ interchanged.

For $L_{\max}=L_2$ and $L_0\geq L_1$, applying \eqref{v0v1} and H\"older's inequality gives
\begin{align*}
S(L,N)&\lesssim\|\Pi_\nu(\Pi_\nu u_{0,n}\Pi_\nu u_{1,n})\|_{l_\tau^2L^2}\lesssim L_1^{\frac12}\bigg(\frac{L_0}{N_1}+1\bigg)^{\frac12}N_2^{\frac{d-1}2}\lesssim L_1^{\frac12}\bigg(\frac{L_0^{\frac12}L_2^{\frac12}}{N_1}+\frac{L_2^{\frac12}}{N_1}\bigg)^{\frac12}N_2^{\frac{d-1}2}\\
&\lesssim N_0^{-\varepsilon}N_1^{-\frac12-\varepsilon}N_2^{\frac{d-1}{2}-\varepsilon}L_0^{\frac14}L_1^{\frac12}L_2^{\frac14+\frac32\varepsilon}\lesssim N_0^{-\varepsilon}N_1^{-\frac12-\varepsilon}N_2^{\frac{d-1}2-\varepsilon}(L_0L_1L_2)^{\frac38+\varepsilon}.
\end{align*}
since $L_1\leq L_0$ and $0<\varepsilon\leq1/8$. This proves the required bound in this subcase. The same bound for $L_0\leq L_1$ follows by interchanging $L_0$ and $L_1$ in the minimum/maximum form of \eqref{v0v1}. Combining these estimates and restoring the three $l_\tau^2L^2$ norms proves \eqref{slnnonres}.
\end{proof}

\begin{lemma}\label{lemres}
Let $S(L,N)$ be as in \eqref{defSLN}. If $N_0\lesssim N_1\sim N_2$, then
\begin{equation}\label{slnres}
S(L,N)\lesssim N_0^{s_2-\varepsilon}(L_0L_1L_2)^{b-\varepsilon}\|v_{0,n}^{L_0N_0}\|_{l_\tau^2L^2}\|v_{1,n}^{L_1N_1}\|_{l_\tau^2L^2}\|v_{2,n}^{L_2N_2}\|_{l_\tau^2L^2},
\end{equation}
where $\varepsilon>0$ is sufficiently small.
\end{lemma}

\begin{proof}
Set $u_{j,n}=v_{j,n}^{L_jN_j}$ for $j=0,1,2$ and normalize $\|u_{j,n}\|_{l_\tau^2L^2}=1$ by homogeneity. Applying \eqref{disclinf}, H\"older's inequality and Sobolev embedding yields the bound
$$
S(L,N)\lesssim\|\Pi_\nu u_{0,n}\|_{l_\tau^\infty L^\infty}\|u_{1,n}\|_{l_\tau^2 L^2}\|u_{2,n}\|_{l_\tau^2 L^2}\lesssim N_0^{\frac d2+\varepsilon}L_0^{\frac12+\varepsilon}.
$$
Moreover, by Cauchy-Schwarz inequality and \eqref{v1v2}, one obtains
$$
S(L,N)\lesssim\|u_{0,n}\|_{l_\tau^2 L^2}\|\Pi_\nu(\Pi_\nu u_{1,n}\Pi_\nu\overline{u_{2,n}})\|_{l_\tau^2 L^2}\lesssim N_0^{\frac d2-1+\varepsilon}L_1^\frac12L_2^\frac12.
$$
Taking the geometric mean of these bounds with weights $\theta$ and $1-\theta$ gives
$$
S(L,N)\lesssim N_0^{\frac d2-1+\theta+\varepsilon}L_0^{\frac\theta2+\theta\varepsilon}L_1^{\frac12-\frac\theta2}L_2^{\frac12-\frac\theta2}.
$$
Now we choose $\theta$ and $\varepsilon>0$ appropriately depending on $s$ and $b$ so that above estimate yields \eqref{slnres}. Note that $s_0=\frac d2-1$, this requires
$$
\theta+2\varepsilon<s_1,\qquad \frac\theta2+\theta\varepsilon+\varepsilon<b,\qquad \frac\theta2-\varepsilon>\frac12-b.
$$
When $s_1=s_2-s_0\leq\frac14$, it holds that $\frac12-\frac12 s_1<b<\frac12$; we choose
$$
\theta\in(1-2b,s_1),\qquad\varepsilon<\min\{b-\frac12+\frac12\theta,\frac12s_1-\frac12\theta\}
$$
to achieve \eqref{slnres}. When $s_1>\frac14$, it holds that $\frac38<b<\frac12$; we choose
$$
\theta=\frac14,\qquad \varepsilon<\min\{b-\frac38,\frac12s_1-\frac18\}
$$
to attain \eqref{slnres}. This concludes the proof.
\end{proof}

\subsubsection*{Proof of \eqref{multdisc1}--\eqref{multdisc4} for $d\geq2$.}

\begin{proof}
We first prove \eqref{multdisc1} by establishing its equivalent trilinear formulation as in \eqref{multeqv}:
\begin{equation}\label{multeqv1}
S=\Big|\tau\sum_n\int_{\mathbb{T}^d}\Pi_\nu v_{0,n}\Pi_\nu\overline{v_{1,n}}\Pi_\nu v_{2,n}\,dx\Big|\lesssim\|v_{0,n}\|_{X_{2,\tau}^{s_2,b}}\|v_{1,n}\|_{X_{1,\tau}^{-s_2-\frac12,b}}\|v_{2,n}\|_{X_{1,\tau}^{s_2+\frac12,b}}.
\end{equation}
The dyadic decomposition and the triangle inequality give $S\leq S_5+S_6+S_7$, where
$$
S_5=\sum_{\substack{L,N\\N_1\geq8N_2}}S(L,N),\qquad
S_6=\sum_{\substack{L,N\\N_2\geq8N_1}}S(L,N),\qquad
S_7=\sum_{\substack{L,N\\N_1<8N_2,\;N_2<8N_1}}S(L,N).
$$
The last region satisfies $N_1\sim N_2$ and $N_0\lesssim N_1$ whenever the trilinear interaction
is nonzero.

For the non-resonant terms, choose $0<\varepsilon<b-\frac38$ sufficiently small for Lemma~\ref{lemnonres} and fix $r\geq s_0=\frac d2-1$. Writing
$$
W_r(L,N)=\|v_{0,n}^{L_0N_0}\|_{X_{2,\tau}^{r,b}}\|v_{1,n}^{L_1N_1}\|_{X_{1,\tau}^{-r-\frac12,b}}\|v_{2,n}^{L_2N_2}\|_{X_{1,\tau}^{r+\frac12,b}},
$$
we retain the frequency relations when comparing the spatial weights before summing the blocks.

On the support of $S_5$, Lemma~\ref{lemnonres} yields
$$
S_5\lesssim\sum_{L,N}N_2^{s_0}(N_0N_1N_2)^{-\varepsilon}(L_0L_1L_2)^{-(b-\frac38-\varepsilon)}W_0(L,N)\lesssim\sum_{L,N}(N_0N_1N_2)^{-\varepsilon}(L_0L_1L_2)^{-(b-\frac38-\varepsilon)}W_r(L,N).
$$
since $N_0\sim N_1$. This implies by Cauchy-Schwarz that
$$
S_5\lesssim\Big(\sum_{L,N}(N_0N_1N_2)^{-2\varepsilon}(L_0L_1L_2)^{-2(b-\frac38-\varepsilon)}\Big)^\frac12\Big(\sum_{L,N}W_r(L,N)^2\Big)^\frac12,
$$
and the factorization of this square sum together with the Littlewood--Paley characterization of the Bourgain norms yields
\begin{equation}\label{s5}
S_5\lesssim\|v_{0,n}\|_{X_{2,\tau}^{r,b}}\|v_{1,n}\|_{X_{1,\tau}^{-r-\frac12,b}}\|v_{2,n}\|_{X_{1,\tau}^{r+\frac12,b}}.
\end{equation}

For $S_6$, the corresponding version of Lemma~\ref{lemnonres} with the two Schr\"odinger frequencies interchanged gives
$$
S(L,N)\lesssim(N_0N_1N_2)^{-\varepsilon}N_1^{s_0+\frac12}N_2^{-\frac12}(L_0L_1L_2)^{\frac38+\varepsilon}\|v_{0,n}^{L_0N_0}\|_{l_\tau^2L^2}\|v_{1,n}^{L_1N_1}\|_{l_\tau^2L^2}\|v_{2,n}^{L_2N_2}\|_{l_\tau^2L^2}.
$$
Since $N_0\sim N_2$ and $N_1\leq N_2$, the spatial weights satisfy
$$
\frac{N_1^{s_0+\frac12}N_2^{-\frac12}}{N_0^rN_1^{-r-\frac12}N_2^{r+\frac12}}\lesssim\Big(\frac{N_1}{N_0}\Big)^r\Big(\frac{N_1}{N_2}\Big)^{r+1}N_1^{s_0-r}\leq1.
$$
The same summation argument consequently proves
\begin{equation}\label{s6}
S_6\lesssim\|v_{0,n}\|_{X_{2,\tau}^{r,b}}\|v_{1,n}\|_{X_{1,\tau}^{-r-\frac12,b}}\|v_{2,n}\|_{X_{1,\tau}^{r+\frac12,b}}.
\end{equation}
Taking $r=s_2$ gives the non-resonant contributions to \eqref{multeqv1}, while $r=s_0$ provides the corresponding estimates needed for \eqref{multdisc3} and \eqref{multdisc4}.

For the resonant term, set $\widehat{s}=s_2+\frac12$ and choose $\varepsilon>0$ sufficiently small as in Lemma~\ref{lemres}. That lemma and the dyadic norm equivalences give
\begin{equation}\label{s71}
\begin{aligned}
S_7&\lesssim\sum_{\substack{L,N\\N_1\sim N_2}}N_0^{s_2-\varepsilon}(L_0L_1L_2)^{b-\varepsilon}
\prod_{j=0}^2\|v_{j,n}^{L_jN_j}\|_{l_\tau^2L^2}\\
&\lesssim\sum_{\substack{L,N\\N_1\sim N_2}}\bigg(\frac{N_1}{N_2}\bigg)^{\widehat{s}}N_0^{-\varepsilon}(L_0L_1L_2)^{-\varepsilon}\|v_{0,n}^{L_0N_0}\|_{X_{2,\tau}^{s_2,b}}\|v_{1,n}^{L_1N_1}\|_{X_{1,\tau}^{-\widehat{s},b}}\|v_{2,n}^{L_2N_2}\|_{X_{1,\tau}^{\widehat{s},b}}\\
&\lesssim\|v_{0,n}\|_{X_{2,\tau}^{s_2,b}}\sum_{N_1\sim N_2}\bigg(\frac{N_1}{N_2}\bigg)^{\widehat{s}}\Big(\sum_{L_1}\|v_{1,n}^{L_1N_1}\|_{X_{1,\tau}^{-\widehat{s},b}}^2\Big)^\frac12\Big(\sum_{L_2}\|v_{2,n}^{L_2N_2}\|_{X_{1,\tau}^{\widehat{s},b}}^2\Big)^\frac12.
\end{aligned}
\end{equation}
The last inequality in \eqref{s71} follows by Cauchy--Schwarz in $(L_0,N_0)$, $L_1$ and $L_2$ since the squared weights $N_0^{-2\varepsilon}$ and $L_j^{-2\varepsilon}$ have convergent dyadic sums.

Since $N_1\sim N_2$, we may write $N_2=2^pN_1$ with $|p|\leq p_0$ for a fixed integer $p_0$. Employing Cauchy--Schwarz again in $N_1$ to give
\begin{equation}\label{s72}
\begin{aligned}
\sum_{N_1\sim N_2}&\bigg(\frac{N_1}{N_2}\bigg)^{\widehat{s}}\Big(\sum_{L_1}\|v_{1,n}^{L_1N_1}\|_{X_{1,\tau}^{-\widehat{s},b}}^2\Big)^\frac12\Big(\sum_{L_2}\|v_{2,n}^{L_2N_2}\|_{X_{1,\tau}^{\widehat{s},b}}^2\Big)^\frac12\\
&\le\sum_{|p|\leq p_0}2^{-p\widehat{s}}\Big(\sum_{L_1,N_1}\|v_{1,n}^{L_1N_1}\|_{X_{1,\tau}^{-\widehat{s},b}}^2\Big)^\frac12\Big(\sum_{\substack{L_2,N_2\\N_2=2^pN_1}}\|v_{2,n}^{L_2N_2}\|_{X_{1,\tau}^{\widehat{s},b}}^2\Big)^\frac12\lesssim\|v_{1,n}\|_{X_{1,\tau}^{-\widehat{s},b}}\|v_{2,n}\|_{X_{1,\tau}^{\widehat{s},b}}.
\end{aligned}
\end{equation}
Combining \eqref{s5}--\eqref{s6} at $r=s_2$ with \eqref{s71}--\eqref{s72} at $\widehat{s}=s_2+\frac12$ establishes \eqref{multeqv1} and hence
$$
\|\Pi_\nu(\Pi_\nu v_n\Pi_\nu w_n)\|_{X_{1,\tau}^{s_2+\frac12,-b}}\lesssim\|v_n\|_{X_{2,\tau}^{s_2,b}}\|w_n\|_{X_{1,\tau}^{s_2+\frac12,b}}.
$$
Similarly, we obtain
$$
\|\Pi_\nu(\Pi_\nu\overline{v_n}\Pi_\nu w_n)\|_{X_{1,\tau}^{s_2+\frac12,-b}}\lesssim\|v_n\|_{X_{2,\tau}^{s_2,b}}\|w_n\|_{X_{1,\tau}^{s_2+\frac12,b}},
$$
which completes the proof of \eqref{multdisc1}.

For \eqref{multdisc3} and \eqref{multdisc4}, the same duality argument reduces the proof to
\begin{align*}
S_i\lesssim\|v_{0,n}\|_{X_{2,\tau}^{s_2,b}}\|v_{1,n}\|_{X_{1,\tau}^{-s_0-\frac12,b}}\|v_{2,n}\|_{X_{1,\tau}^{s_0+\frac12,b}},\quad
S_i\lesssim\|v_{0,n}\|_{X_{2,\tau}^{s_0,b}}\|v_{1,n}\|_{X_{1,\tau}^{-s_0-\frac12,b}}\|v_{2,n}\|_{X_{1,\tau}^{s_2+\frac12,b}}
\end{align*}
for $i=5,6,7$. Note that the non-resonant contributions $S_5$ and $S_6$ follow from \eqref{s5}--\eqref{s6} with $r=s_0$ and the embeddings
$$
X_{2,\tau}^{s_2,b}\hookrightarrow X_{2,\tau}^{s_0,b},\qquad
X_{1,\tau}^{s_2+\frac12,b}\hookrightarrow X_{1,\tau}^{s_0+\frac12,b}.
$$
Taking $\widehat{s}=s_0+\frac12$ in \eqref{s71}--\eqref{s72} gives
$$
S_7\lesssim\|v_{0,n}\|_{X_{2,\tau}^{s_2,b}}\|v_{1,n}\|_{X_{1,\tau}^{-s_0-\frac12,b}}\|v_{2,n}\|_{X_{1,\tau}^{s_0+\frac12,b}},
$$
which is the resonant estimate required for \eqref{multdisc4}.

To obtain the estimate \eqref{multdisc3} for $S_7$, we use $N_0\lesssim N_2$ and $s_2>s_0$ to transfer the additional spatial weight at the block level:
\begin{equation}\label{substi3}
\begin{aligned}
\|v_{0,n}^{L_0N_0}\|_{X_{2,\tau}^{s_2,b}}&\|v_{2,n}^{L_2N_2}\|_{X_{1,\tau}^{s_0+\frac12,b}}\sim N_0^{s_2-s_0}\|v_{0,n}^{L_0N_0}\|_{X_{2,\tau}^{s_0,b}}\|v_{2,n}^{L_2N_2}\|_{X_{1,\tau}^{s_0+\frac12,b}}\\
&\lesssim N_2^{s_2-s_0}\|v_{0,n}^{L_0N_0}\|_{X_{2,\tau}^{s_0,b}}\|v_{2,n}^{L_2N_2}\|_{X_{1,\tau}^{s_0+\frac12,b}}\sim\|v_{0,n}^{L_0N_0}\|_{X_{2,\tau}^{s_0,b}}\|v_{2,n}^{L_2N_2}\|_{X_{1,\tau}^{s_2+\frac12,b}}.
\end{aligned}
\end{equation}
Applying \eqref{substi3} before the Cauchy--Schwarz step in \eqref{s71} with $\widehat{s}=s_0+\frac12$ leaves the residual weights $N_0^{-\varepsilon}(L_0L_1L_2)^{-\varepsilon}$ unchanged. This proves \eqref{multdisc3} for $S_7$ and finally conclude \eqref{multdisc3} and \eqref{multdisc4} for $d\ge2$.
\end{proof}

\subsection{Proof of \eqref{multdisc2} and \eqref{multdisc5} for $d\geq1$.}

\begin{proof}
As in \eqref{multeqv} and \eqref{multeqv1}, we estimate
$$
S=\Big|\tau\sum_n\int_{\mathbb{T}^d}\Pi_\nu v_{0,n}\Pi_\nu\overline{v_{1,n}}\Pi_\nu v_{2,n}\,dx\Big|.
$$
With $I_\tau=[-\pi/\tau,\pi/\tau)$, set
\begin{align*}
S_8=&\sum_{\langle k_1\rangle\ge2\langle k_2\rangle}\int_{I_\tau^2}\widetilde v_{0,n}(\sigma_1-\sigma_2,k_1-k_2)\widetilde v_{1,n}(\sigma_1,k_1)\widetilde v_{2,n}(\sigma_2,k_2)\,d\sigma_1d\sigma_2,\\
S_9=&\sum_{\langle k_1\rangle<2\langle k_2\rangle}\int_{I_\tau^2}\widetilde v_{0,n}(\sigma_1-\sigma_2,k_1-k_2)\widetilde v_{1,n}(\sigma_1,k_1)\widetilde v_{2,n}(\sigma_2,k_2)\,d\sigma_1d\sigma_2,
\end{align*}
where the temporal Fourier transforms are extended periodically. Fourier inversion and the triangle inequality give $S\lesssim S_8+S_9$, while $k_0=k_1-k_2$ implies
$$
\langle k_0\rangle\lesssim\langle k_1\rangle\quad\text{on the region defining }S_8,\qquad
\langle k_0\rangle\lesssim\langle k_2\rangle\quad\text{on the region defining }S_9.
$$
Consequently, for any $\kappa\geq0$, we may bound the integrand of $S_8$ by inserting $\langle k_0\rangle^{-\kappa} \langle k_1\rangle^\kappa$, and that of $S_9$ by inserting $\langle k_0\rangle^{-\kappa}\langle k_2\rangle^\kappa$, up to constants depending on $\kappa$.

The required trilinear bounds are
$$
S\lesssim\|v_{0,n}\|_{X_{2,\tau}^{-s_2-1,b}}\|v_{1,n}\|_{X_{1,\tau}^{s_2+\frac12,b}}\|v_{2,n}\|_{X_{1,\tau}^{s_2+\frac12,b}},\quad
S\lesssim\|v_{0,n}\|_{X_{2,\tau}^{-s_0-1,b}}\|v_{1,n}\|_{X_{1,\tau}^{s_2+\frac12,b}}\|v_{2,n}\|_{X_{1,\tau}^{s_0+\frac12,b}}.
$$
For the first bound, choose $\kappa=s_2+1$ and apply \eqref{multeqv} for $d=1$ and $\kappa=2s_2+1$ \eqref{multeqv1} for $d\ge2$ to $(\langle\nabla\rangle^{-\kappa}v_{0,n},\langle\nabla\rangle^\kappa v_{1,n},v_{2,n})$ on the region defining $S_8$. This gives
$$
S_8\lesssim\|v_{0,n}\|_{X_{2,\tau}^{-s_2-1,b}}\|v_{1,n}\|_{X_{1,\tau}^{s_2+\frac12,b}}\|v_{2,n}\|_{X_{1,\tau}^{s_2+\frac12,b}},
$$
while the exchanged estimate applied to $(\langle\nabla\rangle^{-\kappa}v_{0,n},v_{1,n},\langle\nabla\rangle^\kappa v_{2,n})$ gives the same bound for $S_9$.

For the mixed regularity estimate, $\kappa=1$ and apply \eqref{multeqv} for $d=1$ and $\kappa=s_0+s_2+1$ \eqref{multeqv1} for $d\ge2$ yields
$$
S_8\lesssim\|v_{0,n}\|_{X_{2,\tau}^{-s_0-1,b}}\|v_{1,n}\|_{X_{1,\tau}^{s_2+\frac12,b}}\|v_{2,n}\|_{X_{1,\tau}^{s_0+\frac12,b}},\quad
S_9\lesssim\|v_{0,n}\|_{X_{2,\tau}^{-s_0-1,b}}\|v_{1,n}\|_{X_{1,\tau}^{s_2+\frac12,b}}\|v_{2,n}\|_{X_{1,\tau}^{s_0+\frac12,b}}.
$$
Combining the two frequency regions proves the mixed trilinear bound. Moreover, repeating the argument for the opposite wave modulation gives the corresponding conjugate-product estimates. This proves \eqref{multdisc2} and \eqref{multdisc5} in all dimensions.
\end{proof}

\end{appendices}

{}

\end{document}